\documentclass[11pt]{book}
 \usepackage[usenames,dvipsnames]{color}
 \usepackage{pgf,tikz}
 \usetikzlibrary{decorations.markings}
\usetikzlibrary{arrows,shapes,positioning,calc,shadows}
\tikzstyle arrowstyle=[scale=1]
\tikzstyle directed=[postaction={decorate,decoration={markings,
    mark=at position .5 with {\arrow[arrowstyle]{stealth}}}}]
\usepackage[kerning=true]{microtype}
\usepackage{fullpage}
\usepackage{amsthm}
\usepackage{bm}
\usepackage{shuffle}
\usepackage[applemac]{inputenc}

\usepackage[T1]{fontenc}

\usepackage{graphicx}
\usepackage[mathcal]{euscript}

\usepackage[french]{todonotes}

\usepackage{wrapfig}
\usepackage{amssymb}
\usepackage{amsmath}

\usepackage[colorlinks=true,
citecolor=cyan,urlcolor=blue,linktocpage=true,pagebackref]{hyperref}

 \def\dessin#1#2{\includegraphics[#1]{#2}}

\newtheorem{exer}{\bleu{\bf Exercice}}[chapter]
  
 \definecolor{midgreen}{rgb}{0.,0.7,0.} 
  \definecolor{darkgreen}{rgb}{0.,0.6,0.} 

\def\rouge{\textcolor{red}}
\def\bleu{\textcolor{blue}}
\def\jaune{\textcolor{yellow}}
\def\vertpale{\textcolor{midgreen}}
\def\vert{\textcolor{darkgreen}}

\definecolor{paleblue}{cmyk}{.4,0,0,0}
\definecolor{grispale}{cmyk}{.05,0,0,0.1}
\def\grispale{\textcolor{grispale}}

 \def\carrebleupale#1#2#3{\node[fill=paleblue,shape=rectangle,  inner sep=+1pt, minimum size=.5cm-5\pgflinewidth] at (#1+.5,#2-.5) {\huge$ #3$}}

\def\dx{0.1}
\def\dy{0}
\def\facteur{.8}
 \def\cellule#1#2{ \path [fill=blue] (\dx+#1,\dy+#2) rectangle ($(\dx+#1,\dy+#2)+\facteur*(1,1)$);}
 
 \def\pcellule#1#2{\carrebleupale{#1}{#2}{}}
\def\marquee#1#2{\carrered{#1}{#2}{}}
\def\marquee#1#2{ \path [fill=red] (\dx+#1,\dy+#2) rectangle ($(\dx+#1,\dy+#2)+\facteur*(1,1)$);}

  \def\bras#1#2{\carrevert{#1}{#2}{}}
   \def\bras#1#2{ \path [fill=midgreen] (\dx+#1,\dy+#2) rectangle ($(\dx+#1,\dy+#2)+\facteur*(1,1)$);}

  \def\jambe#1#2{\carrejaune{#1}{#2}{}}
 \def\jambe#1#2{ \path [fill=yellow] (\dx+#1,\dy+#2) rectangle ($(\dx+#1,\dy+#2)+\facteur*(1,1)$);}

 \def\entree#1#2#3{\node[fill=yellow,shape=rectangle,  inner sep=+1pt, minimum size=.5cm-5\pgflinewidth] at (#1+.5,#2-.5) {$ #3$}}

 \newcommand{\cellrowbleu}[2]{
 \foreach \x in {0,...,#2} {
    \pcellule{\x}{#1} ;}}
    
     \newcommand{\cellrow}[2]{
 \foreach \x in {0,...,#2} {
    \jambe{\x}{#1} ;}}
    
    \newcommand{\crochet}[4]{\coordinate (A) at (#1,#2); 
					  \coordinate (B) at ($(A)+(#3,0)$); 
					  \coordinate (C) at ($(B)-(0,#4)$); 
                                            \draw[red,ultra thick] (A) -- (B) -- (C) ;}

   \newcommand{\crochetbleu}[4]{\coordinate (A) at (#1,#2); 
					  \coordinate (B) at ($(A)+(#3,0)$); 
					  \coordinate (C) at ($(B)-(0,#4)$); 
                                            \draw[blue,ultra thick] (A) -- (B) -- (C) ;}

\def\verythinlines{\linethickness{.25pt}}
\def\verythicklines{\linethickness{2pt}}

\newcommand\audessus[2]{\genfrac{}{}{0pt}{}{#1}{#2}}
\def\lettre#1#2{\Big(\begin{matrix}#1\\[-3pt] #2\end{matrix}\Big)}

\def\MarkedCase#1#2{\put(#1){\mcase{#2}}}
\def\mcase#1{\put(0,0.5){\vcarre}\put(0.25,0.27){#1}}
\def\RSK{\ \longleftrightarrow\ }
\def\case#1{\put(0,0.5){\jcarre}\put(0.25,0.27){#1}}
\def\Case#1#2{\put(#1){\case{#2}}}
\newdimen\carrelength
\def\carre{\bleu{\linethickness{\carrelength}\line(1,0){.85}}}
\def\vcarre{\vertpale{\linethickness{\carrelength}\line(1,0){.85}}}

\def\jcarre{\jaune{\linethickness{\carrelength}\line(1,0){.85}}}

\def\verysmallsquarres{\setlength{\unitlength}{2.5mm}
                                      \setlength{\carrelength}{2.2mm}}
\def\smallsquarres{\setlength{\unitlength}{3mm}
                               \setlength{\carrelength}{2.5mm}}
\def\smallersquarres{\setlength{\unitlength}{4mm}
                           \setlength{\carrelength}{3.5mm}}
                           
\def\ligne(#1,#2)#3{{\multiput(#1,#2)(1,0){#3}{\carre}}}
\def\jligne(#1,#2)#3{{\multiput(#1,#2)(1,0){#3}{\jcarre}}}
\def\colonne(#1,#2)#3{{\multiput(#1,#2)(0,1){#3}{\carre}}}
\def\jcolonne(#1,#2)#3{{\multiput(#1,#2)(0,1){#3}{\jcarre}}}

\def\define#1{{\bf #1}}
\makeatletter
\def\@tocline#1#2#3#4#5#6#7{\relax
  \ifnum #1>\c@tocdepth 
  \else
    \par \addpenalty\@secpenalty\addvspace{#2}%
    \begingroup \hyphenpenalty\@M
    \@ifempty{#4}{%
      \@tempdima\csname r@tocindent\number#1\endcsname\relax
    }{%
      \@tempdima#4\relax
    }%
    \parindent\z@ \leftskip#3\relax \advance\leftskip\@tempdima\relax
    \rightskip\@pnumwidth plus4em \parfillskip-\@pnumwidth
    #5\leavevmode\hskip-\@tempdima
      \ifcase #1
       \or\or \hskip 1em \or \hskip 2em \else \hskip 3em \fi%
      #6\nobreak\relax
    \dotfill\hbox to\@pnumwidth{\@tocpagenum{#7}}\par
    \nobreak
    \endgroup
  \fi}
\newtheorem*{rep@theorem}{\rep@title}
\newcommand{\newreptheorem}[2]{
\newenvironment{rep#1}[1]{
\def\rep@title{#2 \ref{##1}}
\begin{rep@theorem}}
{\end{rep@theorem}}}
\makeatother
\def\auteur#1{{\sc #1}}
\def\titreref#1{{\em #1}}
\def\vol#1{{\bf #1}}
\newtheorem{theorem}{\bleu{Theorem}}

\newtheorem{conjecture}{\bleu{Conjecture}}

\newtheorem{prop}{\bleu{Proposition}}

\numberwithin{equation}{section}
\numberwithin{theorem}{section}
\numberwithin{lemma}{section}
\numberwithin{conjecture}{section}
\numberwithin{rmk}{section}
\numberwithin{prop}{section}
\numberwithin{cor}{section}

\newcommand{\qbinom}[2]{\genfrac{[}{]}{0pt}{}{#1}{#2}_q}

\renewcommand{\S}{\mathbb{S}}

\newcommand{\area}{\operatorname{area}}

\newcommand{\Cat}{\operatorname{Cat}}
\newcommand{\Dyck}{\operatorname{Dyck}}
\newcommand{\Park}{\operatorname{Park}}

\def\define#1{\bleu{\bf{#1}}}

\def\MR#1{\href{http://www.ams.org/mathscinet-getitem?mr=#1}{MR#1}}

\def\S{\mathbb{S}}

\def\R{\mathcal{R}}

\def\N{\mathbb{N}}

\def\A{\mathcal{A}}

\def\C{\mathbb{C}}
\def\Q{{\mathbb{Q}}}
\def\Z{\mathbb{Z}}

\def\pref#1{{\rm (\ref{#1})}}

\def\charac{\raise 2pt\hbox{\large$\chi$}}
\def\Id{\mathrm{Id}}

\title{\begin{picture}(0,0)(18,48)\grispale{\thicklines
$$\includegraphics[width=2\textwidth]{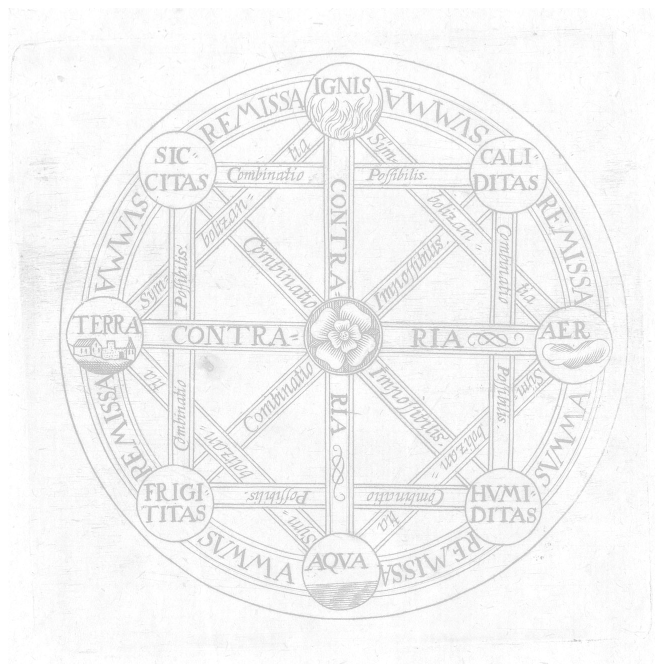}$$}
\end{picture}
\Huge\bf \rouge{Symmetric Functions and Rectangular Catalan Combinatorics}}
\author{\vert{\bf F.~Bergeron}}
  \date{\bleu{\bf \today}}

\begin{document}

\maketitle
 \parskip=0pt
{ \setcounter{tocdepth}{1}\parskip=0pt\footnotesize \tableofcontents}
\parskip=8pt  
\parindent=0pt

\chapter*{\bleu{Introduction}}
\addcontentsline{toc}{chapter}{\bleu{Introduction}}

These notes are intended as a complement for my course in the \href{http://www.risc.jku.at/conferences/aec2018/}{AEC 4th Algorithmic and Enumerative Combinatorics Summer School}, July 30 – August 3, 2018; \href{http://www.risc.jku.at}{RISC}, Hagenberg, Austria. For sure their content is certainly more extensive than what I will realistically be able to cover in the five hours of the course, for which I will assume very little background knowledge. 

I propose many exercises and problems in these notes; but, once again, there are many more of these than what can be solved in the two half-hour periods allocated.  Furthermore, some problems are hard, assuming that one only knows the background in these notes. My aim is for the participants to have something to take home and think about. Beside the exercises that one typically considers in order to become familiar with the notions introduced, I have hence included some problems that may lead to small research projects for the participants that may wish to undertake them. I have also included a wide enough bibliography so that the interested reader is oriented in finding the necessary background, and what is currently happening in the subject.

Much of the notions considered in these notes are interesting to explore using computer algebra tools\footnote{In fact, this is a great way to familiarize oneself with the subject.}, and I give indications of where one may start to learn about these. In particular, there are two special tutorials that were prepared\footnote{See opening paragraphs of chapters 1, 2, and 3.} by graduate students of UQAM: Pauline Hubert and Mélodie Lapointe, on the occasion of a School/Workshop held at CRM in 2017 on \href{http://www.crm.umontreal.ca/2017/Equivariant17/index_e.php}{Equivariant Combinatorics}. Other computer algebra tools are in need of being formatted for general consumption (chapter 4), and this might be one of the outcomes of the school.

Please be lenient with me, since I wrote all of this a bit at the last minute (recycling pieces from other texts of mine). In other words, feel free to underline any typos or mistakes.

\begin{chapter}{\bleu{Combinatorial Background}}
All of this chapter is pretty much classical. For more, see one of the monographs: FB.~\cite{bergeron}, Fulton~\cite{fulton}, or Stanley~\cite{stan_comb}. 
Bruce Sagan also has a nice online text:
\begin{quotation}
\href{http://users.math.msu.edu/users/sagan/Papers/Old/uyt.pdf}{http://users.math.msu.edu/users/sagan/Papers/Old/uyt.pdf}
\end{quotation}
As a fast startup, the web site 
\href{https://en.wikipedia.org/wiki/Young_tableau}{https://en.wikipedia.org/wiki/Young\_tableau}
is  not too bad. Also, see in the Notices of the AMS of February 2007: \href{http://www.ams.org/notices/200702/whatis-yong.pdf}{What is \ldots a Young Tableau?} by Alexander Young (not related). For a fast introduction to partitions in Sage, a nice \href{https://more-sagemath-tutorials.readthedocs.io/en/latest/tutorial-integer-partitions.html}{Sage Partition Tutorial} is available\footnote{https://more-sagemath-tutorials.readthedocs.io/en/latest/tutorial-integer-partitions.html} (It is in part better than the \href{http://doc.sagemath.org/html/en/reference/combinat/sage/combinat/partition.html}{original} for our purpose). It was prepared by  Pauline Hubert and Mélodie Lapointe  (UQAM).

\section{Ferrers diagrams and partitions}\label{Young}

Considering the component-wise partial order on $\N\times \N$, a \define{Ferrers diagram}  $\mu$ (with $n$ cells),  is a cardinal $n$ finite subset of $\N\times \N$ (with pointwise partial order) such that 
     $$(i,j)\leq  (k,l)\quad {\rm and}\quad (k,l)\in\mu\quad  {\rm implies}\quad  (i,j)\in \mu.$$ 
Clearly a Ferrers diagram is characterized by the decreasing integer sequence
$(\mu_1,\mu_2,\ldots,\mu_k)$ with $\mu_i$ denoting the number of cells in the $i^{\rm th}$ row of $\mu$ (reading these row lengths from the bottom to the top). Elements of $\mu$ are called \define{cells}. The \define{conjugate} of a Ferrers diagram $\mu$, denoted by $\mu'$, is the set
   $$\bleu{\mu':=\{(j,i)\ |\ (i,j)\in \mu\}}.$$
Thus, the rows of $\mu'$ are the columns of $\mu$, and vice-versa.

The \define{hook} of a cell $c=(i,j)$ of a Ferrers diagram $\mu$, is the set of cells of $\mu$  that either lie in the same row, to the right of $c$, or lie in the same column, and above $c$. Moreover, $c$ itself belongs to its hook. The corresponding  \define{hook length}  is $ \mathfrak{h}(c)= \mathfrak{h}_{ij}$, is the number of the cells in question. It may be calculated as follows
\begin{equation}
      \bleu{\mathfrak{h}_{ij}:=\mu_{j+1}+\mu_{i+1}'-i-j-1}.
  \end{equation}
The number of cells to the right of $c$ (resp. above) on its row (resp. column) is also called the \define{arm} of $c$ (resp. the \define{leg}) in $\mu$, denoted by $a(c)=a_\mu(c)$ (resp. $\ell(c)=\ell_\mu(c)$). Thus,
   $$\mathfrak{h}(c)=a(c)+\ell(c)+1.$$
 These notions are illustrated in figure~\ref{arm_leg_hook}, with the hook  of the red cell highlighted in yellow and green. The green cells is the arm of $c$, and the yellow ones correspond to  its lef. 
\begin{figure}[ht]
\begin{center}
 \begin{tikzpicture}[thick, scale=0.44]
\cellule07;
\cellule06;\cellule16;
\cellule05;\cellule15;\cellule25;
\cellule04;\cellule14;\cellule24;\jambe34;\cellule44;\cellule54;
\cellule03;\cellule13;\cellule23;\jambe33;\cellule43;\cellule53;\cellule63;\cellule73;
\cellule02;\cellule12;\cellule22;\jambe32;\cellule42;\cellule52;\cellule62;\cellule72;
\cellule01;\cellule11;\cellule21;\marquee31;\bras41;\bras51;\bras61;\bras71;
\cellule00;\cellule10;\cellule20;\cellule30;\cellule40;\cellule50;\cellule60;\cellule70;\cellule80;
\end{tikzpicture}
\end{center}
\vskip-15pt
\caption{A cell and its hook (with arm and leg).}\label{arm_leg_hook}
\end{figure}
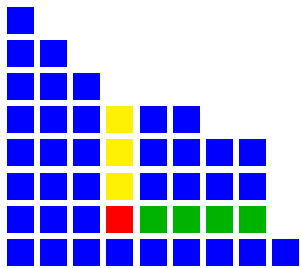
Here, $c=(3,1)$, and its hook length is $ \mathfrak{h}_{21}=8$.
The row lengths of $\mu$ are $(9,8,8,8,6,3,2,1)$, whereas those of $\mu'$ are $(8,7,6,5,5,5,4,4,1)$. The south-west most cell is $(0,0)$.

We are here following the (Cartesian style) right side up  ``French'' convention, rather than the (matrix-style) upside down ``English'' convention, for drawing diagrams. The readers who would like to follow ``Macdonald's advice'' about this should read these notes upside down in a mirror.

To a $n$-cell Ferrers diagram there corresponds a 
\define{partition}\footnote{See \href{https://en.wikipedia.org/wiki/Partition_(number_theory)}{https://en.wikipedia.org/wiki/Partition\_(number\_theory)} for more on partitions.} of the integer $n$. This is simply the decreasing ordered sequence
$(\mu_1,\mu_2,\ldots,\mu_k)$ of row lengths of the diagram. Typically, there is no confusion about also denoting by $\mu$ the  partition associated to a Ferrers diagram $\mu$.  
Each $\mu_i$ is said to be a \define{part} of $\mu$, and  $|\mu|:=\mu_1+\mu_2+\ldots+\mu_k=n$.  We write $\mu\vdash n$ to indicate that
$\mu$ is a partition of $n$. The \define{length}, $\ell(\mu)$, of $\mu$  is simply the number of (non-zero) parts of  $\mu$; hence we have $\ell(\mu')=\mu_1$. 

Partitions are often presented as  \define{words} $\mu=\mu_1\,\mu_2\,\cdots\,\mu_k$,
when all parts are less or equal to $9$. One also considers the \define{empty partition}, denoting it $0$. The \define{partition sets} 
${\rm Part}(n):=\{\mu\ |\ \mu\vdash n\}$, for small $n$, are respectively
 \begin{eqnarray*}
{\rm Part}(0)&=&\{0\}\\ 
{\rm Part}(1)&=&\{1\}\\ 
 {\rm Part}(2)&=&\{2,11\}\\  
 {\rm Part}(3)&=&\{3,21,111\}\\ 
  {\rm Part}(4)&=&\{4,31,22,211,1111\}\\  
  {\rm Part}(5)&=&\{5,41,32,311,221,2111,11111\}\\  
 {\rm Part}(6)&=&\{6,51,42,411,33,321,3111,222,2211,21111,111111\}
\end{eqnarray*}
Another description of partitions consists in writing $\mu=\bm{1}^{d_1}\bm{2}^{d_2}\cdots \bm{j}^{d_j}$, where $d_i$ is the number of parts of size $i$ in $\mu$, and a useful number in the sequel is  
    $$\bleu{z_\mu:=1^{d_1}d_1! 2^{d_2}d_2!\cdots j^{d_j} d_j!}$$
A \define{corner} of $\mu$, is any cell of the form $c=(\mu_i-1,i-1)$ for which if $\mu_i >
\mu_{i+1}$. Corners are exactly the cells that can be removed from the associated Ferrers diagram so that the resulting diagram is also a Ferrers diagram. 
For example, the corners of the partition $\mu={4\,4\,2\,1}$ are
the three dark blue cells in figure~\ref{coins}.
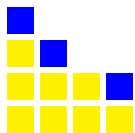
\begin{figure}[ht]
\begin{center}
 \begin{tikzpicture}[thick, scale=0.44]
\cellule03;
\jambe02;\cellule12;
\jambe01;\jambe11;\jambe21;\cellule31;
\jambe00;\jambe10;\jambe20;\jambe30;
\end{tikzpicture}
\end{center}
\vskip-15pt
  \caption{Corners of 4421.}\label{coins}
\end{figure}

If $\nu$ is obtained from $\mu$ by removing one of its corners, we write $\nu\rightarrow  \mu$. The transitive closure of this (covering) relation is a partial order called the \define{Young poset} (see figure~\ref{treillis_de_young}). This simply corresponds to set-inclusion of the corresponding Ferrers diagrams. Since the characteristic property of Ferrers diagram is compatible with union and intersection, the Young poset has the structure of a lattice (as a sub-lattice of the lattice of finite subsets of $\N\times\N$). Hence, for two Ferrers diagrams (or partitions) we have $\nu\rightarrow  \mu$ if and only if $\nu\subset  \mu$, with $\nu\vdash n$ and $\mu\vdash n+1$.

For  $\alpha\subseteq \beta$,  the \define{interval} $[\alpha,\beta]$ is the set  $\{ \mu\ |\ \alpha\subseteq \mu\subseteq \beta\}$. An interesting
special case is the interval $[0,n^k]$ of partitions contained in the \define{rectangular} partition $n^k$. The number of such partitions is easily seen to be equal to the binomial coefficient $\binom{n+k}{k}$. This can be refined to give the classical $q$-analog of the binomial coefficient:
\begin{equation}\label{parta_rect}
  \sum_{\mu\subseteq n^k}q ^{|\mu|} =   \qbinom{n+k}{k}.
 \end{equation}
In general, it is not easy to describe the polynomial 
\begin{equation}\label{q_inclusion}
   P_{\rouge{\mu}}(q):=\sum_{\nu\subseteq \rouge{\mu}}q ^{|\nu|}
\end{equation}
    for a given partition $\rouge{\mu}$.
 In Sage (using $x$ as a variable so that nothing has to be declared) the polynomial in question may be calculated as:
   
  \medskip     
      \centerline{$P_{\rouge{\mu}}(x)$=``\bleu{add($x^k*$Partitions($k$, outer=$\rouge{\mu}$).cardinality() for $k$ in range($n+1$))}''.}
  \medskip

\def\fdots{\reflectbox{$\ddots$}}
\begin{figure}[ht] 
\smallsquarres
$$\begin{matrix}
\begin{picture}(3,4.5)(0,0)
\put(2.2,1.5){\fdots}
\put(-1.3,1.5){$\ddots$}
\put(0,0.5){\carre}
\put(1,0.5){\carre}
\put(2,0.5){\carre}
\end{picture}
&&
\begin{picture}(3,4.5)(-0.5,0)
\put(2.2,2.5){\fdots}
\put(.8,2.5){$\vdots$}
\put(-2,2.5){$\ddots$}
\put(0,0.5){\carre}
\put(1,0.5){\carre}
\put(0,1.5){\carre}
\end{picture}
&&
\begin{picture}(1,4.5)(0,0)
\put(2.2,3.5){\fdots}
\put(-1.3,3.5){$\ddots$}
\put(1,0.5){\carre}
\put(1,1.5){\carre}
\put(1,2.5){\carre}
\end{picture}\\ \\
&
\begin{picture}(2,2)(0,0)
\put(0,1.5){\vector(-1,1){2.5}}
\put(1,1.5){\vector(1,1){2.5}}
\put(-.5,0.5){\carre}
\put(.5,0.5){\carre}
\end{picture}
&&
\begin{picture}(3,2)(0,.5)
\put(.5,2){\vector(-1,1){2.5}}
\put(2,2){\vector(1,1){2.5}}
\put(1,0.5){\carre}
\put(1,1.5){\carre}
\end{picture}\\ 
\\
&&
\begin{picture}(1,2)(0,0)
\put(0,1.5){\vector(-1,1){2.5}}
\put(1,1.5){\vector(1,1){2.5}}
\put(0,0.5){\carre}
\end{picture}\\ 
&&\uparrow\\
&&
\bleu{\bm{0}}
\end{matrix}$$
\caption{Young's lattice.}\label{treillis_de_young}
\end{figure}
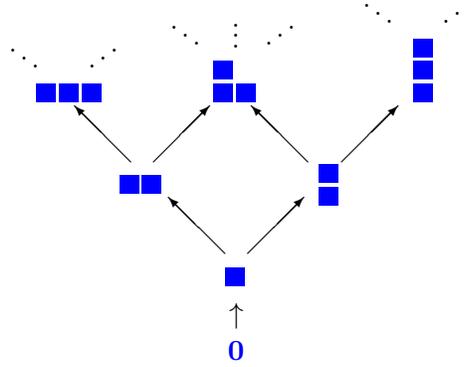

When $\nu$ is contained in $\mu$, we consider the \define{skew partition}, denoted by $\mu/\nu$, having as \define{skew Ferrers diagram} the set difference $\mu\setminus \nu$. 

\subsection*{Dominance order on partitions}
For two partitions $\mu$ and $\lambda$ of $n$, we say that $\mu$ \define{dominates} $\lambda$, and write $\lambda\preceq \mu$, if and only if for all $k$
$$\bleu{\lambda_1 + \lambda_2 + \cdots +
    \lambda_k\leq \mu_1 + \mu_2 + \cdots +\mu_k}.$$
If needed one adds parts $\mu_i=0$ or
$\lambda_i=0$ so that inequalities make sense. The following figure gives the dominance order on
partitions of $n=6$, with an arrow $\mu\rightarrow\lambda$ indicating
that $\mu$ is covered by $\lambda$.

\verysmallsquarres
$$\begin{array}{lclclclclclclclclclclclclclclclclclc}
&&&&&
\begin{picture}(3,2)(-2.3,-0.5)
   \ligne(0,3){1}
   \ligne(0,2){1}
   \ligne(0,1){1}
   \ligne(0,0){3}
\end{picture}
                             &&&
\begin{picture}(3,2)(-2.3,-0.5)
   \ligne(0,1){3}
   \ligne(0,0){3}
\end{picture}
                              \\
 &&&&&\nearrow&\searrow&&\nearrow&\searrow
                              \\
\begin{picture}(1,2.5)(0,2)
   \colonne(0,0){6}
\end{picture}
                              &\rightarrow&
\begin{picture}(2,2.5)(0,1.5)
   \colonne(0,0){5}    \colonne(1,0){1}
\end{picture}
                              &\rightarrow&
\begin{picture}(2,2.5)(0,1)
   \colonne(0,0){4}    \colonne(1,0){2}
\end{picture}
                              &&&
\begin{picture}(3,2.5)(0,0.5)
   \ligne(0,2){1}
   \ligne(0,1){2}
   \ligne(0,0){3}
\end{picture}
                              &&&
\begin{picture}(4,2)(0,0)
   \ligne(0,1){2}
   \ligne(0,0){4}
\end{picture}
                              &\rightarrow&
\begin{picture}(5,2)(0,0)
   \ligne(0,1){1}
   \ligne(0,0){5}
\end{picture}
                              &\rightarrow&
\begin{picture}(6,1)(0,-0.3)
   \ligne(0,0){6}
\end{picture}
                              \\ \\
 &&&&&\searrow&\nearrow&&\searrow&\nearrow\\
&&&&&
\begin{picture}(2,2.2)(-2,0.5)
   \colonne(0,0){3}    \colonne(1,0){3}
\end{picture}
                              &&&
\begin{picture}(4,2.2)(-2,0.5)
   \ligne(0,2){1}
   \ligne(0,1){1}
   \ligne(0,0){4}
\end{picture}
\end{array}$$

This example underlines that the dominance order is not
a total order (although it is for all $n\leq 5$).

\section{Young tableaux, hook length formula, and Kostka numbers}
Let $\bm{d}$ be a finite subset of $\N\times \N$. A \define{tableau} of shape $\bm{d}$, with value in a set $A$ (typically $A=\{1,2,\ldots ,n\}$), is simply a function $\tau:\bm{d}\rightarrow A$. This function is ``displayed'' by filling each cell $c$ of $\bm{d}$ by its value $\tau(c)$, as illustrated in figure~\ref{a_filling}. We mostly consider the case when $\bm{d}$ is a Ferrers diagram $\mu$ (or a skew diagram). 

\setlength{\unitlength}{5mm}
\setlength{\carrelength}{4.3mm}
  \begin{figure}[ht]
  \begin{center}
\begin{picture}(4,3)(0,1)
   \Case{0,3}{4} 
   \Case{0,2}{2} \Case{1,2}{3}
                                                \Case{2,1}{3}\Case{3,1}{8}
   \Case{0,0}{1}                       \Case{2,0}{1}\Case{3,0}{5}
\put(-0.05,-0.05){\vector(1,0){4.5}}
\put(-0.05,-0.05){\vector(0,1){4.5}}
\end{picture}
\end{center}
\caption{A  semi-standard tableau.}
\label{a_filling}
\end{figure}
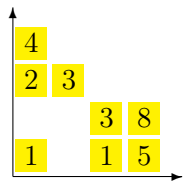
A tableau $\tau$ is said to be \define{semi-standard}  if its entries are non-decreasing along rows (from left to right), and strictly increasing
along columns (from bottom to top) of $\bm{d}$:
    $$i<k\quad {\rm implies}\quad \tau(i,j)\leq \tau(k,j),\qquad {\rm and}\qquad  j<\ell\quad {\rm implies}\quad\tau(i,j)<\tau(i,\ell);$$
 for $(i,j)$, $(k,j)$, and $(i,\ell)$ lying in $\bm{d}$.
 
 A $n$-cell tableau $\tau$ is \define{standard} if it is  a bijective semi-standard tableau with values in $\{1,2,\ldots,n\}$.
Thus, for $\tau$ to be standard we need
$\tau(i,j)< \tau(k,\ell)$ whenever $(i,j)<(k,\ell)$ coordinate-wise. figure~\ref{un_tableau_standard} gives an example of a standard
tableau of shape $431$.
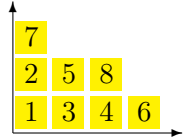
\begin{figure}[ht]
\begin{center}
\begin{picture}(4,3)(0,0)
   \Case{0,2}{7}   
   \Case{0,1}{2}  \Case{1,1}{5}   \Case{2,1}{8}  
   \Case{0,0}{1}  \Case{1,0}{3}   \Case{2,0}{4}  \Case{3,0}{6}
\put(-0.05,-0.05){\vector(1,0){4.5}}
\put(-0.05,-0.05){\vector(0,1){3.5}}
\end{picture}
\end{center}
\caption{A standard tableau.}
\label{un_tableau_standard}
\end{figure}
 
 The \define{reading word}, $\rho(\tau)$, of a tableau $\tau$ is obtained as the row by row reading of the entries of $\tau$. Each row is read from left to right, and rows are read from top to bottom. Thus, the reading word of the tableau in figure~\ref{un_tableau_standard} is: $72581346$. One may check that a partition shaped semi-standard tableau is entirely characterized by its reading word.

A standard tableau of shape $\mu$ corresponds to a
maximal chain 
  $$0=\mu^{(0)}\rightarrow \mu^{(1)}\rightarrow\,\ldots \rightarrow  \mu^{(n)}=\mu,$$ 
where the partition $\mu^{(i)}$ is made out of the cells taking values $\leq i$. For instance, the standard tableau in figure~\ref{un_tableau_standard}
correspond (bijectively) to the maximal chain of figure~\ref{max_chain_tableau}.
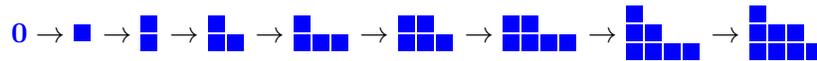
\begin{figure}[ht]
\verysmallsquarres
\begin{center}
$\bleu{\bm{0}}\rightarrow$
\begin{picture}(1,1)(0,-0.5)
   \ligne(0,0){1}
\end{picture}
                         $\rightarrow$
\begin{picture}(1,2)(0,0)
   \ligne(0,1){1}
   \ligne(0,0){1}
\end{picture}
                          $\rightarrow$
\begin{picture}(2,2)(0,0)
   \ligne(0,1){1}
   \ligne(0,0){2}
\end{picture}
                          $\rightarrow$
\begin{picture}(3,2)(0,0)
   \ligne(0,1){1}
   \ligne(0,0){3}
\end{picture}
                          $\rightarrow$
\begin{picture}(3,2)(0,0)
   \ligne(0,1){2}
   \ligne(0,0){3}
\end{picture}
                           $\rightarrow$
\begin{picture}(4,2)(0,0)
   \ligne(0,1){2}
   \ligne(0,0){4}
\end{picture}
                        $\rightarrow$
\begin{picture}(4,2)(0,0.5)
   \ligne(0,2){1}
   \ligne(0,1){2}
   \ligne(0,0){4}
\end{picture}
                           $\rightarrow$
\begin{picture}(4,2)(0,0.5)
   \ligne(0,2){1}
   \ligne(0,1){3}
   \ligne(0,0){4}
\end{picture}
\end{center}
\caption{A maximal chain in Young's lattice}
\label{max_chain_tableau}
\end{figure}  

\subsection{Hook length formula}\label{hook_length}

The number $f^\mu$ of standard tableaux of shape $\mu$ ($\mu\vdash n$)
is given by the Frame-Robinson-Thrall (see \cite{frame}) hook length formula.
\begin{equation}\label{hook_length_formula}
    {f^\mu=\frac{n!}{\prod_{c\in \mu}{\mathfrak{h}(c)}},}
 \end{equation}
By a direct application of this formula, there are exactly $16$ standard tableaux of shape $321$. These are:
$$\begin{array}{ccccccccccc} 
\begin{picture}(3,3)(0,0)
   \Case{0,2}{3}  
   \Case{0,1}{2} \Case{1,1}{5} 
   \Case{0,0}{1} \Case{1,0}{4} \Case{2,0}{6} 
\end{picture}
 &
\begin{picture}(3,2.5)(0,0)
    \Case{0,2}{3}  
   \Case{0,1}{2} \Case{1,1}{6} 
   \Case{0,0}{1} \Case{1,0}{4} \Case{2,0}{5} 
\end{picture}
&
\begin{picture}(3,2.5)(0,0)
   \Case{0,2}{4}  
   \Case{0,1}{2} \Case{1,1}{5} 
   \Case{0,0}{1} \Case{1,0}{3} \Case{2,0}{6} 
\end{picture}
&
\begin{picture}(3,2.5)(0,0)
    \Case{0,2}{4}  
   \Case{0,1}{2} \Case{1,1}{6} 
   \Case{0,0}{1} \Case{1,0}{3} \Case{2,0}{5} 
\end{picture}
& 
\begin{picture}(3,3.5)(0,0)
    \Case{0,2}{4}  
   \Case{0,1}{3} \Case{1,1}{5} 
   \Case{0,0}{1} \Case{1,0}{2} \Case{2,0}{6} 
\end{picture}
&
\begin{picture}(3,3)(0,0)
   \Case{0,2}{4}  
   \Case{0,1}{3} \Case{1,1}{6} 
   \Case{0,0}{1} \Case{1,0}{2} \Case{2,0}{5} 
\end{picture}
 &
\begin{picture}(3,3)(0,0)
   \Case{0,2}{5}  
   \Case{0,1}{2} \Case{1,1}{4} 
   \Case{0,0}{1} \Case{1,0}{3} \Case{2,0}{6} 
\end{picture}
&
\begin{picture}(3,3)(0,0)
   \Case{0,2}{5}  
   \Case{0,1}{2} \Case{1,1}{6} 
   \Case{0,0}{1} \Case{1,0}{3} \Case{2,0}{4} 
\end{picture}
\\ 
\begin{picture}(3,3.5)(0,0)
    \Case{0,2}{5}  
   \Case{0,1}{3} \Case{1,1}{4} 
   \Case{0,0}{1} \Case{1,0}{2} \Case{2,0}{6} 
\end{picture}
 &
\begin{picture}(3,3)(0,0)
    \Case{0,2}{5}  
   \Case{0,1}{3} \Case{1,1}{6} 
   \Case{0,0}{1} \Case{1,0}{2} \Case{2,0}{4} 
\end{picture}
 &
\begin{picture}(3,3)(0,0)
    \Case{0,2}{5}  
   \Case{0,1}{4} \Case{1,1}{6} 
   \Case{0,0}{1} \Case{1,0}{2} \Case{2,0}{3} 
\end{picture}
&
\begin{picture}(3,3)(0,0)
    \Case{0,2}{6}  
   \Case{0,1}{2} \Case{1,1}{4} 
   \Case{0,0}{1} \Case{1,0}{3} \Case{2,0}{5} 
\end{picture}
& 
\begin{picture}(3,3.5)(0,0)
    \Case{0,2}{6}  
   \Case{0,1}{2} \Case{1,1}{5} 
   \Case{0,0}{1} \Case{1,0}{3} \Case{2,0}{4} 
\end{picture}
&
\begin{picture}(3,3)(0,0)
    \Case{0,2}{6}  
   \Case{0,1}{3} \Case{1,1}{4} 
   \Case{0,0}{1} \Case{1,0}{2} \Case{2,0}{5} 
\end{picture}
 &
\begin{picture}(3,3)(0,0)
    \Case{0,2}{6}  
   \Case{0,1}{3} \Case{1,1}{5} 
   \Case{0,0}{1} \Case{1,0}{2} \Case{2,0}{4} 
\end{picture}
 &
\begin{picture}(3,3)(0,0)
    \Case{0,2}{6}  
   \Case{0,1}{4} \Case{1,1}{5} 
   \Case{0,0}{1} \Case{1,0}{2} \Case{2,0}{3} 
\end{picture}
\end{array}$$
It is easy to check that the sum of the hook lengths of a partition $\mu$ is given by the formula 
      $$\sum_{c\in \mu}\mathfrak{h}(c)=n(\mu)+n(\mu')+|\mu|,\qquad {\rm where}\qquad n(\bm{d}):=\sum_{(i,j)\in \bm{d}} j.$$
The classical Robinson-Schensted-Knuth correspondence (see next section)  shows that
   $$\sum_{\mu\vdash n} (f^\mu)^2=n!,$$
and its properties imply that $\sum_{\mu\vdash n} f^\mu$ is the number of involutive permutations. 

\subsection{Kostka numbers}\label{kostka}
The \define{content} $\gamma(\tau)$ of a tableau
$\tau$ is the sequence $\gamma(\tau)=(m_1,m_2,m_3,\ldots)$ 
of \define{multiplicities} of each entry  $i$ in the tableau $\tau$.  For example, the content of the semi-standard tableau
\begin{center}
 \begin{picture}(5,2.5)(0,0)
   \Case{0,2}{4} \Case{1,2}{4} 
   \Case{0,1}{2} \Case{1,1}{2} \Case{2,1}{4}  \Case{3,1}{4}
   \Case{0,0}{1} \Case{1,0}{1} \Case{2,0}{1}  \Case{3,0}{1} \Case{4,0}{2}
\end{picture}
\end{center}
 is $\gamma(\tau)=(4,3,0,4,0,0,\ldots)$. If $\lambda$ and $\mu$ are two partitions of  $n$, we define the \define{Kostka number}
$K_{\lambda,\mu}$ to be the number of semi-standard tableaux of shape
$\lambda$ and content $\mu$. For instance, the  4 possible semi-standard tableaux having content $2211$ and shape  $321$ are 
  $$\begin{matrix}
   \begin{picture}(5,3)(0,0)
   \Case{0,2}{4}  
   \Case{0,1}{2} \Case{1,1}{3}  
   \Case{0,0}{1} \Case{1,0}{1} \Case{2,0}{2}   
\end{picture} &
   \begin{picture}(5,3)(0,0)
   \Case{0,2}{3}  
   \Case{0,1}{2} \Case{1,1}{4}  
   \Case{0,0}{1} \Case{1,0}{1} \Case{2,0}{2}   
\end{picture} &
   \begin{picture}(5,3)(0,0)
   \Case{0,2}{4}  
   \Case{0,1}{2} \Case{1,1}{2}  
   \Case{0,0}{1} \Case{1,0}{1} \Case{2,0}{3}   
\end{picture} &
   \begin{picture}(5,3)(0,0)
   \Case{0,2}{3}  
   \Case{0,1}{2} \Case{1,1}{2}  
   \Case{0,0}{1} \Case{1,0}{1} \Case{2,0}{4}   
\end{picture} 
\end{matrix}
   $$
The content of a semi-standard tableau is
$\bm{1}^n$ if and only if the tableau is standard. Hence $K_{\lambda,1^n}=f^\lambda$. We may show that 
$K_{\lambda,\lambda}=1$, and that $K_{\lambda,\mu}\not=0$ forces $\mu\preceq \lambda$, in dominance order (you are asked to show this in one of the exercises).

The above observations give a somewhat more natural characterization of the dominance order. Indeed, we have $ \lambda\succeq \mu$ if and only if $K_{\lambda,\mu}\not=0$. Moreover we see that the matrix $(K_{\lambda,\mu})_{\lambda,\mu\vdash n}$
is upper triangular when partitions are sorted in any decreasing order which is a linear 
 extension (why?) of the dominance order. This immediately implies that the matrix $(K_{\lambda,\mu})_{\lambda,\mu\vdash n}$ 
is easily invertible. For $n=4$, the Kostka matrix is
\begin{equation}\label{kostka4}
  \begin{array}{llll} & \begin{array}{lllll} \quad 4 &\  31 &\  22 &\  211 & 1111 \end{array}\\[4pt]
      &  \left(\begin{array}{llllll}
                 1 \quad &  1 \quad  &  1 \quad  &  1 \quad & 1 \\
                 0 & 1 & 1 & 2 & 3 \\
                 0 & 0 & 1 & 1 & 2 \\ 
                 0 & 0 & 0 & 1 & 3\\
                 0 & 0 & 0 & 0 & 1\\ 
            \end{array}\right) \end{array}
    \end{equation}
    
\section{Robinson-Schensted-Knuth}
Let $\A$ and $\mathcal{B}$ be two ordered sets, and consider \define{lexicographic} two row matrices:
\begin{equation}
    w=\begin{pmatrix}
              b_1& b_2&\ldots &b_n\cr
              a_1& a_2& \ldots & a_n
          \end{pmatrix}\qquad \hbox{such that}\ {\rm  either}\ (b_i<b_{i+1})\ {\rm or}\ (b_i=b_{i+1}\ {\rm and}\ a_i\leq a_{i+1}),
 \end{equation}   
 with the $a_i$ in  $\A$, and the $b_i$ in $\mathcal{B}$.
We denote by $\varepsilon$ the \define{empty} case ({\it i.e.}: $n=0$).

The RSK correspondence (see Knuth~\cite{knuth} for more on this) $w\RSK (P,Q)$, 
associates bijectively a pair $(P(w),Q(w))=(P,Q)$ of same shape semi-standard tableaux to each such lexicographic matrix. The entries of $P(w)$ are the $a_i$'s, and those of $Q(w)$ are the $b_i$'s. 
The correspondence is recursively defined as follows.
\begin{enumerate}
\item  First, we start with setting $\varepsilon\RSK (0,0)$.
\item  For a nonempty 
                 $$w=\begin{pmatrix}
              b_1& b_2&\ldots &b_{k-1} & \bleu{b}\cr
              a_1& a_2& \ldots & a_{k-1} &\rouge{a}
          \end{pmatrix}\qquad {\rm set}\qquad v:=\begin{pmatrix}
              b_1& b_2&\ldots &b_{k-1} \cr
              a_1& a_2& \ldots & a_{k-1} 
          \end{pmatrix}$$
       and recursively sets $v\RSK (P',Q')$.
  \item Then construct 
            $$P:=(P'\leftarrow \rouge{a})$$ 
by \define{tableau insertion} (see subsection below)  of $a$ into $P'$. The tableau $Q$ is then obtained by adding $\bleu{b}$ to $Q'$
in the position corresponding to the unique cell by which the shape of $P$ differs from that of $P'$.
\end{enumerate}
This may be displayed as follows
$$
\left(\begin{picture}(6,2)(0,1.3)
      \Case{0,2}{3} \Case{1,2}{5} \Case{2,2}{5}
      \Case{0,1}{2} \Case{1,1}{3} \Case{2,1}{3}
      \Case{0,0}{1} \Case{1,0}{1} \Case{2,0}{1} \Case{3,0}{2} \Case{4,0}{3}  \Case{5,0}{3}
   \end{picture}\,,\,
\begin{picture}(6,2)(0,1.3)
     \Case{0,2}{3} \Case{1,2}{4} \Case{2,2}{4}
     \Case{0,1}{2} \Case{1,1}{2} \Case{2,1}{2}
     \Case{0,0}{1} \Case{1,0}{1} \Case{2,0}{1} \Case{3,0}{1} \Case{4,0}{2}  \Case{5,0}{2}
\end{picture}\right)
\leftarrow \lettre{6}{1}\ = \ 
\left(\begin{picture}(6,2.5)(0,2)
      \MarkedCase{0,3}{5}
      \Case{0,2}{3} \MarkedCase{1,2}{3} \Case{2,2}{5}
      \Case{0,1}{2} \MarkedCase{1,1}{2} \Case{2,1}{3}
      \Case{0,0}{1} \Case{1,0}{1} \Case{2,0}{1} \MarkedCase{3,0}{1} \Case{4,0}{3}  \Case{5,0}{3}
   \end{picture}\,,\,
\begin{picture}(6,2.5)(0,2)
     \MarkedCase{0,3}{6}
     \Case{0,2}{3} \Case{1,2}{4} \Case{2,2}{4}
     \Case{0,1}{2} \Case{1,1}{2} \Case{2,1}{2}
     \Case{0,0}{1} \Case{1,0}{1} \Case{2,0}{1} \Case{3,0}{1} \Case{4,0}{2}  \Case{5,0}{2}
\end{picture}\right)
$$
In the ``classical'' RSK correspondence, $w$ takes the form
   $$\begin{pmatrix}1&2&\cdots & n\cr
            a_1 & a_2& \cdots & a_n
        \end{pmatrix}$$
 and it is then customary to identify $w$ with ``word'' $a_1a_2\cdots a_n$.
Since  the $b_i$ are all distinct, it follows that the tableau $Q$  is actually a standard tableau. When $a_1a_2\ldots a_n$ is a permutation of $\{1,2,\ldots,n\}$,  the tableau $P$ is also standard. This establishes a bijection between permutations in $\S_n$ and pairs of  standard tableaux of same shape. As a by-product, we see that
  \begin{equation}\label{dim_decomp_reg_Sn}
         n!=\sum_{\mu\vdash n} (f^\mu)^2,
   \end{equation}
since $(f^\mu)^2$ counts the number of pairs of standard tableaux of shape $\mu$.

We consider the \define{inverse\/} $w^{-1}$ of $w$:
    $$w^{-1}=\mathrm{lex}\begin{pmatrix}a_1& a_2& \ldots & a_k\cr
              b_1& b_2&\ldots &b_k \end{pmatrix},$$
with ``$\mathrm{lex}$'' standing for the increasing lexicographic reordering. 
Thus with
  $$ w=\begin{pmatrix} 1 & 1 & 1 & 2 & 2 & 3 & 3 & 4 & 4 & 4 \\ 3 & 4 & 5 & 1
 & 4 & 1 & 4 & 1 & 2 & 2
    \end{pmatrix} $$ 
we have
 $$w^{-1}= \begin{pmatrix}  1 & 1 & 1 & 2 & 2 & 3 & 4 & 4 & 4 & 5 \\ 2 & 3 & 4 & 4
 & 4 & 1 & 1 & 2 & 3 & 1\end{pmatrix}  $$
 One can show  that 
\begin{equation}\label{propinvRSK}
 w\RSK (P,Q) \qquad \mathrm{iff}\qquad w^{-1}\RSK (Q,P)
\end{equation}

 \subsubsection*{\underline{Tableau insertion}}
One constructs as follows the tableau $\tau':=(\tau\leftarrow \rouge{a})$ obtained by \define{insertion} of $\rouge{a}\in \A$ in a semi-standard $\tau:\mu\rightarrow \A$.
\begin{itemize}
\item[$\bullet$] First, if $\tau=0$ is the empty tableau, then  $\tau':=(0\leftarrow \rouge{a})$ is simply the tableau of shape $1$ with value $\rouge{a}$ in its single cell. 
\item[$\bullet$] Otherwise, one inserts $\rouge{a}$ in the first row of $\tau$, and two cases are to be considered depending on how $\rouge{a}$ compares to the entries occurring in this first row of $\tau$.
\begin{enumerate}
\item[(a)] If $\rouge{a}$ is larger or equal to the largest entry of the first row, then $\rouge{a}$ is simply appended  to the end of the first row of $\tau$; 
\item[(b)] Otherwise, $\rouge{a}$ replaces the leftmost entry $\vert{x}$ of the first row which is larger than $\rouge{a}$. One says that $\rouge{a}$ \define{bumps} $\vert{x}$. 
 In this case, $\vert{x}$ is recursively inserted by the same process in the tableau consisting of the rows of $\tau$ from the second to the top. This results in a \define{bump path}, consisting in the successive cells altered by the insertion process, which ends in a new (corner) cell for the underlying shape of $\tau$.
\end{enumerate}
\end{itemize}
In the example above, the bump path corresponding to the insertion of $1$  is highlighted in green in the $P$-tableau, as well as the new corner in the $Q$-tableau.

\section{Charge  and cocharge}\label{charge}
We say that a  standard tableau $\tau$ has a \define{reading descent} at a cell $(i,j)$ if   
$\tau(i,j)+1$ lies in a cell that is to the left of $(i,j)$. For example, the reading descents of the tableau  on the  left-hand side of figure~\ref{fig_destand} 
\begin{figure}[ht]
\begin{center}
\begin{picture}(5,3)(0,0)
   \Case{0,2}{5} \Case{1,2}{9}
   \Case{0,1}{3}  \MarkedCase{1,1}{4}   \MarkedCase{2,1}{8}  
   \Case{0,0}{1}  \MarkedCase{1,0}{2}  \Case{2,0}{6}  \MarkedCase{3,0}{7}
   \put(6,1){$\longmapsto$}
\end{picture}
\qquad \qquad
\begin{picture}(5,3)(0,0)
   \Case{0,2}{2} \Case{1,2}{4}
   \Case{0,1}{1}  \Case{1,1}{1}  \Case{2,1}{3}  
   \Case{0,0}{0}  \Case{1,0}{0}  \Case{2,0}{2}  \Case{3,0}{2} \
\end{picture}
\end{center}
\caption{Minimization of a tableau.\label{fig_destand}}
\end{figure}
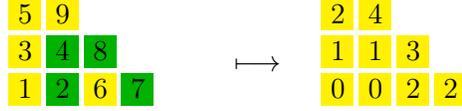
correspond to the cells that contain the values $2$, $4$, $7$ and $8$. 

The   \define{minimization} $\mathfrak{m}(\tau)$, of a $\mu$-shaped semi-standard tableau $\tau$, is the ``minimal'' standard tableau of shape $\mu$ which mimics the reading descent pattern of $\tau$. It is obtained by successively reading off the entries of $\tau$ from $1$ to $n$, and filling the corresponding cells of $\mu$ with entries that stay constant as long as one goes eastward, but rise up by $1$ otherwise. This process may readily be extended to semi-standard tableaux, reading equal entries from left to right.

The \define{co-charge} $\mathrm{coch}(\tau)$ of a standard tableau $\tau$ is simply the sum of the entries of $\mathfrak{m}(\tau)$. The \define{charge} $\mathrm{ch}(\tau)$ of a semi-standard tableau $\tau$ is simply $n(\lambda)-\mathrm{coch}(\tau)$, so that charge is not truly worth an independent description.   The maximal value for the co-charged statistic, for a shape $\lambda$ tableau, is $n(\lambda)$. This maximal value is realized by the standard tableau $\tau$, whose entries are $\tau(i,j)=j-1+\sum_{k<i} \lambda'_i$. Many formulas involving the enumeration of tableaux admit a charge/co-charge refinement. For instance, we have (see \cite{macdonald})
 \begin{equation}\label{cocharge_hook}
     f^\mu(q):=\sum_{\lambda(\tau)=\mu} q^{\mathrm{coch}(\tau)} =
        q^{n(\mu)} \frac{[n]_q!}{\prod_{c\in\mu} (1-q^{\mathfrak{h}(c)})},
  \end{equation}
  where the sum is over the set of standard tableaux of shape $\mu$.  One gets back the usual hook length formula by taking the limit as $q$ tends to $1$. 
  We will encounter later the \define{Kostka-Foulkes polynomials}:
     $$K_{\lambda\mu}(q)=\sum_\tau q^{\mathrm{ch}(\tau)},$$
 with sum over the set of semi-standard tableaux of shape $\lambda$ and content $\mu$. The $q$-Kostka matrix for $n=4$ is:
$$\begin{pmatrix} 
1 & q & q^{2} & q^{3} & q^{6} \\
0 & 1 & q & q^{2} + q & q^{5} + q^{4} + q^{3} \\
0 & 0 & 1 & q & q^{4} + q^{2} \\
0 & 0 & 0 & 1 & q^{3} + q^{2} + q \\
0 & 0 & 0 & 0 & 1
\end{pmatrix}$$

\section{Exercises and problems}
\begin{exer}\rm
Writing $[n]_q$ for the \define{$q$-integer} $\bleu{1+q+\ldots + q^{n-1}}$, one considers the \define{$q$-factorial} 
    $$\bleu{[n]_q!:=[1]_q[2]_q \cdots [n]_q},$$
with $[0]_q!$ set to be equal to $1$. One then defines the \define{$q$-binomial coefficient}:
    $$\bleu{ \textstyle{\qbinom{n}{k}}:=\frac{[n]_q!}{[k]_q!\, [n-k]_q!}}.$$
\begin{enumerate}\itemsep=4pt
\item[(a)] Show that we have the $q$-Pascal-triangle recursion
    $$\textstyle{ \qbinom{n+1}{k}}=q^k\qbinom{n}{k}+\qbinom{n}{k-1},$$
  with appropriate initial conditions. Conclude that $\textstyle{\qbinom{n}{k}}$ is a polynomial in $q$ with positive integer coefficients.
 \item[(b)] Prove the $q$-analog of Newton's theorem:
 \begin{equation}\label{q_newton}
  \textstyle{  \prod_{k=0}^{n-1} (1+q^k\,z)=
       \sum_{k=0}^n q^{k(k-1)/2} \qbinom{n}{k} z^k}.
   \end{equation}
   \end{enumerate}
 \end{exer}
\begin{exer}\rm
  Show that formula~\pref{parta_rect} holds, by proving that the left-hand side satisfies the same recursion as the right-hand side. 
 \end{exer}
\begin{exer}\rm Show that $P_{{\mu}}(q)=P_{{\mu'}}(q)$, see definition~\pref{q_inclusion}, for all partition $\mu$.
 \end{exer}
\begin{exer}\rm Let 
\begin{equation}\label{defn_staircase}
\bleu{\delta_n:=(n-1,n-2,\cdots, 2,1)}
\end{equation}
 be the \define{$n$-staircase} partition (see figure below), and consider the $q$-enumeration of partition sitting inside $\delta_n$:
    $$\mathcal{C}_n(q):=\sum_{\mu\subseteq \delta_n} q^{|\mu|}.$$
Using the ``first return'' to diagonal decomposition illustrated in the following figure: \bigskip
\begin{center}
\thicklines
\verysmallsquarres 
\begin{picture}(8,8)(1,0)
\put(0,-.5){$\underbrace{\hskip28pt}_j$}
\put(-5,1.5){$n-j\left\{\rule{0cm}{20pt}\right.$}
\put(0,4){
   \jligne(0,4){1}
   \jligne(0,3){2}
  \ligne(0,2){1}\jligne(1,2){2}
  \ligne(0,1){2}\jligne(2,1){2}
\put(-0.1,0.5){
\put(0,4){\rouge{\circle*{.2}}}
\put(0,4){\rouge{\line(0,-1){2}}}
\put(0,2){\rouge{\line(1,0){1}}}
\put(1,2){\rouge{\line(0,-1){1}}}
\put(1,1){\rouge{\line(1,0){1}}}
\put(2,1){\rouge{\line(0,-1){1}}}
\put(2,0){\rouge{\line(1,0){2}}}   }}
\put(4,0){
  \jligne(0,4){1}
   \jligne(0,3){2}
  \ligne(0,2){1}\jligne(1,2){2}
  \ligne(0,1){3}\jligne(3,1){1}
  \ligne(0,0){3}\jligne(3,0){2}
\put(-0.1,-0.5){
\put(0,5){\rouge{\line(0,-1){2}}}
\put(0,3){\rouge{\line(1,0){1}}}
\put(1,3){\rouge{\line(0,-1){1}}}
\put(1,2){\rouge{\line(1,0){2}}}
\put(3,2){\rouge{\line(0,-1){2}}}
\put(3,0){\rouge{\line(1,0){2}}}
\put(5,0){\rouge{\circle*{.2}}}  }}
 \multiput(0,0)(0,1){5}{\ligne(0,0){4}}
\put(3.9,5.6){\rotatebox{210}{$\rightsquigarrow$} {\small \underline{first\ return}}}
\end{picture} \hskip80pt
\begin{picture}(8,8)(0,0)
\put(.5,4.5){
  \put(0,4){\vcarre}
   \jligne(0,3){1}\put(1,3){\vcarre}
  \ligne(0,2){1}\jligne(1,2){1}\put(2,2){\vcarre}
  \ligne(0,1){2}\jligne(2,1){1}\put(3,1){\vcarre}
  \put(-0.1,0.5){
\put(0,3){\rouge{\circle*{.2}}}
\put(0,3){\rouge{\line(0,-1){1}}}
\put(0,2){\rouge{\line(1,0){1}}}
\put(1,2){\rouge{\line(0,-1){1}}}
\put(1,1){\rouge{\line(1,0){1}}}
\put(2,1){\rouge{\line(0,-1){1}}}
\put(2,0){\rouge{\line(1,0){1}}}   }}
\put(5,0){
  \jligne(0,4){1}
   \jligne(0,3){2}
  \ligne(0,2){1}\jligne(1,2){2}
 \ligne(0,1){3}\jligne(3,1){1}
  \ligne(0,0){3}\jligne(3,0){2}\put(-0.1,-0.5){
\put(0,5){\rouge{\line(0,-1){2}}}
\put(0,3){\rouge{\line(1,0){1}}}
\put(1,3){\rouge{\line(0,-1){1}}}
\put(1,2){\rouge{\line(1,0){2}}}
\put(3,2){\rouge{\line(0,-1){2}}}
\put(3,0){\rouge{\line(1,0){2}}}
\put(5,0){\rouge{\circle*{.2}}}  }}
 \multiput(0,0)(0,1){5}{\ligne(0,0){4}}
\put(0,-.5){$\underbrace{\hskip28pt}_j$}
\put(-5,1.5){$n-j\left\{\rule{0cm}{20pt}\right.$}
\end{picture}
\end{center}
\bigskip
prove the recursion:
  \begin{equation}\label{catalan_qrec}
      \mathcal{C}_n(q)=\sum_{j=1}^{n} q^{j(n-j)} \mathcal{C}_{j-1}(q) \mathcal{C}_{n-j}(q),
  \end{equation}    
  hence that $\mathcal{C}_n(q)$ lies in $\N[q]$. Can you see that this implies that
      $$ \mathcal{C}_n(1)=\textstyle{\frac{1}{n+1}\binom{2n}{n}}\qquad \hbox{(Catalan numbers)}?$$
 \end{exer}
\begin{exer}\rm{\bf (explore)} For $k<\ell$ and $n<m$, find an expression  for the $q$-enumeration of partitions $\mu$ such that $n^k\subseteq \mu\subseteq n^\ell$.
 \end{exer}
\begin{exer} \rm{\bf (a bit harder)} Prove that
     \begin{equation}\label{q_cat1}
   \textstyle{ C_n(q):=\frac{1}{[n+1]_q}\,\qbinom{2n}{n}\qquad \hbox{lies in}\qquad  \N[q]}.
   \end{equation}
   At least check that $\mathcal{C}_n(q)\not= C_n(q)$ in general. 
 \end{exer}
\begin{exer} \rm{\bf (open)} For any $a\leq b,c\leq d$ in $\N$ (no restriction between $b$ and $c$), such that $ad=bc$, show that
     \begin{equation}\label{q_foulkes}
   \textstyle{ \qbinom{b+c}{b}-\qbinom{a+d}{a} \in \N[q]}.
   \end{equation}
   You may check this in Sage. Can you prove it for some simple families of cases?
 \end{exer}
\begin{exer} \rm{\bf (open exploration)}  Find all the pairs of partitions of $n$ such that either the difference $P_{\rouge{\mu}}(q)-P_{\bleu{\nu}}(q)$ or $P_{\bleu{\nu}}(q)-P_{\rouge{\mu}}(q)$ lies in $\N[q]$. Experiments suggest that this is the case very often, so that it may be more interesting to find when it fails to hold. Observe that 
 exercise (d) reduces the number of cases to consider. It may be interesting the evident fact that when both $P_{\rouge{\mu}}(q)-P_{\bleu{\nu}}(q)$ and $P_{\bleu{\nu}}(q)-P_{\vert{\rho}}(q)$ lie in $\N[q]$, then so does $P_{\rouge{\mu}}(q)-P_{\vert{\rho}}(q)$.
 \end{exer}

\begin{exer} \rm In this exercise, we assume that we are in the context of the classical RSK correspondence (so that $b_i=i$). Show that for all semi-standard tableau $\tau$ we have
 $ (\varepsilon\leftarrow \rho(\tau)) =\tau$
with $\rho(\tau)$ standing for the reading word of $\tau$. For example
$$\left(\begin{picture}(2,1.5)(0,1)
    \Case{0,1}{2}
    \Case{0,0}{1} \Case{1,0}{3}
    \end{picture}
   \leftarrow\
     2\,1\,1\,3\,2\,4\,1\,2\ \right) =\
\begin{picture}(3,2)(0,1.2)
   \Case{0,2}{3} \Case{1,2}{3}
   \Case{0,1}{2}  \Case{1,1}{2}  \Case{2,1}{2}  \Case{3,1}{4}
   \Case{0,0}{1}  \Case{1,0}{1}  \Case{2,0}{1}  \Case{3,0}{1} \Case{4,0}{2}
\end{picture}$$
\end{exer}

\begin{exer}\label{exer_strip} \rm Consider the skew partition obtained as the difference $\mu/\lambda$
\begin{center}\smallersquarres
\begin{picture}(5,2)(0,0)
   \ligne(0,2){1}
   \jligne(0,1){2}\ligne(2,1){2}
   \jligne(0,0){4}\ligne(4,0){2}
\end{picture}
\end{center}
with $\mu$ equal to the shape of $(\tau\leftarrow v)$, for some semi-standard  tableau $\tau$ of shape $\lambda$:
 $$\left(\ 
 \begin{picture}(4,1.5)(0,1)
    \Case{0,1}{2}  \Case{1,1}{5} 
    \Case{0,0}{1}  \Case{1,0}{3}  \Case{2,0}{4}  \Case{3,0}{6}
\end{picture}\ \leftarrow\
1\,1\,2\,3\,5\right)\quad =\quad
\begin{picture}(3,2)(0,1.2)
     \MarkedCase{0,2}{5}
    \Case{0,1}{2}  \Case{1,1}{3} 
                          \MarkedCase{2,1}{4}  \MarkedCase{3,1}{6} 
    \Case{0,0}{1}  \Case{1,0}{1}  \Case{2,0}{1}  \Case{3,0}{2} 
                          \MarkedCase{4,0}{3}   \MarkedCase{5,0}{5} 
\end{picture}
$$ 
Show that there is a bijection between pairs $(\tau,v)$, with $\tau$ semi-standard tableau of shape $\lambda$, and $v$ non-decreasing sequence of $k$ intergers; and pairs $(\tau',\lambda)$, where $\tau'$ is a semi-standard tableau of shape $\mu$ such that $\mu/\lambda$ is an \define{horizontal
$k$ strip} (a skew tableau with $k$ cells, none of which lies on top of the other). 
\end{exer}

 \begin{exer}\rm  
The \define{pointwise sum}, $\mu+\lambda$, of partitions $\mu$ and $\lambda$, is the partition whose parts are $\mu_i+\lambda_i$ (adding zero parts if needed).
For $\mu$ having at most $n$ parts, find a multinomial formula for the number of standard tableaux of skew-shape $(\mu+1^n)/\mu$ (which is a vertical strip with $n$-cells, as illustrated in figure~\ref{vert_strip}).
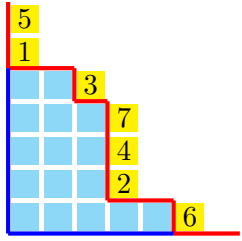
\begin{figure}[ht]
\begin{center}
 \begin{tikzpicture}[thick,scale=.44]
\entree07{5};
\entree06{1};
\pcellule05;\pcellule15;\entree25{3};
\pcellule04;\pcellule14;\pcellule24;\entree34{7};
\pcellule03;\pcellule13;\pcellule23;\entree33{4}; 
\pcellule02;\pcellule12;\pcellule22;\entree32{2}; 
\pcellule01;\pcellule11;\pcellule21;\pcellule31;\pcellule41;\entree51{6};
\draw[red,ultra thick] (0,7) to (0,5) ;
\draw[red,ultra thick] (0,5) to (2,5) ;
\draw[red,ultra thick] (2,5) to (2,4) ;
\draw[red,ultra thick] (2,4) to (3,4) ;
\draw[red,ultra thick] (3,4) to (3,1) ;
\draw[red,ultra thick] (3,1) to (5,1) ;
\draw[red,ultra thick] (5,1) to (5,0) ;
\draw[red,ultra thick] (5,0) to (7,0) ;
\draw[blue,ultra thick] (0,5) to (0,0) ;
\draw[blue,ultra thick] (0,0) to (5,0) ;
\end{tikzpicture}
\end{center}
\caption{A standard tableau of shape $(\mu+1^n)/\mu$, with $\mu=53332$ and $n=7$.}
\label{vert_strip}
\end{figure}
\end{exer}

 \begin{exer}\rm 
Show that the Kostka matrix $(K_{\lambda\mu})_{\lambda,\mu\vdash n}$ is upper unitriangular, if partitions are ordered with any linear extension of the dominance order. If you use the \define{lexicographic} order, {\it e.g.}: $4,\ 31,\ 22,\ 211,\ 1111$ for $n=4$, you get: 
  $$ \begin{pmatrix}
1 & 1 & 1 & 1 & 1 \\
0 & 1 & 1 & 2 & 3 \\
0 & 0 & 1 & 1 & 2 \\
0 & 0 & 0 & 1 & 3 \\
0 & 0 & 0 & 0 & 1
\end{pmatrix}$$
\end{exer}

 \begin{exer}\rm Show that, for $\mu\vdash n$,
     $$f^\mu=n!\, \det\left(a_{ij}\right),\qquad {\rm where}\qquad a_{ij}={1}/{(\mu_i+j-i)}.$$
  \end{exer}
  
 \begin{exer}\rm {\bf(open)} Explore the value of $\sum_{\mu\vdash n} (f^\mu)^k$, for $k\geq 2$.
\end{exer}

\end{chapter}


\begin{chapter}{\bleu{``Classical'' Symmetric functions}}
All of the following (and much more) may be found in Macdonald's book~\cite{macdonald}, Sagan's~\cite{sagan}, or Stanley's~\cite{stan_comb}. See also~\cite{manivel} for links to algebraic geometry. We mostly use Macdonald's notation.

The usual Sage tutorial for symmetric functions is available \href{http://doc.sagemath.org/html/en/reference/combinat/sage/combinat/sf/sfa.html}{here}.
An improved  \href{https://more-sagemath-tutorials.readthedocs.io/en/latest/tutorial-symmetric-functions.html}{Sage Symmetric Function Tutorial}\footnote{https://more-sagemath-tutorials.readthedocs.io/en/latest/tutorial-symmetric-functions.html} was
prepared by  Pauline Hubert and Mélodie Lapointe (UQAM). It is more adapted to our presentation. In Maple, one typically uses John Stembridge's package \href{http://www.math.lsa.umich.edu/~jrs/maple.html\#SF}{SF}.

\section{Basic notions}
We write $\bm{x}$ for the list (vector) $n$ variables $(x_1,x_2,\ldots, x_n)$, and use ``vector'' notation for monomials:
      $$\bm{x}^{\bm{a}}=x_1^{a_1} x_2^{a_2} \cdots x_n^{a_n},\qquad {\rm for}\qquad \bm{a}=(a_1,a_2,\ldots, a_n)\in\N^n.$$      
The \define{degree} of  $\bm{x}^{\bm{a}}$ is $\deg(\bm{x}^{\bm{a}})=|\bm{a}|=a_1+a_2+\ldots+a_n$.  Thus for a scalar $q$, we have $(q\cdot \bm{x})^{\bm{a}}= q^d\,(\bm{x}^{\bm{a}})$ with $d= \deg(\bm{x}^{\bm{a}})$. For certain, $q\cdot \bm{x}:=(qx_1,qx_2,\ldots,q x_n)$.
It is well known that there are $\binom{n+d-1}{d}$ monomials of degree $d$
in $n$ variables. The ring of polynomials in $n$-variables, $\mathcal{R}:=\Q[\bm{x}]$, is \define{graded} by degree
   $$\mathcal{R}=\bigoplus_{d\in\N} \mathcal{R}_d,\qquad {\rm with}\qquad \mathcal{R}_d=\Q\{\bm{x}^{\bm{a}}\ |\ \deg(\bm{x}^{\bm{a}})=d\}.$$
 The $\binom{n+d-1}{d}$-dimensional $\Q$-vector space $\mathcal{R}_d$ is the $d$-degree \define{homogeneous component} of $\mathcal{R}$. Equivalently, we have
    $$p(\bm{x})\in \mathcal{R}_d,\qquad {\rm iff}\qquad p(q\,\bm{x})=q^d\,p(\bm{x}).$$
 We then say that $p(\bm{x})$ is \define{homogeneous} of degree $d$.
Any polynomial $p(\bm{x})$ is uniquely decomposed as a sum of homogeneous polynomials, {\it i.e.}: 
     $$p(\bm{x})=p_0(\bm{x})+p_1(\bm{x})+p_2(\bm{x})+\ldots+p_n(\bm{x}),$$
where $n$ is the degree of $p(\bm{x})$, and $p_d(\bm{x})$ lies in $\mathcal{R}_d$. We say that $p_d(\bm{x})$ is the degree $d$ \define{homogeneous component} of $p(\bm{x})$.

\section{From symmetric polynomials to symmetric functions}
Recall that \define{symmetric polynomials} are polynomials $f(\bm{x})$ in $\mathcal{R}$ such that, for any permutation
$\sigma$ of the set $\{1,2,\ldots,n\}$ (or $\sigma\in\S_n$), we have
    $$f(x_{\sigma(1)},x_{\sigma(2)},\ldots,x_{\sigma(n)})= f(x_1,x_2,\ldots,x_n).$$
 The set of symmetric polynomials forms a graded subring of $\mathcal{R}$, denoted by $\Lambda$, and $\Lambda=\bigoplus_{d\in\N} \Lambda_d$.
 Indeed, It is clear that the sums and products of symmetric polynomials are also symmetric polynomials. Moreover, if $f$ is a symmetric polynomial, so are all of its homogeneous components. In other words, $\Lambda_d=\mathcal{R}_d\cap \Lambda$.
A linear basis of $\Lambda_d$ is formed by the set of (non-vanishing) \define{monomial symmetric} polynomials
    $$\bleu{m_\lambda(\bm{x}):=\sum_{\bm{a}\in \mathfrak{r}(\lambda)}  \bm{x}^{\bm{a}}},\qquad \lambda\vdash d$$
with $\bm{a}$ varying in the set $\mathfrak{r}(\lambda)$ of rearrangements of the length $n$ vector $(\lambda_1,\ldots,\lambda_k,0,\ldots,0)$. For example
      $$m_{211}(x_1,x_2,x_3)= x_1^2x_2x_3+x_1x_2^2x_3+x_1x_2x_3^2,$$
Observe that our definition forces $m_\lambda=0$ when  $\ell(\lambda)>n$. As we will see, this makes it so that identities within $\Lambda_d$ are simpler to describe  when $n>d$.  Moreover, they specialize with no surprises when one restricts the number of variables. For this reason, and after some initial care, it is customary to send $n$ to infinity, and consider a denumerable set of variables $\bm{x}=(x_1,x_1,x_3,\ldots). $
 The resulting ``polynomials'' in infinitely many variables are then called \define{symmetric functions}. In fact, most of the time we will forego actual mention of the variables, and say for example that $\{m_\lambda\}_{\lambda\vdash d}$ forms a linear basis of $\Lambda_d$. Thus the dimension of $\Lambda_d$ is the number of partitions of $d$.

Among the interesting bases of  $\Lambda_d$, one first considers the following three ``multiplicative'' bases, all indexed by partitions $\mu=(\mu_1,\mu_2,\ldots, \mu_k)$ of $d$:
\begin{itemize}
\item[$\bullet$] The \define{complete homogeneous} symmetric functions are defined to be the product $h_\mu:=h_{\mu_1}h_{\mu_2}\cdots h_{\mu_k}$, with simply indexed $h_d$ defined as:
                          $$\bleu{h_d:=\sum_{\lambda\vdash d} m_\lambda}.$$
For example, we have $h_1=m_1$; 
\begin{align*}
 &h_2=m_{11}+m_2, &&h_{11}=h_1^2=2\,m_{11}+m_2;\\
 &h_3=m_{111} + m_{21} + m_{3},&&h_{21}=3\,m_{111} + 2\,m_{21} + m_{3},&&h_{111}=6\,m_{111} + 3\,m_{21} + m_{3};
 \end{align*}
 and
\begin{align*}
 &h_4=m_{4}+ m_{31} + m_{22} + m_{211} +m_{1111}  ,\\
& h_{31}=m_{4} + 2\,m_{31} +2\,m_{22} + 3\,m_{211} +4\,m_{1111}    \\
 &h_{22}=m_{4} + 2\,m_{31} + 3\,m_{22} + 4\,m_{211} + 6\,m_{1111} ,\\
 &h_{211}= m_{4} + 3\,m_{31} + 4\,m_{22} + 7\,m_{211} + 12\,m_{1111},\\
 &h_{1111}=m_{4} + 4\,m_{31} + 6\,m_{22} + 12\,m_{211} + 24\,m_{1111}.
\end{align*}
 \item[$\bullet$] The \define{elementary} symmetric functions defined as the product \bleu{$e_\mu:=e_{\mu_1}e_{\mu_2}\cdots e_{\mu_k}$}, with simply indexed $e_d$ defined as 
 \bleu{$e_d:= m_{1^d}$}. Thus $e_1=m_1=h_1$, implying that $e_{1^d}=h_{1^d}$. Other values (not already explicitly given) are $e_{21}= 3\,m_{111} + m_{21}$, and
\begin{align*}
  &e_{211}=m_{31}  + 2\,m_{22} +  5\,m_{211} + 12\,m_{1111} ,\\
  &e_{22}= \qquad\qquad m_{22} + 2\,m_{211} + 6\,m_{1111},\\
  &e_{31}=\qquad \qquad\qquad\qquad  m_{211} + 4\,m_{1111}.
  \end{align*}
\item[$\bullet$] The \define{power sum} symmetric functions defined as the product \bleu{$p_\mu:=p_{\mu_1}p_{\mu_2}\cdots p_{\mu_k}$}, with simply indexed $p_d$ defined as 
\bleu{$p_d:= m_d$}. Again, $p_1=m_1=h_1=e_1$, so that $p_{1^d}=h_{1^d}=e_{1^d}$. Other values are $p_{21}= m_{21} + m_{3}$, and
\begin{align*}
&p_{31}=m_{4} + m_{31},\\
 &p_{22}=m_{4} + 2\,m_{22},\\
 &p_{211}= m_{4} + 2\,m_{31}  + 2\,m_{22}  + 2\,m_{211}.\\
 \end{align*}
\end{itemize}
There are nice combinatorial descriptions for all the transition matrices between the four bases above ({\it i.e.}, $m_\mu$, $h_\mu$, $e_\mu$, and $p_\mu$). In some instances these come out of the following identities, best coined in terms of generating series in an extra variable $z$. First, the above definition is equivalent to (why?)
\begin{equation}\label{defn_H}
  \bleu{H(z):=\sum_{d\geq 0} h_d\,z^d =\prod_{i\geq 1} \frac{1}{1-x_i\,z}}\qquad {\rm and}\qquad 
  \bleu{E(z):=\sum_{d\geq 0} e_d\,z^d =\prod_{i\geq 1} {1+x_i\,z}}
  \end{equation}
It is thus immediate that $   H(z)\,E(-z)=1$, and a simple calculation gives
\begin{equation}\label{lienHP}
{\rm a)}\  H(z)=\exp\big(P(z) \big),\qquad {\rm and}\qquad 
  {\rm b)} \ E(z)=\exp\big(-P(-z)\big),
\end{equation} 
if one sets \bleu{$P(z):=\sum_{d\geq 1} p_d\,z^d/d$}. Comparing coefficients of $z^d$ in the above identities gives:
  \begin{align}\label{lien_h_e}
     {\rm a)}\ \sum_{k=0}^d (-1)^k\,h_{d-k}e_k=0,&&{\rm b)}\  h_d=\sum_{\mu\vdash d}\frac{p_\mu}{z_\mu},&& {\rm c)}\  e_d=\sum_{\mu\vdash d}(-1)^{d-\ell(\mu)} \frac{p_\mu}{z_\mu},
    \end{align}
 We can recursively solve (\ref{lien_h_e}.a) either to write the $h_k$ in terms of $e_j$, or vice-versa. 
Also, from~(\ref{lienHP}.a) we get $H'(z )={P'(z )\,H(z )}$, from which it follows that
\begin{equation} 
      d\,h_d=\sum_{r=1}^d p_r\,h_{d-r}.
  \end{equation} 
Using Cramer's rule to solve this system of equations for the indeterminate $p_k$, one finds the following expansion in terms of the $h_d$'s:
\begin{equation}\label{p_en_h}
  p_k=(-1)^{k-1}\,\det\left|   \begin{matrix}h_1&1&0&\ldots&0\\ 
                                    2\,h_2&h_1&1&\ldots &0\\ 
                                    3\,h_3&h_2&h_1&\ldots &0\\ 
                                    \vdots &\vdots&\vdots&\ddots&\vdots\\ 
                                   k\,h_k &h_{k-1}& h_{k-2}&\ldots&h_1  
                                   \end{matrix}\right|.
\end{equation}
Similarly, we get $d\,e_d=\sum_{r=1}^d (-1)^{r-1} p_r e_{d-r}$ from~(\ref{lienHP}.b), and derive
\begin{equation}\label{p_en_e}
p_k=\det\left|   \begin{matrix}e_1&1&0&\ldots&0\\ 
                                    2\,e_2&e_1&1&\ldots &0\\ 
                                    3\,e_3&e_2&e_1&\ldots &0\\ 
                                    \vdots &\vdots&\vdots&\ddots&\vdots\\ 
                                   k\,e_k &e_{k-1}& e_{k-2}&\ldots&e_1  \end{matrix}\right|
 \end{equation}
The \define{omega} involution is the linear and multiplicative operator that sends $p_d$ to $(-1)^{d-1}p_d$. Clearly, for $\mu\vdash d$, 
     $$\bleu{\omega(p_\mu)=(-1)^{d-\ell(\mu)}p_\mu}.$$
 From the above, one sees easily that
 $\omega(h_d)=e_d$, or equivalently that $\omega(e_d)=h_d$ (since $\omega^2=\Id$).

\section{Schur functions, a combinatorial approach}\label{sec_schur_comb}
The Schur function $s_\lambda(\bm{x})$ is combinatorially defined as:
   \begin{equation}\label{schur_combin}
       \bleu{s_\lambda(\bm{x}):=\sum_{\tau} \bm{x}_\tau},
    \end{equation}
with the sum being over all semi-standard tableaux 
$\tau:\lambda\rightarrow \N$,
and  $\bm{x}_\tau:=\prod_{c\in \lambda} x_{\tau(c)}$.
 This definition is naturally extended to skew shapes $\lambda/\mu$ to get \define{skew Schur} polynomials $s_{\lambda/\mu}(\bm{x}):=\sum_{\tau} \bm{x}_\tau$,
with the sum over semi-standard tableaux of shape $\lambda/\mu$. 
If one assumes (as is shown in the exercises) that $s_\lambda(\bm{x})$ is a symmetric polynomial, then Equation (\ref{schur_combin})  implies that
 \begin{equation}\label{s_en_m}
      s_\lambda(\bm{x})=\sum_{\mu\vdash n} K_{\lambda\mu}\,m_\mu(\bm{x}),
 \end{equation}
In particular, we always have $s_{(n)}=h_n$, and $s_{1^n}=e_n$. One may show (see exercise) that $\omega s_\mu=s_{\mu'}$.

Often, it is natural to consider a tableau being filled with the variables $x_i$, rather than their indices $i$. With this point of view,  the evaluation monomial is simply the product of all entries of the tableau. Then one may consider $s_\lambda(\bleu{\bm{x}}+\rouge{\bm{y}})$ as the evaluation of the Schur function $s_\lambda$ on the \define{variables} 
 $\bleu{\bm{x}}+\rouge{\bm{y}}=(\bleu{x_1+x_2+x_3+\ldots})+(\rouge{y_1+y_2+y_3+\ldots})$,
 with the convention that all $\rouge{y_k}$ are larger than all $\bleu{x_i}$. It follows that (also see section~\ref{plethysme} on plethysm)
  \begin{equation}\label{schur_sum}
             s_\lambda(\bleu{\bm{x}}+\rouge{\bm{y}})=\sum_{\mu\subseteq \lambda} s_\mu(\bleu{\bm{x}})\,s_{\lambda/\mu}(\rouge{\bm{y}}),
\end{equation}
since any semi-standard tableau $\tau:\lambda\rightarrow \bleu{\bm{x}}+\rouge{\bm{y}}$,
can naturally be separated into the semi-standard tableau corresponding to the sub-partition $\mu\subseteq\lambda$, where only the $\bleu{x_i}$ appear, and the skew semi-standard tableau of shape $\lambda/\mu$, where only the $y_k$ appear.

For the special cases $\lambda=(n)$ and $\lambda=1^n$, we get
\begin{equation}\label{h_somme}
   h_n(\bleu{\bm{x}}+\rouge{\bm{y}})=\sum_{k=0}^n h_k(\bleu{\bm{x}})\,h_{n-k}(\rouge{\bm{y}}),\qquad{\rm and}\qquad 
e_n(\bleu{\bm{x}}+\rouge{\bm{y}})=\sum_{k=0}^n e_k(\bleu{\bm{x}})\,e_{n-k}(\rouge{\bm{y}}).
\end{equation}
Exploiting the RSK correspondence, we find interesting formulas. Indeed, associate to
     $$w=\begin{pmatrix}
      b_1& b_2& \ldots & b_d\\
      a_1& a_2&\ldots &a_d
          \end{pmatrix},$$
the product of degree $d$ monomials  $\bleu{\bm{x}}_{\bm{a}} \rouge{\bm{y}}_{\bm{b}}$, with $\bm{a}=(a_1,\ldots,a_d)$ and $\bm{b}=(b_1,\ldots,b_d)$.
Thus we identify the set of length $d$ lexicographic two row matrices, with the terms of the sum $\sum_{|\bm{a}|=|{\bm{b}|=d}}\bleu{\bm{x}}_{\bm{a}} \rouge{\bm{y}}_{\bm{b}}$. This sum is readily seen to be the same as degree $d$ homogenous component of
   $$\sum_{d\geq 0} h_d(\bleu{\bm{x}}\,\rouge{\bm{y}})=\sum_{i,j\geq 1}\frac{1}{1-\bleu{x_i}\rouge{y_j}},$$
 with $h_d$ expanded in the ``variables'' $\bleu{x_i}\rouge{y_j}$.
 
We now apply the RSK correspondence to each monomial $\bleu{\bm{x}}_{\bm{a}} \rouge{\bm{y}}_{\bm{b}}$ (translated back into a lexicographic two row matrices) to get a pair of same shape semi-standard tableaux $(P,Q)$, for which we have
$\bleu{\bm{x}}_{\bm{a}} \rouge{\bm{y}}_{\bm{b}}=\bleu{\bm{x}}_P\rouge{\bm{y}}_Q$.
We may rewrite this sum $\sum_{\bm{a},_{\bm{b}}}\bleu{\bm{x}}_{\bm{a}} \rouge{\bm{y}}_{\bm{b}}$ as $\sum_{P,Q} \bleu{\bm{x}}_P\rouge{\bm{y}}_Q$,
with the sum now being over the pairs of semi-standard tableaux of same shape. It follows that we have the identities
   \begin{equation}\label{nouyau_schur}
        h_k(\bleu{\bm{x}}\rouge{\bm{y}}) = \sum_{\lambda\vdash k} s_\lambda(\bleu{\bm{x}})
            s_\lambda(\rouge{\bm{y}}).
 \qquad{\rm and}\qquad         e_k(\bleu{\bm{x}}\rouge{\bm{y}}) = \sum_{\lambda\vdash k} s_{\lambda'}(\bleu{\bm{x}})
            s_\lambda(\rouge{\bm{y}}).
     \end{equation}
applying $\omega$ to the $\bleu{\bm{x}}$ component, to get the second from the first.

Now, summing up both sides of (\ref{nouyau_schur}) for $n$, we get Cauchy-Littlewood's formula
\begin{equation} \label{cauchy_littlewood}
\prod_{i,j\geq 1}\frac{1}{1-\bleu{x_i}\,\rouge{y_j}} =\sum_\lambda s_\lambda(\bleu{\bm{x}})s_\lambda(\rouge{\bm{y}}).
\end{equation}
If we restrict this argument to matrices of the form
     $ w=\begin{pmatrix}
      1& 2&\ldots &k\\[-3pt]
       a_1& a_2& \ldots & a_k
          \end{pmatrix}$,
and consider the evaluation monomial of $w$ to be $\bleu{\bm{x}}_{\bm{a}}$, and the evaluation monomial of a pair $(P,Q)$ to be $\bleu{\bm{x}}_P$.  We get the formula (see exercises)
   \begin{equation}\label{frob_regulard}
       h_1^n=\sum_{\lambda\vdash n} f^\lambda\, s_\lambda(\bleu{\bm{x}})
     \end{equation}
with $f^\lambda$ standing for the number of standard tableaux of shape $\lambda$.   

Another interesting formula  (see Stanley~\cite[Proposition~7.19.11]{stan_comb}) says that
\begin{equation}\label{principal_spec}
  s_{\lambda/\mu}(1,q,q^2,q^3,\ldots)=\frac{K_{\mu,1^n}(q)}{(1-q)(1-q^2)\cdots (1-q^n)}.
\end{equation}

\section{Dual basis and Cauchy kernel}\label{sec_dual}
On $\Lambda$, we consider the \define{Hall scalar product} characterized by
 \begin{equation}\label{defn_scalar}
     \bleu{ \langle p_{\lambda},p_{\mu} \rangle := \begin{cases}
   z_{\lambda}, & \text{if}\ \lambda=\mu,\\
   0, & \text{if}\ \lambda\not=\mu.
                                   \end{cases}}
  \end{equation}
Which is readily seen to be preserved by $\omega$, {\it i.e.}:
$\langle \omega(f),\omega(g)\rangle =\langle f,g\rangle$. We may calculate directly that
    $$ \langle h_n,h_n \rangle= \langle e_n,e_n \rangle=\sum_{\mu\vdash n}\frac{1}{z_\mu}.$$ 
In fact, this sum is equal to $1$ (see exercise). 

The Schur functions may also be obtained via the Gram-Schmidt orthogonalization  process applied to the basis of monomial symmetric $\{m_\lambda\}_{\lambda\vdash n}$, written in increasing lexicographic order of partitions. Indeed, the Schur polynomials are uniquely determined\index{Schur polynomials!orthogonal} by the following two properties:
\begin{equation}\label{charac_scalar}
\begin{array}{llll}
   (1)&  \langle s_\lambda,s_\mu \rangle =0,\quad \mathrm{whenever}\quad \lambda\not=\mu,\\[6pt]
   (2)& \displaystyle s_\lambda=
       m_\lambda+\sum_{\mu\prec \lambda} c_{\lambda\mu} m_\mu,
\end{array}
\end{equation}
for some coefficients $c_{\lambda\mu}$ to be calculated so that the first property holds.
Although it is not immediate, this makes it so that 
\begin{equation}\label{schur_normal}
     \langle s_\mu,s_\mu \rangle =1,
\end{equation}
for all partitions $\mu$.
One way to ``see'' this is through the following discussion.
We start with the simple computation
  \begin{eqnarray}\prod_{i,j\geq 1}\frac{1}{1-\bleu{x_i}\,\rouge{y_j}} &=&  \exp\sum_{i,j} -\log(1-\bleu{x_i}\,\rouge{y_j}) \nonumber\\
   &=&  \exp\sum_{i,j} \sum_{k\geq 1} (\bleu{x_i}\,\rouge{y_j})^k/k \nonumber\\
   &= & \exp\sum_{k\geq 1} p_k(\bleu{\bm{x}})\,p_k(\rouge{\bm{y}})/k \nonumber\\
   &=&\sum_\lambda p_\lambda(\bleu{\bm{x}})\frac{ p_\lambda(\rouge{\bm{y}})}{z_\lambda}.\label{dual_power}
\end{eqnarray}
The left-hand side of (\ref{cauchy_littlewood}) and (\ref{dual_power}) is 
the \define{Cauchy kernel}\index{Cauchy kernel}, henceforth denoted by $ \Omega(\bleu{\bm{x}}\rouge{\bm{y}})$.  This is a symmetric series in the variables ``$x_iy_j$'' (as will become clearer in section~\ref{plethysme}). Its degree $n$ homogeneous component (with $x_iy_j$ being considered as having degree $1$) is
\begin{equation}\label{hn_xy}
    h_n(\bleu{\bm{x}}\rouge{\bm{y}})=\sum_{\lambda\vdash n} p_\lambda(\bleu{\bm{x}})\frac{ p_\lambda(\rouge{\bm{y}})}{z_\lambda}.
   \end{equation}
Various similar formulas characterize  \define{dual basis pairs}  $\{u_\lambda\}_{\lambda}$ and $\{v_\mu\}_{\mu}$. These are pairs of bases such that
\begin{equation}\label{scal_dual}
   \langle u_\lambda,v_\mu\rangle=\begin{cases}
       1, & \text{if}\ \lambda=\mu,\\
       0, & \text{otherwise}.
       \end{cases}
  \end{equation}
The relevant statement here is that the set conditions  (\ref{scal_dual}) is equivalent to the single identity
\begin{equation}\label{noyeau_dualite}
  \Omega(\bleu{\bm{x}}\rouge{\bm{y}})=\sum_{n\geq 0}   h_n(\bleu{\bm{x}}\rouge{\bm{y}})=\sum_\lambda u_\lambda(\bleu{\bm{x}})\,v_\lambda(\rouge{\bm{y}}).
\end{equation}
You are asked to show this in an exercise.
In particular, identities (\ref{cauchy_littlewood}) and (\ref{dualite_homogene}) are equivalent to $ \langle s_\lambda,s_\mu\rangle =\delta_{\lambda,\mu}$, and $\langle m_\lambda,h_\mu\rangle =\delta_{\lambda,\mu}$, where $\delta_{\lambda,\mu}$ is the usual\footnote{Which is so usual that one always seems to be required to mention the fact.}  \define{Kronecker delta function}\index{Kronecker!delta function}. 
In other words,  $\{s_\lambda\}_\lambda$ is a self dual orthonormal basis. Exploiting the duality between $\{m_\lambda\}_\lambda$ and $\{h_\mu\}_\mu$, and equality (\ref{s_en_m}), we get
\begin{equation}\label{h_en_s}
   h_\mu=\sum_\lambda K_{\lambda,\mu} s_\lambda
\end{equation}
Comparing this with Equation~(\ref{s_en_m}) we notice that the index of summation is now $\lambda$ rather than $\mu$.

\section{Product, Kronecker product, and Schur positivity}
Since the product of two Schur functions is also a symmetric function, it may be expanded in terms of the Schur basis. In formula, for $\mu\vdash m$ and $\nu\vdash n$, we have
\begin{equation}\label{LR_coeff}
    s_\mu s_\nu= \sum_{\lambda\vdash (m+n)} c_{\mu\nu}^{\lambda} s_\lambda.
 \end{equation}
It is a fact that the coefficients $ c_{\mu\nu}^{\lambda}$, known as the \define{Littlewood-Richardson coefficients}, are all positive integers\footnote{They are hard to calculate, and in fact closely tied to deep questions in Geometric Complexity Theory: a.k.a. the study of the (VP vs VNP) problem.}. They also occur in the Schur expansion of skew Schur functions, {\it i.e.}: for $\lambda\vdash (m+n)$ and $\mu\vdash m$, we have
\begin{equation}\label{Skew_LR}
    s_{\lambda/\mu}= \sum_{\nu\vdash n} c_{\mu\nu}^{\lambda} s_\lambda.
 \end{equation}
Taken together, \pref{LR_coeff} and~\pref{Skew_LR} imply that, for all $\mu$, $\nu$ and $\lambda$, one has
 \begin{equation}
    \langle s_\mu s_\nu,s_\lambda\rangle = \langle s_\nu,s_{\lambda/\mu}\rangle.
 \end{equation}
This makes it natural to consider, for any symmetric function $f$, the \define{skewing operator}, denoted by $f^\perp$, which is simply the adjoint of the multiplication operator $f\cdot(-)$, with respect to the scalar product. In formula, since the set of Schur functions forms a linear basis, it is characterized by the property
 \begin{equation}
    \bleu{\langle f\cdot  s_\nu,s_\lambda \rangle = \langle s_\nu,f^\perp s_\lambda\rangle}.
 \end{equation}
In particular,  $s_\mu^\perp s_\lambda = s_{\lambda/\mu}$.

There is a nice combinatorial rule (not presented here) due to Littlewood and Richardson to express the $c_{\mu\nu}^{\lambda}$ as enumerators for certain semi-standard tableaux. Two special cases are the following \define{Pieri formulas}, where all coefficients are either $0$ or $1$:
\begin{align}\label{h_pieri}
   h_k \, s_\mu =\sum_{\theta} s_\theta,\qquad ({\rm Pieri})&& e_k \, s_\mu =\sum_{\theta} s_\theta\qquad (\hbox{dual Pieri}),
 \end{align}
with the sum being over  partitions $\theta$ for which $\theta/\mu$ is a $k$-cell {horizontal strip} (see exercise~\ref{exer_strip}) in the first formula,  and over those for which $\theta/\mu$ is a $k$-cell  \define{vertical strip} (no two cells in the same row)  in the second. One may use these Pieri rules to ``calculate'' any $s_\mu$ (see exercise).

The \define{Kronecker product} of symmetric functions, denoted by $f*g$, is the unique bilinear product (commutative and associative) such that
  \begin{equation}\label{def_kronecker}
     \bleu{p_\lambda * p_\mu:=\begin{cases}
           z_\lambda\,p_\lambda, & \text{if}\ \lambda=\mu, \\
           0, & \text{otherwise}.
\end{cases} }
\end{equation}
Observe that 
  \begin{equation}\label{kronecker}
     h_n * f:=\begin{cases}
           f, & \text{if}\ f\in\Lambda_n, \\
           0, & \text{otherwise};
         \end{cases} \qquad {\rm and}\qquad  
   e_n * f:=\begin{cases}
           \omega(f), & \text{if}\ f\in\Lambda_n, \\
           0, & \text{otherwise}.
         \end{cases}
\end{equation}
Hence, $H=\sum_n h_n$ acts as identity for this product; and multiplication by $E=\sum_n e_n$ corresponds to $\omega$. We easily check that 
$\omega(f)* g= f* \omega(g)$.  Calculating $s_\mu*s_\lambda$ is hard.

Littlewood showed  the nice identity 
 \begin{equation}\label{kronecker_littlewood}
    (s_\alpha s_\beta) * s_\theta = \sum_{\lambda,\mu} c_{\lambda\mu}^\theta
         (s_\alpha*  s_\lambda)\, (s_\beta* s_\mu),
   \end{equation}
 with the sum being over the pairs $\lambda,\mu$ such that: $c_{\lambda\mu}^\theta$ is not zero, $\lambda \vdash |\alpha|$, and $\mu\vdash |\beta|$. 

A common feature of all of these is that they have expansions with positive integer coefficients in terms of the Schur function basis. When this is the case, we say that we have \define{Schur positivity}\index{Schur!positivity}, or that the function in question is \define{Schur positive}\index{Schur-positive}. Illustrating with instances of operations discussed above, we have
\begin{align*}
    &s_{21}^\perp\, s_{4321}=s_{{43}}+2\,s_{{421}}+s_{{4111}}+2\,s_{{331}}+2\,s_{{322}}+2\,s_{{3211}}+s_{{2221}},\\
    &s_{3111}*s_{321}=s_{51}+2\,s_{42}+2\,s_{411}+s_{33}+4\,s_{321}+2\,s_{3111}
                       +s_{222}+2\,s_{2211}+s_{21111},\\
    &s_{21}s_{32}=s_{53}+s_{521}+s_{44}+2\,s_{431}+s_{422}+s_{4211}+s_{332}+s_{3311}+s_{3221}.
\end{align*} 
The analogous \underline{stronger} notions of \define{$h$-positive} and \define{$e$-positive} are also considered. In view of~\pref{s_en_m}, it is clear that we have monomial positivity for any Schur positive functions. The global result (see exercises) is that
\begin{theorem}
 For any Schur positive symmetric functions $f$ and $g$, the symmetric functions $f+g$, $f\,g$, $f*g$, and $f^\perp g$ are also Schur positive. Any $h$-positive (or $e$-positive) symmetric function is also Schur positive.
\end{theorem}

In fact, there is a nice theorem that makes apparent that Schur positivity is a rare phenomenon. Foregoing a technical description\footnote{One must make precise the sampling space, see exercises.}, it may be stated as follows.
\begin{theorem}[B-Patrias-Reiner 2017]\label{thm_BPR}
The probability of choosing at random (with uniform distribution) a Schur positive symmetric function in the set of homogeneous symmetric functions of degree $n$ which are monomial positive, is equal to
    \begin{equation}
        \frac{1}{\prod_{\mu\vdash n}\sum_{\lambda\vdash n}  K_{\lambda\mu}},
    \end{equation}
  with $K_{\lambda\mu}$ the Kostka numbers.
\end{theorem}
For $n=1,2,3,\ldots$, the first terms of the sequence of denominators are
  $$1,2,9,560,480480,1027458432000,2465474364698304960000,\ldots $$
making the point that Schur positivity is indeed rare. There is a similar statement making the point that $h$-positivity (or $e$-positivity) is rare among Schur positivity functions.

For symmetric functions with coefficients in $\Q(q,t)$ (such as we will consider in the sequel), it is natural to say that we have Schur positivity when functions have Schur expansions with coefficients in $\N[q,t]$.

\section{Transition matrices}
The matrices that express the possible changes of basis between the 6 fundamental bases that have been discussed may be described as follows,  borrowing figure~\ref{fig_transition} from Macdonald~\cite[page 104]{macdonald}. Each oriented edge encodes a transition matrix between bases, with bases labeling the vertices. Thus the arrow labeled $L=(L_{\mu\lambda})$ stands for the matrix describing how write the $p_\mu$'s in terms of the $h_\lambda$, 
$K$ stands for the Kostka matrix, and $J$ for the matrix that corresponds to the involution $\omega$. As usual we write $M^{\scriptscriptstyle\mathrm{Tr}}$ for the transposed  of  the matrix $M$. By composition of the red arrows (or their inverse), all unlabeled edges may be filled in to complete the graph.
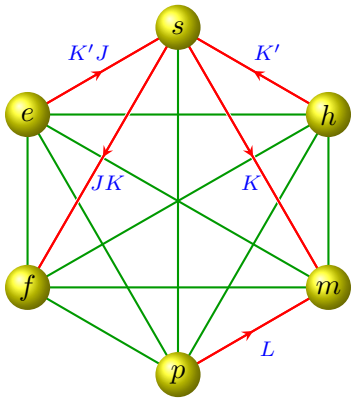
\begin{figure}[ht]\setlength{\unitlength}{4mm}
\begin{center}
\begin{tikzpicture}
     \coordinate (H) at (2,1.15) ; \coordinate (S) at (0,2.31) ;
    \coordinate (E) at (-2,1.15) ; \coordinate (F) at (-2,-1.15) ;
    \coordinate (P) at (0,-2.31) ;  \coordinate (M) at (2,-1.15) ;
    \draw[red,thick,->] (H) -- (S) node[midway,above] {\bleu{$\scriptstyle \quad K'$}};
    \draw[red,thick,->] (E) -- (S) node[midway,above] {\bleu{$\scriptstyle K'J\quad $}};
\draw[darkgreen,thick] (H) to (E) ; \draw[darkgreen,thick] (E) to (M) ; 
\draw[darkgreen,thick] (S) to (P) ; \draw[darkgreen,thick] (P) to (H) ; 
\draw[darkgreen,thick] (H) to (F) ;
\draw[darkgreen,thick] (P) to (F) ; 
\draw[darkgreen,thick] (P) to (E) ; \draw[darkgreen,thick] (F) to (M) ; 
\draw[darkgreen,thick] (F) to (E) ; \draw[darkgreen,thick] (M) to (H) ; 
\draw[red,directed,thick] (E) to (S) ; \draw[red,directed,thick] (H) to (S) ; 
    \draw[white,ultra thick] (S) to (M) ; 
          \draw[red,thick,->] (S) -- (M) node[midway,below=.1] {\bleu{$\scriptstyle K\,$}};
          \draw[red,directed,thick] (S) to (M) ; 
    \draw[white,ultra thick] (F) to (S) ; 
          \draw[red,thick,->] (S) -- (F) node[midway,below=.1] {\bleu{$\scriptstyle \  JK$}};
          \draw[red,directed,thick] (S) to (F) ; 
   \draw[white,ultra thick] (P) to (M) ; 
          \draw[red,thick,->] (P) -- (M) node[midway,below] {\bleu{$\scriptstyle\quad L$}};
           \draw[red,directed,thick] (P) to (M) ; 
\shade[ball color=yellow] (H) circle (0.3);\node at (H) {$h$};
 \shade[ball color=yellow]  (E) circle (0.3); \node at (E) {$e$};
 \shade[ball color=yellow]  (P) circle (0.3); \node  at (P) {$p$}; 
 \shade[ball color=yellow] (S) circle (0.3);\node  at (S) {$s$}; 
 \shade[ball color=yellow]   (F) circle (0.3); \node  at (F) {$f$};
 \shade[ball color=yellow]  (M) circle (0.3); \node at (M) {$m$};
\end{tikzpicture}
\end{center}
\vskip-15pt
\caption{Transition matrices.}
\label{fig_transition}
\end{figure}

For the $n=3$ case, a systematic presentation (excluding the Schur basis, and the forgotten one)  is given by the following table where  all expressions along a given line are equal.
\renewcommand{\arraystretch}{1.2}
$$  \begin{tabular}{|c|c|c|c|}
    \hline $m_3$  & $e_1^3-3\,e_1\,e_2+3\,e_3$ & $h_1^3-3\,h_1\,h_2+3\,h_3$ & $p_3$ \\[3pt]
    \hline $m_{21}$ & $e_1\,e_2-3\,e_3$ & $-2\,h_1^3+5\,h_1\,h_2-3\,h_3$ & $p_1\,p_2-p_3$  \\[3pt]
    \hline $m_{111}$  & $e_3$ & $h_1^3-2\,h_1\,h_2$ & $\frac16(p_1^3-3\,p_1\,p_2+2\,p_3)$ \\[3pt]
    \hline $m_{21}+3\,m_{111}$ & $e_{1}e_{2}$ & $h_1^3-h_1\,h_2$  & $\frac12(p_1^3-p_1\,p_2)$\\[3pt]
    \hline $m_3+3\,m_{21}+6\,m_{111}$  & $e_{1}^3$ & $h_1^3$ & $p_1^3$  \\[3pt]
    \hline  $m_3+m_{21}+m_{111}$ & $e_1^3-2\,e_1\,e_2+e_3$ & $h_{3}$ & $\frac16(p_1^3+3\,p_1\,p_2+2\,p_3)$ \\[3pt]
    \hline  $m_3+2\,m_{21}+3\,m_{111}$ & $e_1^3-e_1\,e_2$ & $h_{1}\,h_{2}$ & $\frac12(p_1^3+p_1\,p_2)$ \\[3pt]
    \hline  $m_3+m_{21}$ & $e_1^3-2\,e_1\,e_2$ & $-h_1^3+2\,h_1\,h_2$ & $p_{1}\,p_{2}$  \\[3pt]
     \hline
        \end{tabular}$$ \renewcommand{\arraystretch}{1}

\section{Jacobi-Trudi determinants}\label{sec_jacobi_trudi}
\define{Jacobi-Trudi's formula} (resp. its dual) gives an explicit expansion of the Schur functions in terms of complete homogeneous functions (resp. elementary). We have
\begin{equation}\label{jacobitrudi}
   s_{\mu}=\det(h_{{\mu_i}+j-i})_{1\leq i,j\leq n},\qquad \hbox{and resp.}\qquad s_{\mu'}=\det(e_{{\mu_i}+j-i})_{1\leq i,j\leq n},
 \end{equation}   
with $h_k=e_k=0$ whenever $k<0$.
This may be nicely deduced from the tableaux definition of Schur functions  using the technology of Lindstr\"om~\cite{lindstrom}  and Gessel-Viennot~\cite{gessel_viennot}. It is articulated around an interpretation of  Schur functions as configurations of non-intersecting north-east lattice paths in $\Z\times \Z$, which correspond to tableaux as illustrated in figure~\ref{ex_config}.

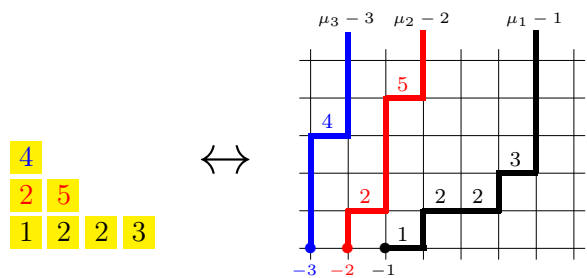
\begin{figure}[ht]
\begin{center}
\begin{picture}(0,0)(7,6)
   \Case{0,2}{\bleu{4}}
   \Case{0,1}{\rouge{2}} \Case{1,1}{\rouge{5}}
   \Case{0,0}{{1}} \Case{1,0}{{2}} \Case{2,0}{{2}} \Case{3,0}{{3}}
\end{picture}
\end{center}
\begin{center}
\begin{picture}(4,5)(-3,0)
\verythinlines
\multiput(-0.3,0)(0,1){6}{\line(1,0){7.6}}
\multiput(0,-0.3)(1,0){8}{\line(0,1){5.6}}
\put(-0.5,-0.7){\tiny \bleu{$-3$}}
\put(0.5,-0.7){\tiny \rouge{$-2$}}
\put(1.6,-0.7){\tiny {$-1$}}
\put(0.2,6){\tiny {$\mu_3-3$}}
\put(2.2,6){\tiny {$\mu_2-2$}}
\put(5.2,6){\tiny {$\mu_1-1$}}
\verythicklines
\put(0,0.02){\bleu{\circle*{.3}}}
\put(0,0.01){\bleu{\line(0,1){3.05}}}
\put(-0.03,3){\bleu{\line(1,0){1.1}}}
\put(0.3,3.2){\scriptsize{\bleu{4}}}
\put(1,2.95){\bleu{\line(0,1){2.8}}}
\put(0.98,0.02){\rouge{\circle*{.3}}}
\put(.99,-0.05){\rouge{\line(0,1){1.1}}}
\put(.95,1){\rouge{\line(1,0){1.1}}}
\put(1.3,1.2){\scriptsize{\rouge{2}}}
\put(2.01,0.95){\rouge{\line(0,1){3.1}}}
\put(2,4){\rouge{\line(1,0){1.05}}}
\put(2.3,4.2){\scriptsize{\rouge{5}}}
\put(3,4){\rouge{\line(0,1){1.8}}}
\put(1.98,0.02){{\circle*{.3}}}
\put(1.99,0){{\line(1,0){1.1}}}
\put(2.3,0.2){\scriptsize{{1}}}
\put(3,-0.05){{\line(0,1){1.1}}}
\put(2.95,1){{\line(1,0){2.1}}}
\put(3.3,1.2){\scriptsize{{2}}}
\put(5,0.95){{\line(0,1){1.1}}}
\put(4.95,2){{\line(1,0){1.1}}}
\put(4.3,1.2){\scriptsize{{2}}}
\put(5.3,2.2){\scriptsize{{3}}}
\put(6,2){{\line(0,1){3.8}}}
\put(-3,2){\huge$\leftrightarrow$}
\end{picture}
\end{center}
\caption{Non-crossing configuration associated to a semi-standard tableau.}
\label{ex_config}
\end{figure}

\section{Plethysm}\label{plethysme}

Plethysm was first defined by Littlewood (see \cite{macdonald}), as an operation between symmetric functions.
More generally, we consider here a symmetric function $f$ as an operator that transforms a rational expression $A$ into an expression $f[A]$.
Let $A$ and $B$ be rational fraction expressions in variables $\bm{x}=(x_1,x_2,\ldots)$, and consider symmetric functions $f$ and $g$ in $\Lambda$.
\define{Plethysm} is characterized by the following ``calculation'' rules, with $\alpha$, and $\beta$ denoting scalars:
\medskip
\hrule\medskip
\centerline{\underline{Rules for Plethysm}}
\bleu{\begin{enumerate}
\item  $(\alpha f+\beta g)[A]=\alpha\, f[A]+\beta\, g[A]$,
\item $(f\cdot g)[A]=f[A]\cdot g[A]$,
\item $p_k[\alpha A+\beta B]=\alpha\, p_k[A]+\beta\, p_k[B]$, for any power sum $p_k$,
\item $p_k[ A\cdot B]=p_k[A]\cdot p_k[B]$, 
\item $p_k[ A/B]=p_k[A]/p_k[B]$, 
\item $p_k[x]=x^k$, for $x$ a variable,
\item $p_k[1]=1$.
\end{enumerate}}
\hrule
Since any symmetric function $f$ may be expanded as a polynomial in the $p_k$, the first two rules clearly reduce the calculation of $f[A]$ to instances of $p_k[A]$. Then, $A$ being a rational expression in the $x_i$'s, the last five rules may be used to finish any needed calculation. \define{Words of caution} are necessary here. First, one needs to be clear about the \underline{distinction} between variables and scalars. Indeed, on variables $z$ one applies the rule $x\mapsto x^k$, whereas scalars are left invariant under plethysm. Second, as $f[-]$ is an operator it needs \underline{not commute} with other operators, such as specialization of the variables. 

We may consider Formula~\pref{schur_sum} as an instance of a plethystic calculation. Also, it follows from the rules that
\begin{align}
&\sum_{n\geq 0}h_n[A]\,z^n=\exp\Big({\sum_{k\geq 1}\displaystyle
          {p_k[A]\,}\frac{z^k}{k}}\Big),\qquad {\rm and} \label{h_plethysm} \\
& \sum_{n\geq 0}e_n[A]\,z^n=\exp\Big({\sum_{k\geq 1}\displaystyle
         (-1)^{k-1}  {p_k[A]\,}\frac{z^k}{k}}\Big),         
 \end{align}
 so that we may calculate many instances of plethysm in one stroke.
It is then helpful to remember that, considering $z$ as a variable, we have $f[z\,A]=z^d\,f[A]$, whenever $f$ is homogeneous of degree $d$.

Observe that a special case of plethysm may be considered an operation between symmetric functions (considering the case when $A$ lies in $\Lambda$). In this situation one also writes $f\circ g$ for $f[g]$. This is an associative operation, and the above rules imply that 
    \begin{equation}
        p_k\circ p_j=p_{kj}.
    \end{equation}
The ``unique'' homogeneous symmetric function of degree $1$, say $p_1$ ($=m_1=h_1=e_1=s_1$), acts as  identity for plethysm. In the context of plethysm, it is also very handy to consider variables sets as formal sums, writing 
    $$\bm{x}=x_1+x_2+x_3+\ldots+x_n+\ldots$$
 Indeed, we then have $f[\bm{x}]=f(x_1,x_2,\ldots)$, since $p_k[x_1+x_2+\ldots]=x_1^k+x_2^k+\ldots=p_k(\bm{x})$.  In the sequel, it will be convenient to use the \define{star} notation
    \begin{equation}
        f^*(\bm{x}):=f\Big[ \frac{\bm{x}}{1-q}\Big].
    \end{equation}

 \subsection*{\underline{Some examples}}
 \begin{enumerate}
 \item[$\bullet$] Calculating $f[x_1+x_2+\ldots+x_k]$  corresponds to the specialization of symmetric functions to a finite number of variables. We already know that $m_\mu[x_1+x_2+\ldots+x_k]$ vanishes when $\ell(\mu)>k$, and this is also the case for $s_\mu$.
 
  \item[$\bullet$] Calculating $f[-1]$ comes as our first ``surprise''.  Indeed, using~\pref{h_plethysm} we get:
	\begin{align*}
		&\sum_{n\geq 0}h_n[-1]\,z^n=\exp\Big({\sum_{k\geq 1}\displaystyle {p_k[-1]\,}\frac{z^k}{k}}\Big)
		                                           =\exp\Big(\sum_{k\geq 1}\displaystyle {-\frac{z^k}{k}}\Big)
		                                           =1-z\\
		&\sum_{n\geq 0}e_n[-1]\,z^n=\exp\Big({\sum_{k\geq 1}\displaystyle (-1)^{k-1} {p_k[-1]\,}\frac{z^k}{k}}\Big)
		                                           =\exp\Big(\sum_{k\geq 1}\displaystyle {\frac{(-z)^k}{k}}\Big)
		                                           =\frac{1}{1+z},
	\end{align*}
so that  
   $$h_n[-1]=\begin{cases}
     1, & \text{if}\ n=0,\\
     -1, & \text{if}\ n=1,\\
     0, & \text{otherwise},
\end{cases}\qquad {\rm and}\qquad
e_n[-1]=\begin{cases}
     1, & \text{if}\ n\ \text{is even},\\
     -1,  & \text{otherwise}.
\end{cases}$$
Using the Jacobi-Trudi identities, we may deduce that $s_\mu[-1]=0$ for all partitions of $n\geq 3$, except when $\mu=1^n$. Using the fact that $p_k[-\bm{x}]=-p_k(\bm{x})$, it may likewise be seen that
  \begin{equation}\label{moins_x}
     s_\mu[-\bm{x}]= (-1)^n s_{\mu'}(\bm{x}).
  \end{equation}
  \item[$\bullet$] For the plethysm $f[\bm{x}-\bm{y}]$, mixing Formula~\pref{h_somme} with~\pref{moins_x} we find that
  \begin{equation}\label{x_moins_y}
       h_n[\bm{x}-\bm{y}]= \sum_{k=0}^n (-1)^{n-k}h_k(\bm{x})\,e_{n-k}(\bm{y}).
    \end{equation}
More generally, for all partition $\mu$, we have
   \begin{equation}\label{Schur_moins_x}
       s_\mu[\bm{x}-\bm{y}]= \sum_{\nu\subseteq \mu} (-1)^{n-k} s_{\mu/\nu}(\bm{x})\,s_\nu(\bm{y}).
    \end{equation}

  \item[$\bullet$] For any partition $\mu$ of $n$, one may show that for a positive scalar $k$:
      \begin{align}\label{schur_eval_n}
         &s_\mu[k]=\prod_{(i,j)\in\mu}\frac{k+j-i}{\mathfrak{h}_{ij}}.
      \end{align} 
 In fact, this is the number of semi-standard tableaux of shape $\mu$ with values in $\{1,2,\ldots, k\}$.
One may also show that:
      \begin{align}
         &{\rm (a)}\ s_\mu\big[{\textstyle{\frac{1-q^k}{1-q}}}\big]=q^{n(\mu)} \prod_{(i,j)\in\mu}\frac{1-q^{k+j-i}}{1-q^{\mathfrak{h}_{ij}}}, \qquad{\rm and}\label{schur_eval_qn}\\
         &{\rm (b)}\ s_\mu^*(1)=s_\mu\big[{\textstyle{\frac{1}{1-q}}}\big]=\frac{q^{n(\mu)}}{\textstyle \prod_{(i,j)\in\mu}{1-q^{\mathfrak{h}_{ij}}}}.\label{specialization_standard}
      \end{align} 
  
  \item[$\bullet$] We have (see exercises)
\begin{equation}\label{plethysme_schur_1-u}
  \frac{s_\mu[1-u]}{1-u}=\begin{cases}(-u)^k, & \text{if}\ \mu=(n-k, 1^k),\\
                         0, & \text{otherwise},
                     \end{cases}
\end{equation}

\item[$\bullet$] We calculate that
  $$\sum_{n\geq 1} e_n[(1-q)(1-t)] z^n =  {\frac {z \left( 1-q \right)  \left( 1-t \right) }{ \left( 1+t\,z
 \right)  \left( 1+q\,z \right) }}
,$$
which implies that
  \begin{equation}\label{plet_en_1_t_1_q}
      \frac{(-1)^{n-1}}{(1-q)(1-t)}\, e_n[(1-q)(1-t)]  =  [n]_{q,t},
   \end{equation}
 where we use  the \define{$(q,t)$-integer} notation 
   $$\bleu{ [n]_{q,t}:=\frac{q^n-t^n}{q-t}=q^{n-1}+q^{n-2}t+\ldots+q\,t^{n-2}+t^{n-1}}.$$
Similarly, we may see that
   \begin{equation}\label{hstar_1}
     h_n^*(1)= h_n\Big[\frac{1}{1-q}\Big]=\prod_{k=1}^n \frac{1}{1-q^k}
   \end{equation}
 Observe that
 \begin{equation}\label{hstar_factorial}
       \frac{1}{h_n^*(1)}= (1-q)^n\,[n]_q!
 \end{equation} 

\item[$\bullet$] Starting with
 \begin{equation}\label{calcul_plet_hn}
    h_n^*({\bm{x}})=
     \sum_{\mu\vdash n} \frac{1}{z_\mu}\prod_{i= 1}^{\ell(\mu)} 
            \frac{p_{\mu_i}(\bm{x})}{1-q^{\mu_i}},
\end{equation}
together with~\pref{principal_spec}, and the Cauchy kernel formula, we may calculate that
 \begin{align}
       \frac{h_n^*(\bm{x})}{h_n^*(1)} 
       		&=\sum_{\mu\vdash n} \frac{s_\mu^*(1)}{h_n^*(1)}\, s_\mu(\bm{x}),\\
		&=\sum_{\mu\vdash n} K_{\mu,1^n}(q)\, s_\mu(\bm{x}).
 \end{align}
 \end{enumerate}
This last $\N[q]$-coefficient symmetric polynomial is an important instance of the ``Combinatorial Macdonald polynomials'', denoted by $H_\mu(q;\bm{x})$, that will arise in the next chapter. Thus, we have
\begin{equation}\label{defn_Hn}
    H_{n}(q;\bm{x}) :=  \frac{h_n^*(\bm{x})}{h_n^*(1)} .
\end{equation}
 Some examples are:
 \begin{align*}
&H_{1}= s_{1},\\
&H_{2}= s_{2} + q\,s_{11},\\
&H_{3}= s_{3}+ \left(q^{2} + q\right)s_{21} + q^{3}s_{111} ,\\
&H_{4}=s_{4} + \left(q^{3} + q^{2} + q\right)s_{31} + \left(q^{4} + q^{2}\right)s_{22} +  \left(q^{5} + q^{4} + q^{3}\right)s_{211} +  q^{6}s_{1111},\\
&H_{5}=s_{5} + \left(q^{4} + q^{3} + q^{2} + q\right)s_{41} + \left(q^{6} + q^{5} + q^{4} + q^{3} + q^{2}\right)s_{32},\\
&\qquad\qquad + \left(q^{6} + q^{5} + q^{4} + q^{3} + q^{2}\right)s_{32}  + \left(q^{7} + q^{6} + 2 q^{5} + q^{4} + q^{3}\right)s_{311},\\
&\qquad\qquad +  \left(q^{8} + q^{7} + q^{6} + q^{5} + q^{4}\right)s_{221}  + \left(q^{9} + q^{8} + q^{7} + q^{6}\right)s_{2111} +    q^{10}s_{11111}. 
\end{align*}
 
 \section{Some ties with geometric complexity theory}
As explained in ~\cite{Mulmuley}, fundamental  \define{$\bm{P} \not= \bm{NP}$} questions of complexity theory may be ``reformulated'' (this has to be clarified) in the context of calculations of the Schur function expansions of either $s_\mu s_\nu$, $s_\mu * s_\nu$, or of $s_\mu\circ s_\nu$. The aim is to prove typical complexity statements about explicit calculations of the Littlewood-Richardson coefficients $c_{\mu\nu}^\lambda$ in the expansion~\pref{LR_coeff}:
\begin{displaymath}
    s_\mu s_\nu= \sum_{\lambda\vdash (m+n)} c_{\mu\nu}^{\lambda} s_\lambda,
 \end{displaymath}
 or similar expansions for Kronecker product, or plethysm.
Questions regarding these include showing that $c_{\mu\nu}^{\lambda}\not=0$, for  $\mu$, $\nu$, and $\lambda$ satisfying some constraints.  The Littlewood-Richardson rule does give a procedure, but it is hard to apply when partitions involved become large. The following is chief among results in this context.

\begin{theorem}[Knutson-Tao, De Loera-McAllister] The problem of deciding the nonvanishing of $c_{\mu\nu}^{\lambda}$ is in $\bm{P}$, this is to say that it can be solved in time that is polynomial in the bit-lengths of $\mu$, $\nu$ and $\lambda$.
\end{theorem} 
      
There are several open conjecture in this area, which is currently very active. An overview may be found in \cite{burgisser}, and more complete presentations in~\cite{Landsberg,Landsberg2}. You may find there an explanation of how the Foulke's conjecture plays a role in this story.

 \section{Counting with symmetric functions}
There are several classical formulas and identities from enumerative combinatorics that may be obtained using symmetric functions.
Indeed, we can get interesting integers (or $q$-integers) out of a symmetric function $f$ (typically homogeneous), in the following way.
First, one may calculate some scalar product $\langle f, g\rangle$ choosing a ``special'' $g$, such as
   $$ \langle f, e_n\rangle,\qquad \langle f, h_n\rangle,\qquad \langle f, h_k e_{n-k}\rangle,\qquad{\rm or}\qquad \langle f, s_1^k\rangle.$$
 For example,
    $$\langle h_1^n,h_1^n\rangle = n!, \qquad \langle h_1^n,h_ke_{n-k}\rangle= \binom{n}{k},$$
One may also evaluate $f$ plethysticaly at an integer or $q$-integer:
   $$f[n],\qquad{\rm or}\qquad f[1+q+\ldots+q^{n-1}].$$
For instance, formulas~\pref{schur_eval_n} and~\pref{schur_eval_qn} give:
\begin{align}
   &e_k[n]=\binom{n}{k}, &&h_k[n]=\binom{n+k-1}{k},\label{formule_binome}\\
   &e_k[1+q+\ldots+q^{n-1}]=q^{\binom{k}{2}}\,\qbinom{n}{k},&&h_k[1+q+\ldots+q^{n-1}]=\qbinom{n+k-1}{k}
\end{align}
It follows that we get combinatorial identities out of any symmetric function identity, and this applies to series. To illustrate,
consider what happens to the series in~\pref{defn_H} when one assumes that there are $n$ variables $x_i$, which are all set to equal $1$.
We get, using~\pref{formule_binome},
   \begin{align}
        \sum_{k\geq 0}\binom{n}{k} z^k =(1+z)^k, && {\rm and}\qquad  \sum_{k\geq 0}\binom{n+k-1}{k} z^k =\left(\frac{1}{1-z}\right)^n.
    \end{align}
similarly, with $x_i=q^{i-1}$, we find that
    \begin{align}
        \sum_{k\geq 0}q^{\binom{k}{2}}\,\qbinom{n}{k} z^k =\prod_{i=1}^n (1+q^{i-1}z), && {\rm and}\qquad  
        \sum_{k\geq 0}\qbinom{n+k-1}{k}  z^k =\prod_{i=1}^n \left(\frac{1}{1-q^{i-1} z}\right).
    \end{align}
Exploiting remarks of our first chapter, it may be seen that the last formula counts partitions with some restriction (which is it?), taking into account the number of parts that they have.

The \define{theory of species} (for which the typical reference is~\cite{BLL}) furnishes a very general systematic approach to the part of enumerative combinatorics that considers constructions on finite sets, and the corresponding Pólya theory. Recall that Pólya theory gives general tools for the enumeration of \define{isomorphism type of structures} (a.k.a. structures on unlabelled objects). See chapter 6 of \href{http://bergeron.math.uqam.ca/wp-content/uploads/2015/04/Comb_alg1.pdf}{Bergeron/Combinatoire\_Algebrique} (in french), or \href{https://en.wikipedia.org/wiki/Combinatorial_species}{wiki/Combinatorial\_species} for simple introductions. 

To each species of combinatorial objects there corresponds a symmetric function series, known as its Pólya-Joyal cycle index series. Via a combinatorial calculus on species, which allows the construction of more complex species out of simpler ones, one gets symmetric function series identities (involving, sums, product, plethysm, and Kronecker product). This makes it possible to obtain both the enumeration of labelled and unlabelled structures, simply from specializations of formulas entirely coined in terms of symmetric functions series, in which the ``basic'' variables are power-sums.

Examples of typical species include: graphs, connected graphs, oriented graphs, permutations, cyclic permutations, derangements, orders, lists, $k$-subsets, set partitions, trees of various kinds, etc. In the algebra of (combinatorial constructions on) species, one may solve equations, consider Lagrange inversion, and interpret combinatorially many other classical calculations. For sure, all of this corresponds to interesting calculations on symmetric functions series.

\section{Exercises and problems}\label{Probabililty_Schur}
\begin{exer}\rm 
\begin{enumerate}
\item[(a)] Adapt RSK to show that
	\begin{equation}
		(x_1+x_2+\ldots + x_k)^n = \sum_{\mu\vdash n} f^\mu\,s_\mu(\bm{x}),
	\end{equation}
	with $f^\mu$ the number of standard tableaux of shape $\mu$.
\item[(b)] Find the monomial functions expansion of the above formula (as a linear combination of the $m_\mu$), and check that it is essentially Newton's binomial formula.
\end{enumerate}
\end{exer}

\begin{exer}\rm For any partition $\mu$, let us set
\begin{equation}\label{defn_pi_mu}
    \pi_\mu(\bm{x})=\pi_{\mu_1}(\bm{x})\pi_{\mu_2}(\bm{x})\cdots \pi_{\mu_k}(\bm{x}),
 \end{equation}
  where $\pi_n(\bm{x})$ is the following alternating sum of Schur functions indexed by hooks:
\begin{equation}\label{deff_pi_n}
     \pi_n(q,t;\bm{x}):=s_{1^n}(\bm{x})+\ldots +\big({-1}/{q\,t}\big)^{k-1}\!s_{(k,1^{n-k})}(\bm{x}) +\ldots+\big({-1}/{q\,t}\big)^{n-1}\!s_{n}(\bm{x}).
  \end{equation}
\begin{enumerate}
\item[(a)] Exploit \pref{lienHP} to check that $P(z)=-E'(-z)H(z)$, and use the Pieri rule to deduce from this the Schur expansion of $p_n$.
\item[(b)] Check that $\pi_n$ specializes to $(-1)^{n-1}p_n$ whenever $qt=1$.
\item[(c)] Conclude that the set $\{ \pi_\mu\}_{\mu\vdash d}$ constitutes a basis of $\Lambda_d$.
\end{enumerate}
\end{exer}

\begin{exer}\rm The aim here is to prove that  
    \begin{equation}\label{prob_p1p2}
          \sum_{k=0}^{\lfloor n/2\rfloor} p_1^{n-2\,k} p_2^k\quad\hbox{ is Schur-positive},
     \end{equation}
for all $n\geq1$. 
Recall that $p_2=s_2-s_{11}$. Consider the series 
    $${\frac {1}{ \left( 1-p_1\,z \right)  \left( 1-p_2\,{
z}^{2} \right) }}=1+p_1\,z+ \left( {p_1}^{2}+p_2 \right) {z}^{2}+ \left( {p_1}^{3}+p_1\,p_2
 \right) {z}^{3}+ \left( {p_1}^{4}+p_2\,{p_1}^{2}+{
p_2}^{2} \right) {z}^{4}+\ldots 
$$
\begin{enumerate}
\item[(a)] Show that
  \begin{align*}
       & {\frac {1}{ \left( 1-p_1\,z \right)  \left( 1-p_2\,{
z}^{2} \right) }} =\frac{1+s_1z}{2\,{s_{{11}}} {z}^{2}}
 \left( \frac{1}{1- \left( s_{{2}}+s_{{11}} \right) {z}^{2}}- 
         \frac{1}{1- \left( s_{{2}}-s_{{11}} \right) {z}^{2}} 
   \right).
  \end{align*}
  \item[(b)] Show that the above expression may be written without division by $2\,s_{11}z^2$.
  \item[(c)] Conclude that \pref{prob_p1p2} is Schur positive.
\end{enumerate}
\end{exer}

\begin{exer}\rm 
\begin{enumerate}
  \item[(a)] Assuming the Pieri rules, prove that $h$-positivity (or $e$-positivity) implies Schur positivity, and that the product of a $h$-positive symmetric function with a Schur positive one is also Schur positive. 
   \item[(b)] Show that one may recursively calculate the Schur functions using the Pieri rule (or dual Pieri rule). Compare your approach to the Jacobi-Trudi formula.
   \item[(c)] Exploit the above to prove that $\omega s_\mu=s_{\mu'}$.
 \end{enumerate}
\end{exer}

 \begin{exer}\rm 
\begin{enumerate} 
  \item[(a)] With scalars in $\mathbb{R}$, consider the set of symmetric functions of the form
                    $$\sum_{\mu\vdash n} c_\mu\,m_\mu,\qquad{\rm with}\qquad 0\leq c_\mu\leq 1,$$
             such that $\sum_{\mu\vdash n}c_\mu=1$. One may thus think that this is a simplex in $\mathbb{R}^{p(n)}$, where $p(n)$ is the number of partitions of $n$.
             Show that the Schur positive symmetric functions form a convex subset of this simplex. 
    \item[(b)] Defining the probability of Schur positivity as the quotient of the volume of this convex subset, by the volume of the simplex, prove theorem~\ref{thm_BPR}.
  \item[(c)] {\bf (explore)} Find a formula for the probability of being $e$-positive among Schur positive symmetric functions of degree $n$. Likewise for $h$-positivity.
\end{enumerate}
\end{exer}

\begin{exer}\rm   
\begin{enumerate}
\item[(a)] Show that the determinant of the \define{Hankel} matrix\footnote{See \cite{kratt1, kratt2, kratt} for many interesting related identities.}  
     $$\det( h_{i+j+k})_{0\leq i,j\leq n-1}$$ 
    is equal (up to a sign)  to a single Schur function (What is the sign? What is the shape?). For example, for $n=3$ and $k=1$, we have the Hankel matrix
     $$\Big(h_{i+j+1}\Big)_{0\leq i,j\leq 2}=\begin{pmatrix} h_1 & h_2 & h_3\\ h_2 & h_3 & h_4\\ h_3 & h_4 & h_5\end{pmatrix}.$$
  \item[(b)] Same question for $\det( e_{i+j+k})_{0\leq i,j\leq n-1}$.
 \end{enumerate}
\end{exer}

 \begin{exer}\rm 
\begin{enumerate}
\item[(a)] Extend formula (\ref{s_en_m}) to the context of skew Schur functions.
Prove that $s_{\lambda/\mu}(\bm{x})$ is invariant under the exchange of the variables $x_i$ and $x_{i+1}$, by constructing for each $\lambda/\mu$-shape semi-standard tableau $\tau$  a new semi-standard tableau of same shape but in which the respective number of occurrences of $i$ and $i+1$ are exchanged.
Conclude that $s_{\lambda/\mu}(\bm{x})$ is symmetric.
\item[(b)] There is an analog of formula~\pref{schur_sum} for skew Schur functions. Can you figure out what it is?
\item[(c)] There is an analog of the Jacobi-Trudi formula for skew Schur functions. Can you figure out what it is?
\end{enumerate}
\end{exer}

\begin{exer}\rm 
\begin{enumerate}
\item[(a)] Use the (dual) Pieri formula, and the fact that the Schur functions are orthonormal, to show that for all partition $\mu$ of $n$, we have $\langle e_\mu,e_n\rangle =1$.
\item[(b)] Recall that we have set $\langle p_\lambda,p_\mu\rangle =0$, whenever $\lambda\not=\mu$. Hence $\langle f,g\rangle =0$ whenever $f$ and $g$ are both homogeneous, but of different degree.  Conclude that,
     \begin{equation}
         \langle f,E\rangle =\sum_\mu c_\mu\qquad {\rm if}\qquad f(\bm{x})=\sum_\mu c_\mu\, e_\mu(\bm{x}).
     \end{equation}
where, as before, $E:=\sum_n e_n$.
\end{enumerate}
\end{exer}

\begin{exer}\rm \begin{enumerate}
\item[(a)] By a direct calculation, show that 
\begin{equation}\label{dualite_homogene}
 \prod_{i,j\geq 1}\frac{1}{1-{x_i}\,{y_j}}=\sum_\lambda h_\lambda({\bm{x}})m_\lambda({\bm{y}}). 
\end{equation}
\item[(b)] Show that the bases $\{u_\lambda\}_\lambda$ and $\{v_\mu\}_\mu$ are dual, if and only if the bases $\{\omega\,u_\lambda\}_\lambda$ and $\{\omega\,v_\mu\}_\mu$ are also dual.
\item[(c)] Conclude that the set of $f_\mu:=\omega\, m_\mu$, for $\mu\vdash n$, (aka the \define{forgotten} basis) constitutes a basis which is dual to the $e_\lambda$.
\item[(d)] The \define{dual Cauchy Kernel} $\Omega'(\bm{xy})$ is 
       \begin{equation}
           \bleu{\Omega'(\bm{xy}):=\prod_{i,j\geq 1}(1+{x_i}\,{y_j})}
       \end{equation}
  Check that its degree $n$ component is
       \begin{align}
           e_n(\bm{xy}) & = \sum_{\mu\vdash n} (-1)^{n-\ell(\mu)} p_\mu(\bm{x})\frac{p_\mu(\bm{y})}{z_\mu}\\
                                 & = \sum_{\mu\vdash n} s_\mu'(\bm{x})\,s_\mu(\bm{y})\label{dual_schur}\\
                                 & = \sum_{\mu\vdash n} e_\mu(\bm{x})\,m_\mu(\bm{y})\\
                                 & = \sum_{\mu\vdash n} h_\mu(\bm{x})\,f_\mu(\bm{y})\label{dualite_elementaire}.
      \end{align}
\end{enumerate}
\end{exer}

\begin{exer}\rm
Show the equivalence of (\ref{scal_dual}) and (\ref{noyeau_dualite}) by expanding $p_\lambda$ and $p_\mu/z_\mu$ in the given basis $u_\lambda$ and $v_\mu$: to get $p_\lambda=\sum_\rho a_{{\lambda}{\rho}}\, v_\rho$, and $p_\mu/z_\mu=\sum_\nu b_{{\mu}{\nu}}\, v_\nu$.
and that we must have $ \delta_{\lambda,\mu}=\langle p_\lambda,p_\mu/z_\mu\rangle=\sum_{\rho,\nu} a_{\rho\lambda}\,\langle u_\lambda,v_\nu\rangle\,b_{\nu\mu}$.  In matrix form, this becomes $A\,Z\,B^{\rm Tr}={\rm Id}$,
where $A:=(a_{\lambda\rho})$, $B:=(b_{{\mu}{\nu}})$,  and $Z:=(\langle u_\lambda,v_\mu\rangle)$.
Check that statement (\ref{scal_dual}) is equivalent to $Z={\rm Id}$, and thus it is also equivalent to $ A\,B^{\rm Tr}={\rm Id}$. Using the fact that
  $$\sum_\lambda u_\lambda(\bm{x})\,v_\lambda(\bm{y})= \sum_\lambda p_\lambda(\bm{x})\,\frac{p_\lambda(\bm{y})}{ z_\lambda}.$$
and that the $p_\lambda(\bm{x})\,p_\lambda(\bm{y})/ z_\lambda$'s are linearly independent, 
conclude that  $\sum_\lambda a_{\lambda\rho}\,b_{\lambda\nu}=\delta_{\rho,\nu}$,
is equivalent to (\ref{noyeau_dualite}). 
\end{exer}

\begin{exer}\rm \begin{enumerate}
\item[(a)] Using~\pref{hstar_factorial} and \pref{specialization_standard}, show that for all $\mu$, the expression $s_\mu^*(1)/h_n^*(1)$ is a polynomial in $q$.
\item[(b)] Give an expression for  $s_\mu^*(1)/h_n^*(1)$ in terms of $[n]_q!$, and the various $[\mathfrak{h}_{ij}]_q$ for $(i,j)\in\mu$.
\item[(c)] Find a nice formula for $e_n^*(1)$.
\end{enumerate}
\end{exer}
 
 \begin{exer}\rm \begin{enumerate}
 \item[(a)] Our aim here is to prove~\pref{plethysme_schur_1-u}. 
 To this end, find the value of each $h_n[1-u]$ by calculating the series
 \begin{align}
 \sum_{n\geq 0}h_n[{1-u}]\,z^n
\end{align}
using~\pref{h_plethysm}, and $p_k[1-u]=1-u^k$.
Conclude using Jacobi-Trudi.
\item[(b)] Show that
 \begin{align}
\frac{h_n(\bm{x})}{1-q}= \sum_{k=0}^{n-1} (-q)^k\, s_{(n-k,1^k)}\Big[\frac{\bm{x}}{1-q}\Big] .
\end{align}
To this end, calculate
  $$ \frac{h_n(\bm{x})}{1-q}=\frac{1}{1-q}\,h_n\Big[ (1-q)\,\frac{\bm{x}}{1-q}\Big],$$
using the Cauchy kernel formula with $\bm{x}=1-q$ and $\bm{y}=\frac{1}{1-q}$, and conclude using~\pref{plethysme_schur_1-u}.
\end{enumerate}
\end{exer}

\begin{exer}\rm 
\begin{enumerate}
\item[(a)]  Show that plethysm of symmetric functions is associative.
\item[(b)] Show that, for any $a$ and $b$,
    \begin{equation}
        (h_a\circ h_b)[1+q]=\qbinom{a+b}{b}.
     \end{equation}
\item[(c)] Prove \define{Hermite reciprocity}, i.e.:
    \begin{equation}
        (h_a\circ h_b)(x,y)= (h_b\circ h_a)(x,y).
     \end{equation}
\item[(d)] Show that, for all $n$
   \begin{equation}
        (h_2\circ h_n)=\sum_{k=0}^{\lfloor n/2\rfloor} s_{(2n-2k,2k)}.
     \end{equation}  
  \item[(e)] Show that,  for all $n$
   \begin{equation}
        (h_n\circ h_2)=\sum_{\audessus{\mu\vdash 2n}{\mu\ {\rm has\ even\ parts}}}^{\lfloor n/2\rfloor} s_{\mu}.
     \end{equation}  
  \item[(f)] {\bf (open)} \define{Foulkes conjecture}, for any $a<b$,
     \begin{equation}
        (h_b\circ h_a)-(h_a\circ h_b) \quad \hbox{is Schur postive}.
     \end{equation}
   \item[(g)] {\bf (open)} \define{Generalized Foulkes conjecture} (due to Vessenes), for any $a\leq b,c\leq d$,
     \begin{equation}\label{gen_Foulkes}
        (h_b\circ h_c)-(h_a\circ h_d) \quad \hbox{is Schur postive}.
     \end{equation}
      \item[(g)] Calculate $s_\mu[1+q]$, and conclude that~\pref{gen_Foulkes} implies~\pref{q_foulkes}.
      \item[(h)] Let $n=ab$. Show that $\langle h_a\circ h_b,p_1^{ab}\rangle $ is the number of set partitions of a $n$-set into $a$ blocks of size $b$.
\end{enumerate}
\end{exer}

\begin{exer}\rm In this exercise we consider plethysm for series like $H(z;\bm{x})=H(z)$ (see~\pref{defn_H}), here underlining the role of $\bm{x}$.  We also write $H$ for $H(1)$, and observe that $H(\bm{x};z)$ may be obtained as the plethysm $H[z\,\bm{x}]$, thinking that $z\,\bm{x}=z\cdot(x_1+x_2+\ldots)$. 
\begin{enumerate}
  \item[(a)] Considering $x$ as a (single) variable, so that $p_k[x]=x^k$, find closed form formulas (not series) for $H[x]$, $P[x]$, and $E[x]$.
  \item[(b)] Let $\bm{x}$ and $\bm{y}$ be two sets of variables. Show that $H[\bm{x}+\bm{y}] = H(\bm{x})\,H(\bm{y})$.
  \item[(c)] Considering  $t$ as a constant,  show that $H[t\bm{x}]=H(\bm{x})^t$.
    \item[(d)] Considering $x$ as a (single) variable, show that 
       $$(H\circ (H-1))[x]=\prod_{j\geq 1} \frac{1}{1-x^j}.$$
    \item[(d)] Consider any symmetric function series of the form 
         $$F( \bm{x} )= p_1(\bm{x}) +\sum_{d\geq 2} f_d(\bm{x}),\qquad {\rm with}\qquad f_d\in \Lambda_d.$$
      Can you prove that $F$ affords an inverse for plethysm? 
 \end{enumerate} 
 \end{exer}
 
\end{chapter}


\begin{chapter}{\bleu{Macdonald symmetric functions and operators}}

For more background and details, see FB.~\cite{bergeron} or Haglund~\cite{haglund_livre}. 
The \href{https://more-sagemath-tutorials.readthedocs.io/en/latest/tutorial-symmetric-functions.html}{Sage Symmetric Function Tutorial} includes examples about Macdonald symmetric functions, and related operators.

\section{Macdonald symmetric functions}
The notions that we plan to explore are closely linked to the theory of Macdonald symmetric functions, and operators for which they are joint eigenfunctions. 

\subsection{Macdonald's original definition}\label{original_macdonald}

Our context here is the ring $\Lambda_{qt}$ of symmetric functions with coefficients in the fraction field $\C(q,t)$ in the parameters, $q$ and $t$, equipped with the \define{scalar product}\index{scalar product}
\begin{equation}\label{scal_mac}
  \bleu{ {\langle p_\lambda,p_\mu \rangle}_{q,t}:=\begin{cases}z_{\lambda}(q,t) & \text{if}\ \lambda=\mu,\\[10pt]
                               0 & \text{otherwise},\end{cases}}
\qquad
{\rm where}
\qquad 
    \bleu{\textstyle z_{\lambda}(q,t):=z_{\lambda}\cdot p_\lambda\left[\frac{1-q}{1-t}\right]}.
\end{equation}
In his original paper of 1988, Macdonald establishes  the existence and uniqueness of symmetric functions (polynomials) $P_\mu={P_\mu(q,t;\bm{x})}$\index{Macdonald symmetric functions} such that
\begin{itemize}\itemsep=3pt\itemindent=20pt
  \item[(1)] $P_\mu=m_\mu+\displaystyle{\sum_{\lambda \prec \mu} \gamma_{\mu \lambda}(q,t)}\, m_\lambda$,\quad
with coefficients $\gamma_{\mu \lambda}(q,t)$ in $\C(q,t)$; and
  \item[(2)] $\langle P_\lambda,P_\mu \rangle_{q,t}=0$, whenever $\lambda \not= \mu$. 
 \end{itemize}
Since at $q=t$, the scalar product (\ref{scal_mac}) specializes to the usual scalar product for symmetric functions, we get  
 $P_{\mu}(q,q;\bm{x})=s_{\mu}(\bm{x})$.
The Hall-Littlewood symmetric functions are obtained by setting $q=0$, and the Jack symmetric functions by setting $q=t^\alpha$ ($\alpha\in\R$ and $\alpha>0$) and letting $t\rightarrow 1$. Zonal symmetric functions are just the special case when $\alpha=2$; and we have
$$P_{\mu}(q,1;\bm{x})=m_{\mu},\quad \mathrm{and}\quad P_{\mu}(1,t;\bm{x})=e_{\mu'}.$$

\subsection{The ``combinatorial'' renormalization}
We are interested here in the \define{combinatorial} renormalization of the Macdonald symmetric functions, denoted by $H_\mu$, as introduced by Garsia.  One sets
\begin{equation}
    H_\mu(q,t;\bm{x}):=P_{\mu}\left[ q,t^{-1};\frac{\bm{x}}{1-t} \right]
                                  \prod_{c\, \in\, \mu} (q^{a(c)}-t^{\ell(c)+1}),
 \end{equation}
where $a(c)=a_\mu(c)$ and $\ell(c)=\ell_\mu(c)$ denote respectively the arm length and the leg length of the cell $c$ of $\mu$, see section~\ref{Young}.  After a long process, it has been shown that the Schur expansion of the functions $H_\mu(q,t;\bm{x})$ has coefficients in $\N[q,t]$. We have already encountered $H_n$ in~\pref{defn_Hn}. There is an interesting combinatorial description of the $H_\mu$ given in~\cite{haglund}. Small explicit values are as follows: 
have 
\begin{align*}
   &H_2=s_2+q\,s_{11},\\
   &H_{11}=s_2+t\,s_1;\\[4pt]
   &H_3=s_{{3}}+ \left( {q}^{2}+q \right) s_{{21}}+{q}^{3}s_{{111}},\\
   &H_{21}=s_{{3}}+ \left( q+t \right) s_{{21}}+q\,t\,s_{{111}},\\
   &H_{111}=s_{{3}}+ \left( {t}^{2}+t \right) s_{{21}}+{t}^{3}s_{{111}};\\[4pt]
   &H_4=s_4+(q^3+q^2+q)\,s_{31}+(q^4+q^2)\,s_{22}+(q^5+q^4+q^3)\,s_{211}+q^6\,s_{1111},\\
   &H_{31}=s_4+(q^2+q+t)\,s_{31}+(q^2+q\,t)\,s_{22}+(q^3+q^2\,t+q\,t)\,s_{211}+q^3\,t\, s_{1111},\\
   &H_{22}= s_4+(q\,t+q+t)\,s_{31}+(q^2+t^2)\,s_{22}+(q^2\,t+q\,t^2+q\,t)\,s_{211}+q^2\,t^2\,s_{1111},\\
   &H_{211}=s_4+(q+t+t^2)\,s_{31}+(q\,t+t^2)\,s_{22}+(q\,t+q\,t^2+t^3)\,s_{211}+q\,t^3\, s_{1111},\\
   &H_{1111}=s_4+(t+t^2+t^3)\,s_{31}+(t^2+t^4)\,s_{22}+(t^3+t^4+t^5)\,s_{211}+t^6\,s_{1111}.
\end{align*}
The coefficients $K_{\lambda,\mu}(q,t)$ of the expansion $H_\mu(q,t;\bm{x})=\sum_{\lambda\vdash n}
       K_{\lambda,\mu}(q,t)\, s_\lambda(\bm{x})$,
are said to be \define{$q,t$-Kostka polynomials}. We have the specializations 
   \begin{align}  
       &K_{\lambda,\mu}(1, 1)=f^\mu, && K_{\lambda,\mu}(0,1)=K_{\lambda,\mu}\\
       &H_\mu(0,0;\bm{x})=s_n,  && H_\mu(0,1;\bm{x})=h_\mu,  && H_\mu(1,1;\bm{x})=s_1^n,
  \end{align}
The next properties that stand out are the symmetries 
  \begin{align}
&H_{\mu'}(q,t;\bm{x})=H_\mu(t,q;\bm{x}),&& H_\mu(q,t;\bm{x})=q^{n(\mu')}t^{n(\mu)}\,  \omega\!\left(H_{\mu}(1/q,1/t;\bm{x})\right),\\
&K_{\lambda,\mu}(q,t)=K_{\lambda,\mu'}(t,q),&&K_{\lambda,\mu}(q,t)=q^{n(\mu')}t^{n(\mu)}\, K_{\lambda',\mu}(1/q,1/t)\label{K_omega}.
\end{align}
These symmetries are made apparent in the $q,t$-Kostka matrix given here for $n=4$:
$$\begin{pmatrix}
 1&{q}^{3}+{q}^2+q&{q}^{4}+{q}^2&{q}^{5}+{q}^{4}+{q}^{3}&{q}^{6}\\
1&{q}^2+q+t&{q}^2+qt&{q}^{3}+{q}^2t+qt&{q}^{3}t\\
1&qt+q+t&{q}^2+{t}^2&{q}^2t+q{t}^2+qt&{t}^2{q}^2\\
1&q+{t}^2+t&qt+{t}^2&q{t}^2+qt+{t}^{3}&q{t}^{3}\\ 
1&{t}^{3}+{t}^2+t&{t}^{4}+{t}^2&{t}^{5}+{t}^{4}+{t}^{3}&{t}^{6}
\end{pmatrix}.$$
Their specialization $t=1$, are multiplicative (see exercises), {\it i.e.}: 
\begin{equation}\label{special_t1}
  H_\mu(q,1;\bm{x})=H_{\mu_1}(q,1;\bm{x})\cdots H_{\mu_k}(q,1;\bm{x}),
\end{equation}
and thus characterized by the simple functions $H_n(q;\bm{x})=H_n(q,1;\bm{x})=h_n^*(\bm{x})/h_n^*(1)$ (see formula~\pref{defn_Hn}).
A second interesting specialization is at $t=1/q$. We then get the following Schur expansion
  \begin{equation}\label{specialization_qinv}
 H_{\mu}(q,1/q;\bm{x})=q^{-n(\mu)}
                           \prod_{c\, \in\, \mu} (1-q^{a(c)+\ell(c)+1})\, s_{\mu}^*(\bm{x}).
  \end{equation}
Other formulas are obtained by specializing $\bm{x}$ (using plethystic substitutions). For instance, we have
  \begin{equation}\label{1_u}
 H_\mu[q,t;1-u]=\prod_{(i,j)\in\mu} (1-q^it^j\,u),
\end{equation}
which implies that 
    \begin{equation}\label{K_hook}
        \langle H_\mu,s_{(n-k,1^k)}\rangle  = e_k[B_\mu-1]\qquad{\rm where}\qquad B_\mu=B_\mu(q,t):=\sum_{(i,j)\in \mu} q^it^j.
    \end{equation}
   Since the cell $(0,0)$ always belongs to $\mu$, the polynomial $B_\mu(q,t)-1$ lies in $\N[q,t]$.
In particular,
\begin{equation}\label{Hn_binom}
    \langle H_n,s_{(n-k,1^k)}\rangle = q^{\binom{k}{2}}\qbinom{n-1}{k}.
\end{equation} 
 An analog of the Cauchy-formula for the $H_\mu$ is:  
\begin{equation}\label{H_mu_cauchy}
   e_n\left[\frac{\bm{x}\bm{y}}{ (1-q)(1-t)}\right]= \sum_{\mu\vdash n}\frac{H_\mu(\bm{x})\, H_\mu(\bm{y})}{w_\mu(q,t)}
\end{equation}
where  
  $$w_\mu(q,t):=\prod_{c\in\mu}(q^{a(c)}-t^{\ell(c)+1})(t^{\ell(c)}-q^{a(c)+1)}.$$
In particular,
\begin{equation}\label{H_mu_cauchy_x}
   e_n\left[\frac{\bm{x}}{ (1-q)(1-t)}\right]= \sum_{\mu\vdash n}\frac{H_\mu(\bm{x})}{w_\mu(q,t)}
\end{equation}
We may also characterize the $H_\mu$ as being the unique solution of the system of equations:
  \begin{equation}
 \begin{array}{llll}
    \mathrm{(i)}\ {\displaystyle  \langle s_\lambda(\bm{x}), H_\mu[q,t;(1-q)\,\bm{x}]\rangle=0},\qquad{\rm if}\qquad  \lambda\not\succeq\mu,\\[6pt]
    \mathrm{(ii)}\ {\displaystyle   \langle s_\lambda(\bm{x}), H_\mu[q,t;(1-t)\,\bm{x}]\rangle=0},\qquad{\rm if}\qquad   \lambda\not\succeq\mu',\ \mathrm{and}\\ [6pt]
    \mathrm{(iii)}\ {\displaystyle   \langle s_n(\bm{x}), H_\mu(q,t;\bm{x})\rangle=1},
      \end{array}
 \end{equation} 
It is also useful to have the expansion
\begin{align}
 &e_n(\bm{x})=(1-q)(1-t) \sum_{\mu} \frac{\Pi_\mu B_\mu}{\mathfrak{h}_\mu \mathfrak{h}_{\mu}'}\, H_{\mu}(\bm{x}),\qquad {\rm with}\qquad \Pi_\mu:=\prod_{\audessus{(i,j)\in\mu,}{(i,j)\not=(0,0)}}(1-q^{i}t^{j});
 \end{align}
and
\begin{align}
 &(-1)^{n-1}p_n(\bm{x})=(1-t^n)(1-q^n)\sum_{\mu}
    \frac{\Pi_\mu}{\mathfrak{h}_{\mu}\mathfrak{h}_{\mu}'}\, H_{\mu}(\bm{x}).
  \end{align}

 \subsection{Macdonald scalar product}
The \define{$\rouge{\star}$-scalar product} is defined on the power-sum basis by
    \begin{displaymath}
      \bleu{\langle p_\mu,p_\lambda\rangle_{\rouge{\star}} }:= \begin{cases}
     \bleu{Z_\mu(q,t)}  & \text{if}\ \bleu{\mu=\lambda}, \\
     \bleu{0} & \text{otherwise},
\end{cases} 
    \end{displaymath}
  where $ Z_\mu(q,t)$ is defined by
  \begin{displaymath}
     \bleu{Z_\mu(q,t):=(-1)^{|\mu|-\ell(\mu)}\,p_\mu[(1-q)\,(1-t)]\,z_\mu}.
  \end{displaymath}
Observe that our two scalar products are linked by the relation (to be checked on power-sums)
\begin{displaymath}
 {\langle f,g\rangle = \langle f,\omega\, g^{\rouge{\star}}\rangle_{\rouge{\star}} }, 
 \end{displaymath}
where 
    \begin{displaymath} \bleu{g^{\rouge{\star}}(\bm{x}):=g\!\left[\frac{\bm{x}}{(1-t)(1-q)}\right]}.\end{displaymath}
 The Macdonald polynomials form an orthogonal family for this scalar product. More precisely, we have
       \begin{displaymath} {\langle H_\mu,H_\lambda\rangle_{\rouge{\star}} =0 }, \qquad {\rm whenever}\qquad {\lambda\not=\mu}.
        \end{displaymath}

 \section{Macdonald eigenoperators}
 Many operators that we will consider share the Macdonald functions $H_\mu$ has eigenfunctions. We have taken the habit of calling them Macdonald \define{eigenoperators} to underline this fact. Thus, these operators are entirely characterized by their eigenvalues for the $H_\mu$. The first such, here denoted by \bleu{$\Delta$}, was essentially introduced by Macdonald himself (up to a reformulation in terms of the renormalized $H_\mu$), and it has eigenvalues as follows
    \begin{equation}\label{define_Delta}
             \bleu{\Delta(H_\mu):= B_\mu\,H_\mu}.
    \end{equation} 
Observe that all these eigenvalues are different (since $B_\mu:=\sum_{(i,j)\in\mu} q^it^j$ characterizes $\mu$). It may be shown that
\begin{align}
   \Delta(e_n)=\sum_{k=1}^{n} \frac{q^k-t^k}{q-t}\, e_{n-k}e_{k}.
\end{align}

\subsection{The Nabla operator}
Next, we considered (the author and A.~Garsia, see~\cite{ScienceFiction,nabla}) the \define{operator $\nabla$}, such that
     \begin{equation}\label{define_nabla}
             \bleu{\nabla(H_\mu):= q^{n(\mu')}t^{n(\mu)}\,H_\mu}.
    \end{equation} 
 Its eigenvalues may also be described as $q^{n(\mu')}t^{n(\mu)}=\prod_{(i,j)\in\mu} q^it^j$. 
It has been shown by Haiman (see~\cite{haiman}) that $\nabla(e_n)$ is Schur positive, and that it corresponds to the bigraded Frobenius of the $\S_n$-module of diagonal harmonic polynomials. 
Powers $\nabla^r$, for $r$ in $\N$,  also have nice representation theory interpretations.
It is interesting that we have 
\begin{align}
   &\langle \nabla^r(e_n),e_n\rangle_{q=t=1} =\frac{1}{r\,n+1}\binom{(r+1)n}{n}, 
   		&& \langle \nabla^r(e_n),e_1^n\rangle_{q=t=1} =(r\,n+1)^{n-1},\\
   &\langle \nabla^r(e_n),e_n\rangle_{t=1/q} =\frac{q^{-r\binom{n}{2}}}{[r\,n+1]_q}\qbinom{(r+1)n}{n}, 
		&& \langle \nabla^r(e_n),e_1^n\rangle_{t=1/q} =q^{-r\binom{n}{2}}[r\,n+1]_q^{n-1}. 
\end{align}
 These formulas come from the two following more global ones
 \begin{align}
   &\nabla^r(e_n)_{q=t=1} =\frac{1}{r\,n+1}e_n[(r\,n+1)\bm{x}], \qquad{\rm and}\\
   &\nabla^r(e_n)_{t=1/q} =q^{-r\binom{n}{2}}\frac{1}{[r\,n+1]_q}e_n\big[[r\,n+1]_q\bm{x}\big].
\end{align}
We write $\widetilde{\nabla}$ for the specialization of $\nabla$ at  $t=1$. The operator $\widetilde{\nabla}$ is multiplicative (see exercise), and it is  characterized by the fact that $\widetilde{\nabla} h_\mu^*(\bm{x})= q^{n(\mu')} h_\mu^*(\bm{x})$, since $h_\mu^*(\bm{x})$ is proportional to $H_\mu(q;1,\bm{x})$. The specialization of $\nabla$ at  $t=1/q$ is denoted by $\widehat{\nabla}$, and it is  characterized by the fact that $\widehat{\nabla} s_\mu^*(\bm{x})= q^{n(\mu')-n(\mu)} s_\mu^*(\bm{x})$, since $s_\mu^*(\bm{x})$ is proportional to $H_\mu(q;1/q,\bm{x})$.

It may be shown that (see definition~\ref{defn_staircase} of $\delta_n$)
   \begin{equation}\label{riser_formula}
       \widetilde{\nabla}(e_n)=\sum_{\mu\subseteq \delta_n} q^{\rm area(\mu)}\, s_{(\mu+1^n)/\mu}(\bm{x}),
   \end{equation}
 where \define{${\rm area}(\mu)$} is defined as $\binom{n}{2}-|\mu|$. One may understand this in terms of \define{$q$-Lagrange inversion}. More explicitly, the series 
  \begin{align*}
      \mathcal{F}(q,z;\bm{x})&=\sum_{n\geq 0}  \widetilde{\nabla}(e_n) \,z^{n+1},
   \end{align*}    
satisfies the  functional equation 
  \begin{equation}\label{q_lagrange}
    \mathcal{F} = z\, E\circ_q \mathcal{F},\qquad {\rm with}\qquad E(z)=1+e_1z+e_2z^2+\ldots
  \end{equation}
involving the \define{$q$-composition}\footnote{Careful, this operation is not associative See~\cite{garsia_lagrange,garsia_gessel}.} 
   $$\bleu{f \circ_q g:=\sum_{n\geq 0} f_n\, g(z)\,g(qz)\cdots g(q^{n-1}z)}, \qquad {\rm when}\qquad f(z)=\sum_{n\geq 0} f_n z^n.$$
From this we may deduce (see also Flajolet~\cite[page 330, formula (74)]{flajolet})  that 
 generating series
   \begin{align*}
      F(q;z)&=\sum_{n\geq 0} \langle \widetilde{\nabla}(e_n),e_n\rangle \,z^n,
   \end{align*}    
is given by the formula
    \begin{equation}\label{def_F}
       F(q;z)=\frac{\bm{e}_q(z)}{\bm{e}_q(z/q)},\qquad {\rm where}\qquad \bm{e}_q(z):=\sum_{n\geq 0} \frac{(-z)^n\,q^{n^2}}{(1-q)^n\,[n]_q!}.
  \end{equation}
Beside these interesting facts, our original experimental observation about $\nabla$ was the following. Let us write 
       $$\widehat{s}_\mu:=\Big(\frac{-1}{q\,t}\Big)^{\iota(\mu)}\,s_\mu,$$
     where $\iota(\mu)$ is equal to the sum of the values $\mu_i-i$ which are positive. Then, we have  
 \begin{conjecture}[FB. and A.~Garsia 1999]
   For all $\mu$,  the symmetric function $\nabla(\widehat{s}_\mu)$ is Schur positive.
 \end{conjecture}  
\end{chapter}
This is still widely open. Observe that we can  readily calculate $\widetilde{\nabla}(\widehat{s}_\mu)$ using Jacobi-Trudi. 

\subsection{The \texorpdfstring{$\Delta_f$}{Df}-operators}
In general, for any symmetric function $f$, we set
$$\bleu{\Delta_f(H_{\mu}):=f[B_\mu] H_{\mu}},$$
It is easy to see that, on homogeneous symmetric function of degree $n$, we have\footnote{Careful though: this is not so for other symmetric functions.} $\nabla:=\Delta_{e_{n}}$.
It is also clear that 
 $$ \Delta_{f+g}=\Delta_f +\Delta_g, \qquad{\rm and}\qquad \Delta_{fg}=\Delta_f\circ\Delta_g.$$
 It is known that 
   \begin{equation}\label{lien_pi_n_e_n}
      \nabla(\pi_n)=\Delta_{e_{n-1}} e_n,\qquad {\rm and}\qquad \pi_n=\Delta_{e_1}\widehat{h}_n,
   \end{equation}
  with  $\widehat{h}_n=\widehat{s}_{(n)}$, and thus stands for $(-1/q\,t)^{n-1} h_n$. Many relevant operator identities may be found in~\cite{nabla}, see in particular those listed in (I.12) {\it loc. cit}. These imply that
    \begin{equation}\label{lien_h_n1_pi_n}
     e_1^\perp \nabla(e_n) =\nabla\Big(\sum_{k=0}^n s_{k}(q,t) \,e_{k}e_{n-k}\Big),\qquad {\rm and}\qquad e_1^\perp \nabla(\widehat{h}_n)=\nabla(\pi_n),   
     \end{equation}   
   as well as
     $$(-1)^{n-1}\,M\,\Delta_{e_{n-1}}(p_n[\bm{x}/M]) =  \nabla(\widehat{h}_n),$$
   where $\bleu{M:=(1-q)(1-t)}$.
We also consider $\mathcal{D}_0$ the Macdonald eigenoperator whose eigenvalue for $H_\mu$ is  
 \begin{equation}\label{D_0_eigenvalue}
 1-(1-q)(1-t)B_\mu=1-(1-t)\,(1-q) \sum_{(i,j)\in \mu} q^i t^j.
\end{equation}
The \define{Delta Conjecture} (see~\cite{delta}) proposes an explicit combinatorial formula for $\Delta_{e_k}(e_n)$. See exercises.

\section{Exercises and problems}    
\begin{exer}\rm
\begin{enumerate}
\item[(a)] Show that
\begin{equation}\label{define_Hmu1}
   H_\mu(q,1;\bm{x})=\frac{h_\mu^*(\bm{x})}{h_\mu^*(1)}.
\end{equation} 
\item[(b)] Conclude that
     \begin{equation}
            \widetilde{\nabla}(e_n)=\sum_{\mu\vdash n} f_\mu[1-q]\, q^{n(\mu')} h_\mu^*(\bm{x}).
     \end{equation}
 \item[(c)] Show that
     \begin{equation}
        H_\mu(q,1/q;\bm{x})= \frac{s_\mu^*(\bm{x})}{s_\mu^*(1)}.
     \end{equation}
 \item[(d)] Conclude that
     \begin{equation}
            \widehat{\nabla}(e_n)=\frac{q-1}{q}\sum_{k=1}^n (-1)^{n-k}q^{(n+1)(2k-n)/2}\, s_{(k,1^{n-k})}^*(\bm{x}).
     \end{equation}
\end{enumerate}
\end{exer}

\begin{exer}\rm
The \define{riser sequence} $\rho(\mu)=(r_1,\ldots,r_k)$, of a partition $\mu$ siting inside $\delta_n$, is defined by setting
   $$r_i:=\mu'_{i-1}-\mu'_{i},\qquad ({\rm with}\ \mu'_0:=n).$$
These are the column height differences in $\mu$. 
See that $e_{\rho(\mu)}(\bm{x})=s_{(\mu+1^n)/\mu}(\bm{x})$.

\end{exer}

\begin{exer}\rm  
Recall that the Macdonald polynomial $H_n$ does not depend on $t$ (see~\pref{defn_Hn}). 
\begin{enumerate}
 \item[(a)] Show that $H_b\circ H_a-H_a\circ H_b$ is divisible by $1-q$. To this end, calculate $(H_b\circ H_a)_{q=1}$.
 
 \item[(b)] Calculate the limit, as $q\to 1$, of 
    $$\frac{H_b\circ H_a-H_a\circ H_b}{1-q},$$
  and prove that it is Schur positive when $b>a$.
  
  \item[(c)] Show that 
     $$\lim_{q\to 0} \frac{H_b\circ H_a-H_a\circ H_b}{1-q} =h_b\circ h_a-h_a\circ h_b.$$
     
  \item[(d)] { (Conjecture FB., 2017)} Show that $(H_b\circ H_a-H_a\circ H_b)/(1-q)$ is Schur positive.
\end{enumerate}
\end{exer}

\begin{exer}\rm  
\begin{enumerate}
\item[(a)] Check that 
\begin{equation}
    n(\mu')=\sum_i \binom{\mu_i}{2},
 \end{equation}
and conclude that   
$\widetilde{\nabla}$ is multiplicative.  
\item[(b)] Give a determinantal formula for $\widetilde{\nabla}(\widehat{s}_\mu)$ (of a matrix with entries of the form $\widetilde{\nabla}(e_n)$), and for $\langle \widetilde{\nabla}(\widehat{s}_\mu),e_n\rangle$ (of a matrix with entries of the form $\mathcal{C}_n(q)$).

\end{enumerate}
\end{exer}

\begin{exer}\rm
\begin{enumerate}
\item[(a)] Show that the function $F(q;z)$ in~\pref{def_F} satisfies the difference equation
  \begin{equation}
      F(q;z/q) =\frac{1}{1-z\, F(q;z)}
  \end{equation}
\item[(b)] Deduce from this equality that $\langle \widetilde{\nabla}(e_n),e_n\rangle = \mathcal{C}_n(q)$  (see equation~\pref{catalan_qrec}).
\item[(c)] Compare the above with the specialization of~\pref{q_lagrange} at $q=1$ when one takes the scalar product with $E=E(1)$. 
\item[(d)] Make sense of the statement: 
 \begin{quotation}
 ``The series $\sum_{n\geq 0} \nabla(e_n)_{q=t=1} z^{n+1}$ is the generic Lagrange inverse.''
 \end{quotation}
 by exploiting the fact that the $e_n$ are algebraically independent (so that they may be specialized as one wishes). See~\cite{lenart} for more on this.
\end{enumerate}
\end{exer}

\begin{exer}\rm
\begin{enumerate}
 \item[(a)] Calculate $\Delta_{e_1}(e_n)$, and conclude that
     \begin{equation}
         \frac{1}{n} \Delta_{e_1}(e_n)\Big|_{q=t=1}= \begin{cases}
    \textstyle \sum_{j=0}^{k} e_{k+j+1}e_{k-j},  & \text{if}\ n=2k+1, \\[5pt]
    \textstyle e_{kk}/2+\sum_{j=1}^{k} e_{k+j}e_{k-j},   &  \text{if}\ n=2k.
    \end{cases}
    \end{equation}
 
\item[(b)] ({\bf explore}) Show that
    \begin{equation}
      \frac{1}{(n)_k}  \Delta_{e_k}(e_n)\Big|_{q=t=1}=\sum_{\audessus{\mu\vdash n}{\ell(mu)\leq k}} \frac{1}{d_0!d_1!\cdots d_n!}\  e_\mu
     \end{equation}
where $\mu$ is packed with $0$'s to make it of length $k$, and $d_i=d_i(\mu)$ is the multiplicity of $i$ in $\mu$ (including the $0$ parts). We use here the notation $(n)_k=n(n-1)\cdots (n-k+1)$.	
\end{enumerate}
\end{exer}

\begin{chapter}{\bleu{Rectangular combinatorics, and symmetric functions}}
In recent years, a lot of attention has been given to ``Rational'' Catalan combinatorics (see~\cite{ALW}) associated to pairs $(a,b)$ of co-prime integers, with $a$ and $b$ the proportions of a $(m,n)$-stairshape partition ($\mu\subseteq \delta_{mn}$),  inside which one consider sub-partitions (aka $(a,b)$-Dyck paths). More recently, this has been expanded (see~\cite{open}) to ``Rectangular'' Catalan combinatorics (removing the co-primality condition), considering any pair $(m,n)$ in $\N\times \N$. The greatest common divisor $d=\gcd(m,n)$ playing a special role here,  it is natural to think that $(m,n)=(ad,bd)$, so that the two approaches are linked together.

This is closely tied to the study of $(m,n)$-indexed symmetric functions that arise from ``creation'' operators\footnote{In Theoretical Physics, this is a classical notion which typically builds an object out of the empty one (for us this will be the symmetric function $s_0=1$), thus the name.} on symmetric functions. 
These symmetric functions play a fundamental role in many areas of recent research, such as: Algebraic Combinatorics (Delta Conjecture~\cite{delta}), Representation Theory (graded modules of bivariate harmonic polynomials~\cite{open}), Knot Theory (Khovanov-Rozansky homology of $(m,n)$-torus links~\cite{hogancamp}), Algebraic Geometry (Hilbert scheme of points in the plane~\cite{haiman}), Theoretical Physics (quantum field theory~\cite{drinfeld,schiffmann}), and more.

\section{Rectangular Catalan combinatorics}
In preparation for our upcoming extensions of ``Macdonald operator'' theory, we consider the $(m\times n)$-rectangular extension of the notions of ``{Dyck paths}''
and ``{parking functions}''. In our context, these correspond to sub-partitions $\mu$ of the \define{$(m,n)$-staircase} shape $\delta_{mn}$, and standard tableaux of skew shape $(\mu+1^n)/\mu$. The \define{$(m,n)$-staircase shape}, denoted by $ \delta_{mn}$, (see figure~\ref{mn_staircase}) is set to be
\def\facteur{.95}
\def\dx{0.01}
\def\dy{-.99}
 \begin{figure}[ht]
\begin{center}
 \begin{tikzpicture}[thick,scale=.4]
 \coordinate (NW) at (0,7);
 \coordinate (SE) at (9,0);
 \draw[step=1.0,black,thin] (0,0) grid (9,7);
 \draw[blue] (NW) to (SE) ;
\crochet{0}{7}{0}{1}
\crochet{0}{6}{1}{1}
\crochet{1}{5}{1}{1}
\crochet{2}{4}{1}{1}
\crochet{3}{3}{2}{1}
\crochet{5}{2}{1}{1}
\crochet{6}{1}{1}{1}
 \draw[red,ultra thick] (7,0) to (9,0) ;
\cellrow{6}{0}
\cellrow{5}{1}
\cellrow{4}{2}
\cellrow{3}{4}
\cellrow{2}{5}
\cellrow{1}{6}
\end{tikzpicture}
\qquad\qquad
\begin{tikzpicture}[thick,scale=.4]
 \draw[step=1.0,black,thin] (0,0) grid (4,9);
 \draw[blue] (0,9) to (4,0) ;
\crochet{0}{9}{0}{1}
\crochet{0}{8}{0}{1}
\crochet{0}{7}{0}{1}
\crochet{0}{6}{1}{2}
\crochet{1}{4}{1}{2}
\crochet{2}{2}{1}{2}
 \draw[red,ultra thick] (3,0) to (4,0) ;
\cellrow{6}{0}
\cellrow{5}{0}
 \cellrow{4}{1}
\cellrow{3}{1}
\cellrow{2}{2}
\cellrow{1}{2}
\end{tikzpicture}
\end{center}
\qquad  \vskip-15pt
\caption{The $(9,7)$-staircase: $765321$, and the $(4,9)$-staircase: $332211$.}
\label{mn_staircase} 
\end{figure}
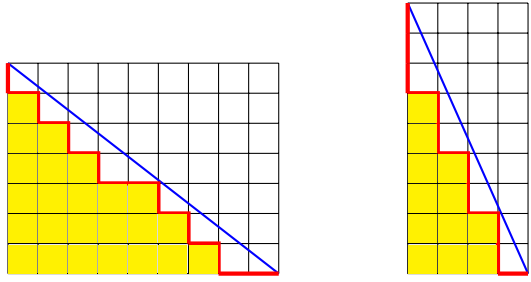
\vskip-15pt
   \begin{equation}
      \bleu{\delta_{mn}:=r_1r_2\cdots r_{n}},\qquad{\rm with}\qquad \bleu{r_k:=\lfloor m\,(n-k)/n\rfloor}.
   \end{equation}

More precisely, we define the \define{$(m,n)$-Dyck path} associated to sub-partition $\mu$ of $\delta_{mn}$,
as the path that goes from $(0,n)$ to $(m,0)$ with either south or east steps following the grid $\N\times \N$, following the boundary of $\mu$. 
In other terms, the cells of $\mu$ are precisely those that lie bellow the Dyck path. We denote by $\mu$ both notions. The \define{area} of a $(m,n)$-Dyck path $\mu$ is set to be
\begin{equation}
   \bleu{\area(\mu):=|\delta_{mn}|-|\mu|},
\end{equation}
hence, it is the number of cells that lie between the path $\mu$ and the path $\delta_{mn}$.  Further along, we will need to encode $\delta_{mn}$ as a monomial in the variables  $\bm{z}=z_0,z_1,z_2,\ldots,z_{m-1}$, setting
\begin{displaymath}
   \bleu{\bm{z}_{mn}:=z_{r_1}z_{r_2}\cdots z_{r_n}},\qquad {\rm for}\qquad \bleu{\delta_{mn}=r_1r_2\cdots r_n},
 \end{displaymath}
 where we are careful to include $0$ parts (observe that the largest possible value of the $r_i$'s is $m-1$).

Given a {$(m,n)$-Dyck path} $\mu$, a \define{$(m,n)$-parking function} $\pi$ on $\mu$ is a bijective labeling of the $n$ vertical steps in the path $\mu$, by the elements of the set $\{1,2,\ldots,n\}$, such that consecutive vertical steps have increasing labels reading bottom up. This is clearly the same as a standard tableau 
   $$\pi:(\mu+1^n)/\mu \rightarrow \{1,2,\ldots,n\},$$
of shape $(\mu+1^n)/\mu$, as is illustrated in figure~\ref{parking}.
 \begin{figure}[ht]
\begin{center}
 \begin{tikzpicture}[thick,scale=.4]
 \coordinate (NW) at (0,7);
 \coordinate (SE) at (9,0);
 \draw[blue,opacity=.5] (NW) to (SE) ;
 \draw[red,ultra thick] (7,0) to (9,0) ;
\draw[step=1.0,black,thin] (0,0) grid (9,7);
 \entree07{5} ;
  \entree06{4} ;
   \entree05{1} ;
    \entree24{3} ;
     \entree33{7} ;
      \entree42{6} ;
       \entree41{2} ;
\crochetbleu{0}{7}{0}{3}
\crochetbleu{0}{4}{2}{1}
\crochetbleu{2}{3}{1}{1}
\crochetbleu{3}{2}{1}{2}
 \draw[blue,ultra thick] (4,0) to (9,0) ;
\cellrowbleu{4}{1}
\cellrowbleu{3}{2}
\cellrowbleu{2}{3}
\cellrowbleu{1}{3}
\end{tikzpicture}\end{center}
\vskip-15pt
\caption{The $(9,7)$-parking function $0420043$ on $4432$.}
\label{parking} 
\end{figure}
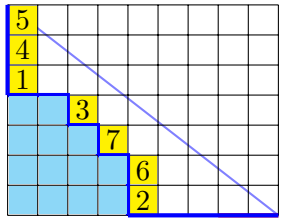

Equivalently, a $(m,n)$-parking function $\pi$  on $\mu$ corresponds to a permutation $\pi_1\pi_2\cdots \pi_n$ of the rows of $\mu$ (packed with $0$'s so that it be of length $n$), with $\pi_i$ equal to the number of cells (of $\mu$) to the left of cells in which $i$ sits. We will also say that the parking function $\pi$ is a \define{$\mu$-parking function}.

We will respectively denote by \define{$\Dyck_{mn}$},  \define{$\Park_{\mu}$}, and by \define{$\Park_{mn}$}, the sets of {$(m,n)$-Dyck paths},  $\mu$-parking functions, and the {$(m,n)$-parking functions}. The cardinalities of these sets are respectively denoted by 
\bleu{$\Cat_{mn}$},  \bleu{$\mathcal{P}_{\mu}$}, and by \bleu{$\mathcal{P}_{mn}$}; and the $q$-enumeration by area of Dyck paths is set to be
\begin{equation}
   \bleu{\Cat_{mn}(q):=\sum_{\mu\subseteq \delta_{mn}} q^{\area(\mu)}}.
   \end{equation}
  For instance, we have
  \begin{eqnarray*}
\Cat_{10}(q)&=&1,\\
\Cat_{21}(q)&=&1,\\
\Cat_{32}(q)&=&1+q,\\
\Cat_{43}(q)&=&1+2\,q+{q}^{2}+{q}^{3},\\
\Cat_{54}(q)&=&1+3\,q+3\,{q}^{2}+3\,{q}^{3}+2\,{q}^{4}+{q}^{5}+{q}^{6},\\
\Cat_{65}(q)&=&1+4\,q+6\,{q}^{2}+7\,{q}^{3}+7\,{q}^{4}+5\,{q}^{5}+5\,{q}^{6}+3\,{q}^{7}+2\,{q}^{8}+{q}^{9}+{q}^{10}.
\end{eqnarray*}

\section{Counting  \texorpdfstring{$(m,n)$}{mn}-Dyck paths and  \texorpdfstring{$(m,n)$}{mn}-parking functions}
Our exposition here will be made clearer if we consider that $(m,n)=(ad,bd)$ with $a$ and $b$ co-prime (relatively prime), so that $d$ is equal to $\gcd(m,n)$. 
\subsection{The co-prime case}
To start with, the enumeration of $(a,b)$-Dyck path (and parking functions) is much easier than the more general case. In a nutshell the number \bleu{$\Cat_{ab}$} of  $(a,b)$-Dyck path
is simply given by the formula
   \begin{equation}\label{cat_ab}
       \Cat_{ab}=\frac{1}{a+b}\binom{a+b}{a}.
    \end{equation}
This may be obtained by a ``classical'' cyclic argument (see exercises), which appears to be due to Dvoretzky-Motzkin (see~\cite{motzkin}), or even earlier to Lukasiewicz. It follows from the same argument that the number \bleu{${\mathcal{P}}_{ab}$} of $(a,b)$-parking functions is given by the formula
   \begin{equation}\label{park_ab}
       \mathcal{P}_{ab}=a^{b-1}.
    \end{equation}
Both of the above formulas may be deduced from (any) of the following (equivalent) more general formulas (see exercises).
   \begin{align}
 &{\rm (a)}& \sum_{\mu\subseteq \delta_{a,b}} s_{(\mu+1^n)/\mu}(\bm{x}) &= \frac{1}{a}\, e_b[a\,\bm{x}],\nonumber \\
 &{\rm (b)}      &&= {\frac{1}{a}\,\sum_{\lambda\vdash b}a^{\ell(\lambda)}\, \frac{(-1)^{b-\ell(\lambda)}} {z_\lambda} }\,p_\lambda(\bm{x}),\nonumber\\   
 &{\rm (c)}     &&={\frac{1}{a}\,\sum_{\lambda\vdash b} \prod_{k\in\mu} \binom{a}{k}\,m_\lambda(\bm{x})}, \label{formules_ab}\\
 &{\rm (d)}      &&={\frac{1}{a}\,\sum_{\lambda\vdash b } s_\lambda'(a)\,s_\lambda(\bm{x})}, \nonumber\\
  &{\rm (e)}     &&={ \sum_{\lambda\vdash n}(-1)^{\ell(\lambda)-1} \frac{(a+1)\cdots (a+\ell(\lambda)-1)}{d_1(\lambda)!\cdots  d_n(\lambda)! }\,h_\lambda(\bm{x})},\nonumber
   \end{align} 
where $d_i(\lambda)$ is the number of parts of size $i$ in $\lambda$ in this last expression.

\subsection{Bizley's formula}
It is a bit harder is to figure out the formula without the co-prime condition. Although this was entirely done in the 1950's by Grossmann and Bizley (see~\cite{bizley}), it remained largely unknown\footnote{Maybe this is in part because the paper was published in a journal for actuaries.} by the community until the early 2000's. 
One finds the following values
\begin{table}[ht]\renewcommand{\arraystretch}{1.5}
\begin{center}
\begin{tabular}{c||c|c|c|c|c|c|c|c|c|c|c|}
$n\setminus m$&1&2&3&4&5&6&7&8&9\\
\hline\hline
1&1&1&1&1&1&1&1&1&1\\ 
\hline
2&1&2&2&3&3&4&4&5&5\\ 
\hline
3&1&2&5&5&7&12&12&15&22\\ 
\hline
4&1&3&5&14&14&23&30&55&55\\ 
\hline
5&1&3&7&14&42&42&66&99&143\\ 
\hline
6&1&4&12&23&42&132&132&227&377\\ 
\hline
7&1&4&12&30&66&132&429&429&715\\ 
\hline
\end{tabular}
\end{center}
\medskip
\caption{Number of $(m,n)$-Dyck paths.}\label{tab1}
\end{table}

The general formulas is 
   \begin{equation}\label{bizley_mn}
      \Cat_{mn} = \sum_{\mu\vdash d} \frac{1}{z_\mu}\prod_{k\in\mu} \frac{1}{a+b} \binom{ak+bk}{ak},
\end{equation} 
with the product over the parts of $\mu$.
In other words, for each fixed $a$ and $b$ co-prime, one has the generating function
\begin{equation}\label{bizley_gen}
   \sum_{d=1}^\infty  \Cat_{(ad,bd)}\ x^d = \exp\!\left(\sum_{j\geq 1} \frac{1}{a+b} \binom{aj+bj}{aj}\, x^j/j\right).
\end{equation} 
\begin{table}[ht]\renewcommand{\arraystretch}{1.5}
\begin{center}
\begin{tabular}{c||c|c|c|c|c|c|c|c|}
$n\setminus m$&1&2&3&4&5&6&7&8\\
\hline\hline 
1&1&1&1&1&1&1&1&1\\
\hline
2&1&3&3&5&5&7&7&9\\ 
\hline
3&1&4&16&16&25&49&49&64\\ 
\hline
4&1&11&27&125&125&243&343&729\\
\hline 
5&1&16&81&256&1296&1296&2401&4096\\ 
\hline
6&1&42&378&1184&3125&16807&16807&35328\\ 
\hline
7&1&64&729&4096&15625&46656&262144&262144\\ 
\hline
\end{tabular}
\end{center}
\medskip
\caption{Number of $(m,n)$-parking functions.}\label{tab2}
\end{table}

For $(m,n)$-parking function enumeration, the analogous formula is
\begin{equation}
     	   \mathcal{P}_{mn}=\sum_{\mu\vdash d} \frac{1}{z_\mu}\,\binom{n}{\mu} \prod_{k\in\mu}\frac{1}{a} (ka)^ {kb-1}.
\end{equation} 
All the above formulas follow from
\begin{equation}\label{e_mn_2}
     	   \sum_{\mu\subseteq \delta_{mn}} s_{(\mu+1^n)/\mu}(\bm{x}) = \sum_{\mu\vdash d} \frac{1}{z_\mu}\,\prod_{k\in\mu} \frac{1}{a}\,e_{kb}[ka\,\bm{x}].
\end{equation} 
The left hand side may be modified in several ways, to account for interesting families of partitions. For instance, we have a similar formula
\begin{equation}\label{h_mn_2}
     	   \sum_{\mu\subseteq \delta_{mn}'} s_{(\mu+1^n)/\mu}(\bm{x}) = (-1)^{n-1}\sum_{\mu\vdash d} \frac{(-1)^{d-\ell(\mu)}}{z_\mu}\,\prod_{k\in\mu} \frac{1}{a}\,e_{kb}[ka\,\bm{x}],
\end{equation} 
where $\mu\subseteq \delta_{mn}'$ corresponds to restricting the summation to ``diagonal avoiding'' partitions $\mu$.

\subsection{Constant term formula}
Let us now set
   \begin{displaymath}
       \bleu{\Omega(\bm{z}):=\prod_{k=0}^m\frac{1}{1-z_k}}.
   \end{displaymath}
We get the $q$-enumeration of $(m,n)$-Dyck path as follows (see exercises).
\begin{prop}\label{prop_constant_simple}
  For any $m$ and $n$, we have the constant term formula
    \begin{equation}\label{form_constant}
           {\Cat_{mn} (q) = \left(\frac{\Omega(\bm{z})}{\bm{z}_{mn}}\, \prod_{i=0}^{m} \frac{1}{1-q\,z_{i+1}/z_{i}}\right)\Big|_{\bm{z}^0},}
    \end{equation}
  where $(-)|_{\bm{z}^0}$ means that we take the constant term with respect to the variables $z_k$.
\end{prop} 

\subsection{Cell rank}
For our upcoming constructions, we consider the notion of \define{$(m,n)$-rank} on cells $(x,y)$ in $\N\times \N$:
  \begin{displaymath}
      \bleu{\rho_{mn}(x,y):=mn -nx-my},
  \end{displaymath}
and outline some of its properties. Clearly, $\rho_{m,n}(x,y)=0$ if and only if $(x,y)$ sits on the line of equation $nx+my =mn$, which crosses the $x$-axis at $m$, and positive rank cells sit below this line.
For example, Figure~\ref{table_rank} presents the $(7,5)$-ranks for cells in the first quadrant.
\begin{figure}[ht]
$$
\begin{array} {rrrrrrrrrr}
\vdots &\vdots&\vdots&\vdots&\vdots&\vdots&\vdots&\vdots\\ \noalign{\medskip}
\rouge{0}&-5&-10&-15&-20&-25&-30&-35&\cdots\\ \noalign{\medskip}
\bleu{7}&\rouge{2}&-3&-8&-13&-18&-23&-28&\cdots\\ \noalign{\medskip}
14&\bleu{9}&\rouge{4}&-1&-6&-11&-16&-21&\cdots\\ \noalign{\medskip}
21&16&\bleu{11}&\bleu{6}&\rouge{1}&-4&-9&-14&\cdots\\ \noalign{\medskip}
28&23&18&13&\bleu{8}&\rouge{3}&-2&-7&\cdots\\ \noalign{\medskip}
35&30&25&20&15&\bleu{10}&\bleu{5}&\rouge{0}&\cdots\end {array}
\begin{picture}(0,0)(17.2,4.6)
\put(0,0){\line(1,0){13.6}}
\put(0,0){\line(0,1){8}}
\put(13.6,0){\line(0,1){8}}
\put(0,8){\line(1,0){13.6}}
\end{picture}
$$
\caption{Example of $(m,n)$-rank ($m=7$ and $n=5$).}\label{table_rank}
\end{figure}
 When $a$ and $b$ are co-prime, all the cells of the $(a\times b)$-rectangle have different $(a,b)$-ranks. We can thus order them (without ambiguity) in decreasing values of their $(a,b)$-rank.

\section{Elliptic Hall algebra operators} 
For any given degree $d$ symmetric function $g(\bm{x})$, that we call a \define{seed}, we construct in this section a family of symmetric function $\bm{g}_{mn}(q,t;\bm{x})$ with $(m,n)=(ad,bd)$, for any $a$ and $b$ co-prime as in the previous section. For the respective seeds $e_d$ and $\widehat{h}_d:=(-1/(qt))^{d-1}h_d$, we respectively get symmetric functions $\bm{e}_{mn}$ and $\widehat{\bm{h}}_{mn}$ that are such that
   $$\bm{e}_{mn}(1,1;\bm{x})=\sum_{\mu\subseteq \delta_{mn}} s_{(\mu+1^n)/\mu}(\bm{x}),\qquad{\rm and}\qquad
      \widehat{\bm{h}}(1,1;\bm{x})=\sum_{\mu\subseteq \delta_{mn}'} s_{(\mu+1^n)/\mu}(\bm{x}).$$
The sum in the second formula is restricted to ``diagonal avoiding'' partitions $\mu$, like at the end of the previous section.
As we will see, the occurence of the parameters $q$ and $t$ turn these symmetric functions in a much refined enumeration of $(m,n)$-parking functions. In this, the parameter $q$ accounts for the area, and the parameter $t$ for the 
number of ``diagonal inversions''. 
\begin{enumerate}
\item 
Our first step is to calculate  the coefficients \bleu{$ c_\mu=c_\mu(q,t)$} of the expansion of $g(\bm{x})$:
\begin{equation}\label{def_coeff_q_n}
    g(\bm{x})=\sum_{\mu\vdash d} \bleu{c_\mu}\, \pi_\mu(\bm{x})
\end{equation}
in the multiplicative basis $\pi_\mu(\bm{x})=\pi_\mu(q,t;\bm{x})$ (see definition~\pref{defn_pi_mu}). 

\item 
Next, for any co-prime pair $(a,b)$ in $\N\times \N$, having constructed (in the next section)  raising degree operators $\bm{Q}_{(a,b)\mu}$ on symmetric functions, we set
   \begin{equation}
        \bleu{\bm{g}_{mn}(q,t;\bm{x}):=\sum_{\mu\vdash d} \bleu{c_\mu}\, \bm{Q}_{(a,b)\mu}(1)},\qquad {\rm where}\qquad (m,n)=(ad,bd).
      \end{equation}
    The operator  $\bm{Q}_{(a,b)\mu}$ raises degrees by $n$. As it is here applied to the symmetric function ``$1$'', we get $\deg (\bm{g}_{mn})=n$.
 \end{enumerate}
We thus have a process 
      $$g_d\to \{\bm{g}_{mn}\}_{(m,n)=(ad,bd)},$$
that maps a symmetric function $g=g_d$ in $\Lambda_d$,  to a family $\{\bm{g}_{mn}\}$ indexed by pairs of the form $(ad,bd)$, with $a$ and $b$ co-prime.

Among of the interesting  identities in this setup (see~\cite{open})  we have the operator equality 
\begin{equation}\label{nabla_property}  \nabla\bm{Q}_{(a,b)\mu}\nabla^{-1}=\bm{Q}_{(a+b,b)\mu}.\end{equation}
Moreover, $\bm{Q}_{(0,1)\mu}(1)=\pi_\mu(\bm{x})$.
It follows that
        \begin{equation}\label{formule_nabla}
      \bm{g}_{0d}(\bm{x})=g_d(\bm{x}),\qquad {\rm and}\qquad       \nabla(\bm{g}_{mn})=  \bm{g}_{m+n,n},
          \end{equation}
hence $\bm{g}_{rd,d}=\nabla^r(g_{d})$. We have thus generalized to the $(m,n)$-rectangular context and to other seeds, our previous ``squarre'' situation corresponding to the study of $\nabla(e_d)$, hence the seed $e_d$.

\section{Construction of the operators}
Our calculations are going to be done inside an operator algebra realization of the elliptic Hall algebra (see~\cite{open}), which is graded\footnote{In fact it is graded over $\Z\times\Z$, but we will only use the ``positive'' part.} over $\N\times \N$. We will not actually need to describe the whole algebra, but use (without further details) the fact that collinearly graded operators commute, as well as property~\pref{nabla_property}. In fact, we only work with the portion of the algebra that is algebraically generated by the two operators: 
\begin{enumerate}
\item  multiplication by $p_1(\bm{x})$, and 
\item the $\mathcal{D}_0$ Macdonald eigenoperator, see~\pref{D_0_eigenvalue}.
\end{enumerate}
To each pair $(m,n)$ in $\N\times\N$, we recursively associate a bracketings of these operators, via the following \define{splitting} procedure  of $(m,n)$.
One may show (see exercise) that there are precisely $d$ positive integer entries matrices such that
\begin{equation}\label{choix_matrice}
    \det \begin{pmatrix} r & s\\ u & v \end{pmatrix} =\gcd(m,n),\qquad {\rm and}\qquad  (m,n)=(r,s)+(u,v).
 \end{equation}
One then recursively sets
   \begin{equation}
        \bleu{ \bm{Q}_{mn}:=\frac{1}{(1-q)(1-t)} \left[ \bm{Q}_{rs},\bm{Q}_{uv}\right]},
        \quad {\rm with}\quad \bleu{\bm{Q}_{01}:=p_1},
        \quad {\rm and}\quad \bleu{\bm{Q}_{10}:=\mathcal{D}_0};  \end{equation} 
         and then finally
\begin{equation}
    \bleu{\bm{Q}_{(a,b)\mu}:=\bm{Q}_{(a\mu_1,b\mu_1)}\cdots \bm{Q}_{(a\mu_k,b\mu_k)}},\qquad {\rm for}\qquad \mu=\mu_1\cdots \mu_k\vdash d.
  \end{equation}
 This is a composition of operators having collinear gradings, hence they commute (making the definition independent of a choice of order).    
From general properties of the elliptic Hall algebra, it may also be shown that this construction does not depend on the choice of a solution for~\pref{choix_matrice}. 
 
  For example, writing $\bleu{M:=(1-q)(1-t)}$ as before, we have
   $$\bm{Q}_{43}=\frac{1}{M^6} [[p_{{1}},\mathcal{D}_{{0}}],[[p_{{1}},\mathcal{D}_{{0}}],[[p_{{1}},\mathcal{D}_{{0}}],\mathcal{D}_{{0}}]]],$$
and
  $$\bm{Q}_{63}=\frac{1}{M^8}[[p_{{1}},\mathcal{D}_{{0}}],[[[p_{{1}},\mathcal{D}_{{0}}],\mathcal{D}_{{0}}],[[[p_{{1}},\mathcal{D}_{{0}}],\mathcal{D}_
{{0}}],\mathcal{D}_{{0}}]]].$$

\section{Interesting seeds}
There are many seeds $g_d$, for which the corresponding $\bm{g}_{mn}$'s all appear to be Schur positive. Many are also tied to interesting combinatorics.
We consider here the three cases $e_d$, $\widehat{h}_d$, and $\pi_d$, with the respective families denoted by $\bm{e}_{mn}$, $\widehat{\bm{h}}_{mn}$, and $\bm{\pi}_{mn}$. In view of~\pref{formule_nabla}, we have
    $$\bm{e}_{rn,n}=\nabla^r(e_n), \qquad \widehat{\bm{h}}_{rn,n}=\nabla^r(\widehat{h}_n),\quad{\rm and}\quad \bm{\pi}_{rn,n}=\nabla^r(\pi_n).$$
 It may be seen that, for all $r$ and $n$, 
 \begin{equation}\label{plus_moins_un}
    \bm{e}_{rn+1,n}= \bm{e}_{rn,n},\qquad {\rm and} \qquad  \bm{e}_{rn-1,n}= \bm{\widehat{h}}_{rn,n}.
 \end{equation}
Observe that in the co-prime situation, one has $\bm{e}_{ab}= \widehat{\bm{h}}_{ab}= \bm{\pi}_{ab}$, illustrating that the more interesting situation is when co-primality does not hold. Part of the conjectures in~\cite{BGLX} is to propose the following explicit combinatorial formula 
  \begin{align}
    &\bm{e}_{mn}(q,t;\bm{x})=\sum_{\mu\subseteq \delta_{mn}} \sum_{\tau \in P_\mu}
q^{{\rm area}(\mu)} t^{{\rm dinv}(\tau)} \bm{x}_\tau,\label{formule_compositional}
  \end{align}
  where $P_\mu$ stands for the set of semi-standard tableaux of skew-shape $(\mu+1^n)/\mu$ (recall that $\mu+\lambda$ denotes the partition obtained by pointwise addition of parts), having values in $\N$. Broadly speaking, ${\rm dinv}$ is a statistics that depends (in some manner not described here) on how entries of cells of  $\tau$ compare, taking into account some ``rank function'' for cells that depends on $m$ and $n$. For more on ``dinv'', see~\cite{HHLRU}. This relates to the ``sweep'' map considered in~\cite{thomas_williams}.

\subsection{Specializations at \texorpdfstring{$t=1$}{t}}
It appears (although this is still to be proven, see~\cite{open}) that, for any symmetric function $g(\bm{x})$  independent of $t$, one has 
     \begin{equation}\bm{Q}_{(a,b)\mu}(g(\bm{x}))\big|_{t=1} = \bm{Q}_{(a,b)\mu}(1)\big|_{t=1}\cdot\,g(\bm{x}).
               \end{equation}
In other words, at $t=1$, the effect of the operator $\bm{Q}_{(a,b)\mu}\big|_{t=1}$ on any $g(\bm{x})$ is to multiply it by the fixed symmetric function 
    $$\bm{\pi}_{(a,b)\mu}(q,1;\bm{x})=\prod_{k\in\mu} \bm{\pi}_{(ak,bk)}(q,1;\bm{x}),$$
  reducing all calculations to that of polynomials in the $\bm{\pi}_{mn}(q,1;\bm{x})$'s. In fact, this ``multiplicativity'' property holds for any other multiplicative basis. For instance, we have that
      $$g(\bm{x})=\sum_{\mu\vdash d} c_\mu\,e_\mu(\bm{x})\qquad {\rm implies}\qquad
       \bm{g}_{mn}(q,1;\bm{x})=\sum_{\mu\vdash d} c_\mu\, \bm{e}_{(a,b)\mu}(q,1;\bm{x}),$$
    for any co-prime $(a,b)$, and $\bm{e}_{(a,b)\mu}(q,1;\bm{x})=\prod_{k\in\mu} \bm{e}_{(ak,bk)}(q,1;\bm{x})$. The individual $\bm{e}_{mn}(q,1;\bm{x})$ afford a combinatorial description (see~\cite{BGLX}) of the form
     \begin{equation}\label{formule_emn_dyck}
             \bm{e}_{mn}(q,1;\bm{x}) =\sum_{\mu \subseteq \delta_{mn}} q^{{\rm area}(\mu)}  s_{(\mu+1^n)/\mu}(\bm{x}).
     \end{equation}
\subsection{Specializations at  \texorpdfstring{$t=1/q$}{t}}
It is also interesting to recall (see~\cite{open}) that
\begin{equation} \label{special_t_1_q}
q^{\alpha(m,n)}\,\bm{\pi}_{mn}(q,1/q;\bm{x})=\frac{[d]_q}{[m]_q}\,e_n\Big[[m]_q\, \bm{x}\Big],
      \end{equation}
 with $d=\gcd(m,n)$, and where $\alpha(m,n)=((n-1)(m-1) +d-1)/2$. 
Further setting $q=1$, we find that
$$\bm{\pi}_{m,n}(1,1;\bm{x})= \frac{d}{m}e_n[m\,\bm{x}],$$
where one considers  $m$, $r$ and $n$ as a constant for the calculation of this plethysm. It follows from~\pref{plus_moins_un} that 
 \begin{displaymath}
\bm{e}_{rn,n}(1,1;\bm{x})= \frac{1}{(r\,n+1)}e_n[(r\,n+1)\,\bm{x}],\quad{\rm and}\quad
\bm{\widehat{h}}_{rn,n}(1,1;\bm{x})= \frac{1}{(r\,n-1)}e_n[(r\,n-1)\,\bm{x}].
\end{displaymath}
These formulas imply that
\begin{align*}
  &\langle \bm{e}_{rn,n}(1,1;\bm{x}),e_1^n\rangle  =(r\,n+1)^{n-1}, \quad 
  &&\langle \bm{e}_{rn,n}(1,1;\bm{x})) ,e_n\rangle =\frac{1}{rn+1}\binom{(r+1)\,n}{n},\\
   &\langle\bm{\widehat{h}}_{rn,n}(1,1;\bm{x})) ,e_1^n\rangle=(r\,n-1)^{n-1}, \quad 
   &&\langle \bm{\widehat{h}}_{rn,n}(1,1;\bm{x})) ,e_n\rangle =\frac{1}{rn-1}\binom{(r+1)\,n-2}{n},\\
 &\langle \bm{\pi}_{mn}(1,1;\bm{x})) ,e_1^n\rangle=d\,m^{n-1}, \quad 
 &&\langle \bm{\pi}_{m,n}(1,1;\bm{x})) ,e_n\rangle =\frac{d}{m+n}\binom{m+n}{n}.
\end{align*}

Moreover, for $m$ and $n$ co-prime we have
 \begin{equation} \label{e_mn_special_t_1_q}
   \bm{\pi}_{mn}(q,1/q;\bm{x})=\frac{(-1)^{n-1}q^{(m+1)(n-1)/2}}{[n]_qe_m[[{n]_q]}} \left(\Delta_{e_m} p_n(\bm{x})\right)_{t=1/q}.
      \end{equation}

\section{What comes next}
The next step is to expand the previous considerations to $3$ parameters. This leads to the study of intervals in the Tamari lattice, or its extension to the $(m,n)$-context. There are several ways to define the \define{Tamari order}\footnote{See for instance \href{https://en.wikipedia.org/wiki/Tamari_lattice}{https://en.wikipedia.org/wiki/Tamari\_lattice}.} on (classical $m=n$) Dyck paths. The simplest may be to consider these as well parenthesized expressions, with $()()\cdots()$ as \define{smallest element}, and \define{associativity} used for covering relation. Thus the \define{largest} element is $(((\cdots)))$. See next subsection for an explicit description in terms of partitions.

There are also many (well-known) ways of constructing a convex polytope, called the \define{associahedron}\footnote{See \href{https://en.wikipedia.org/wiki/Associahedron}{https://en.wikipedia.org/wiki/Associahedron}.} whose vertices an edges correspond to the Hasse diagram of the Tamari lattice.  The Tamari order corresponds to an orientation of this polytope, that may be specified by choosing a minimal vertex (why?). An example with $m=n=4$ corresponds to figure~\ref{tamari3}. 
 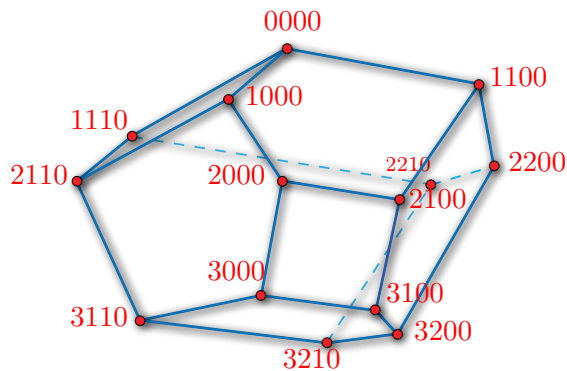
\begin{figure}[ht]
   \setlength{\unitlength}{2.5mm}
\begin{center}
\begin{picture}(0,0)(0,0)
                \put(19,11){$\rouge{{\scriptstyle 2210}}$}  
 \end{picture}
\scalebox{.33}{\includegraphics{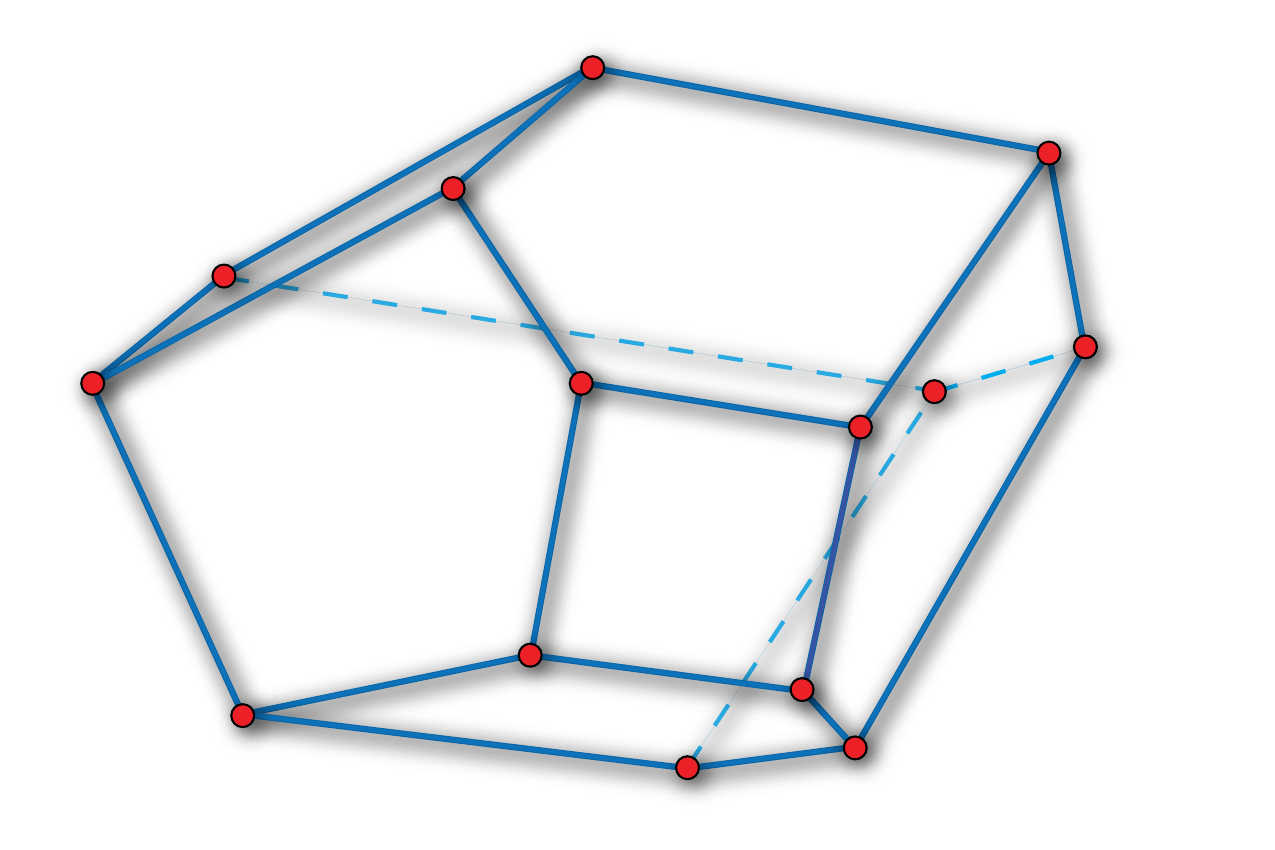}}
\begin{picture}(0,0)(0,0)
               \put(-17,18.5){$\rouge{0000}$}
       \put(-27.3,13.2){$\rouge{1110}$}        \put(-18,14.5){$\rouge{1000}$}      \put(-5,15.5){$\rouge{1100}$}
  \put(-30.5,10.3){$\rouge{2110}$}   \put(-20,10.3){$\rouge{2000}$}  \put(-9.3,9){$\rouge{2100}$}\put(-4,11){$\rouge{2200}$}
                 \put(-20,5.2){$\rouge{3000}$} 
  \put(-27.3,2.7){$\rouge{3110}$}      \put(-10.5,4){$\rouge{3100}$} \put(-9,2){$\rouge{3200}$}
              \put(-16,0.5){$\rouge{3210}$}
\end{picture}
\end{center}
\vskip-8pt
\caption{The $(4,4)$-Tamari poset for $n=4$, bottom element is $3210$.}
\label{tamari3}
\end{figure}

Similarly, the \define{$(m,n)$-Tamari lattice}, \bleu{$\mathcal{T}_{mn}$}, corresponds to an order on the set of $(m,n)$-Dyck paths, with the $m=n$ case corresponding to 
the classical Tamari lattice. Examples include the $(8,4)$-Tamari lattice of figure~\ref{tamari84} (a bit ``exploded'' here to better display its structure). Its smallest element is the one at the forefront (the one with only blue faces adjacent). The blue portion of this figure is a copy of \ref{tamari3}. 
Increasing order corresponds to following shortest paths away from the smallest node. 
Likewise, we have the $(4,6)$-Tamari lattice in figure~\ref{tamari46}. Nice geometrical realizations and combinatorial properties of these posets are studied in~\cite{ceballos}.

\begin{figure}[ht]
\begin{center}
\includegraphics[width=1.8in]{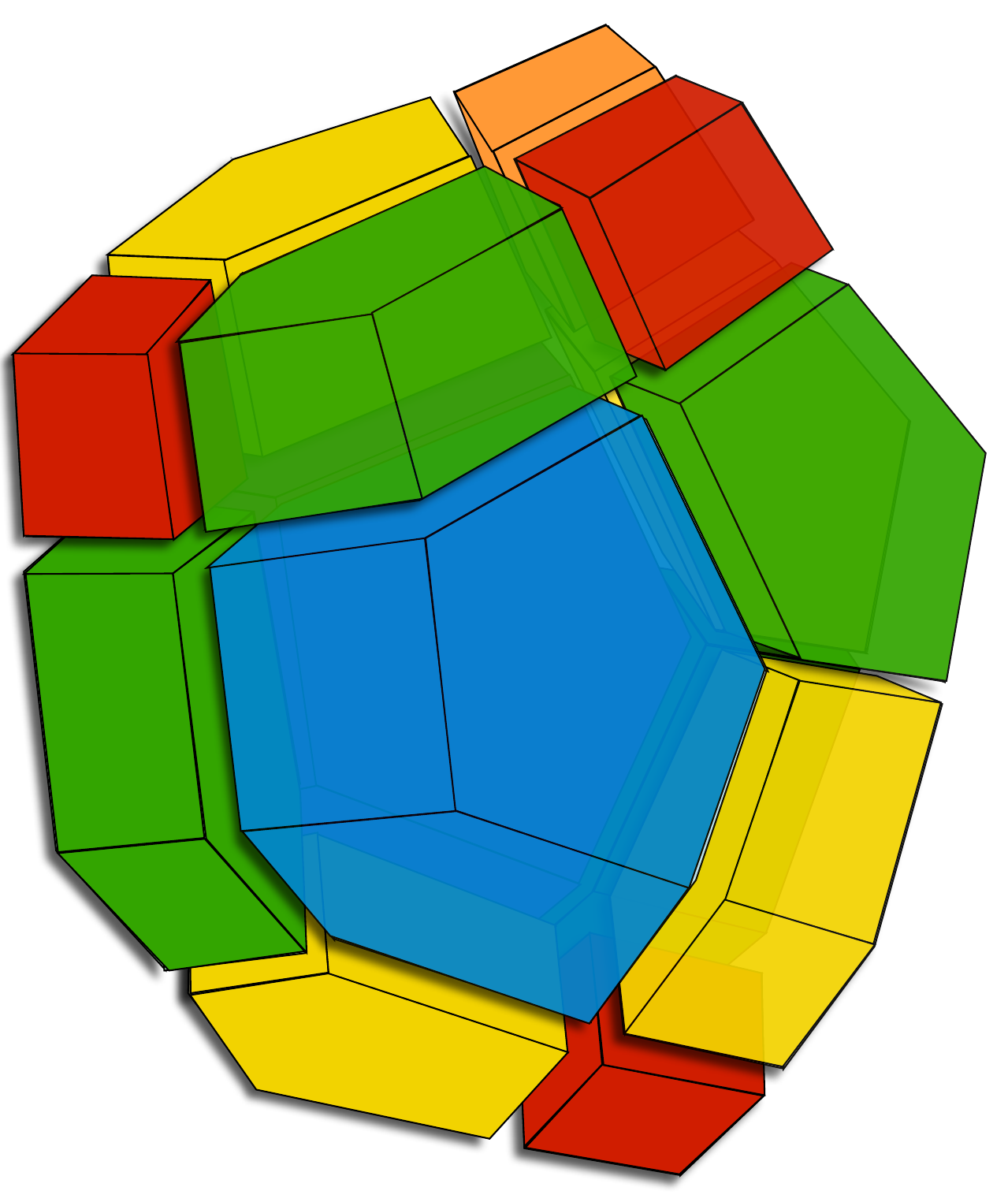}
\caption{The $(8,4)$-Tamari lattice, seen from the bottom.}
\label{tamari84}
\end{center}
\end{figure}

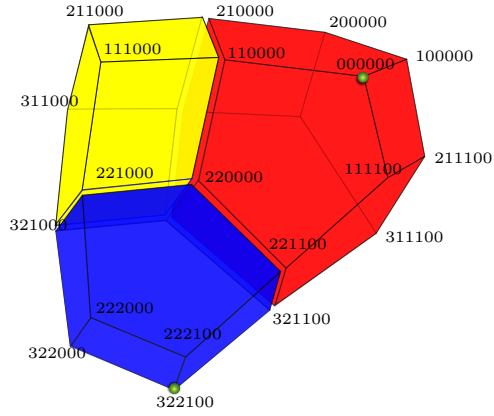
\begin{figure}[ht]
\begin{center}
\includegraphics[width=2in]{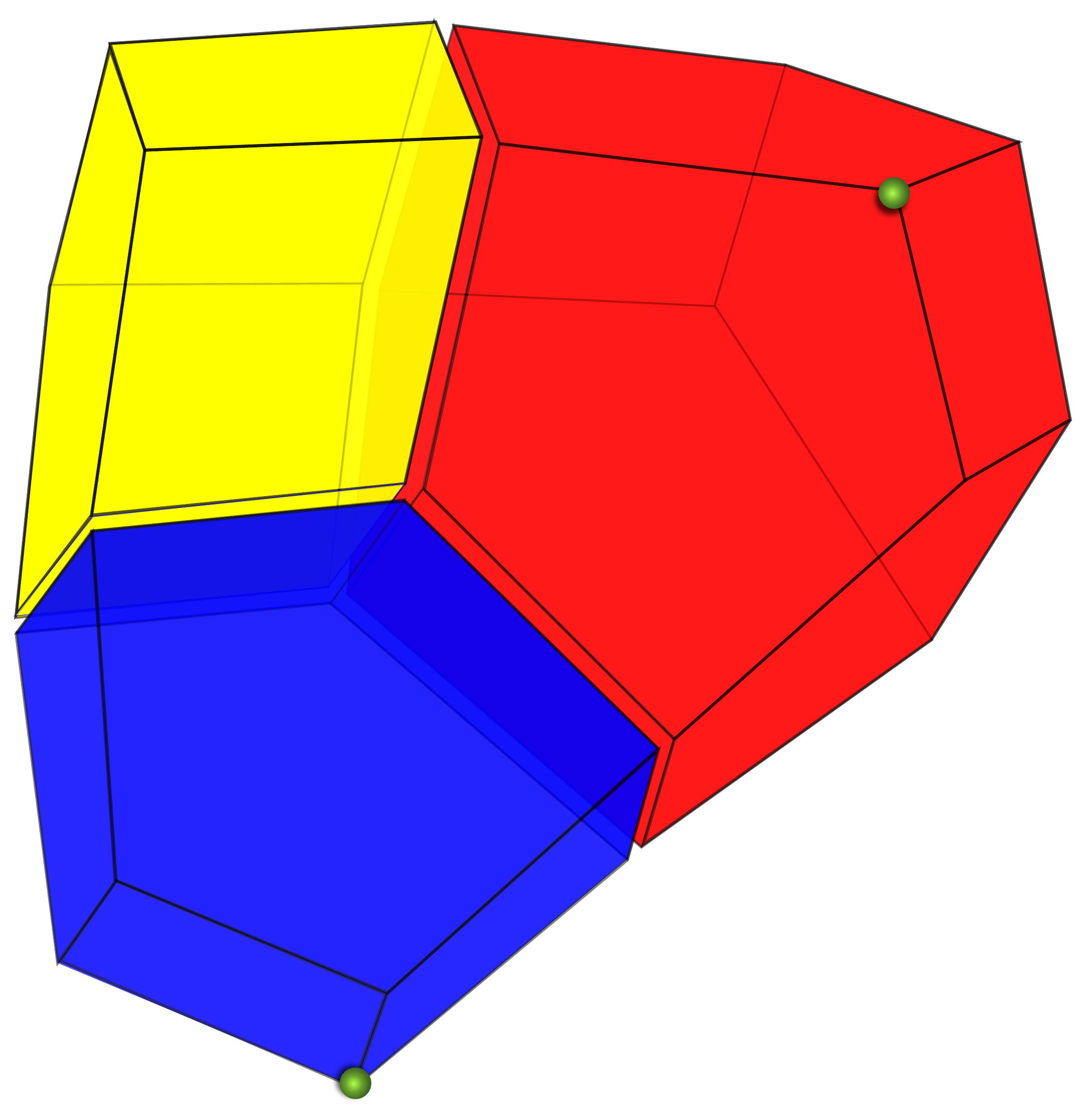}
\begin{picture}(0,0)(0,0)
\setlength{\unitlength}{5mm}
\put(-5.7,9.1){\tiny 110000}
\put(-6,10.2){\tiny 210000}
\put(-3,9.9){\tiny 200000}
\put(-2.8,8.8){\tiny 000000}
\put(-.7,9){\tiny 100000}
\put(-.2,6.3){\tiny 211100}
\put(-6.3,5.8){\tiny 220000}
\put(-2.6,6){\tiny 111100}
\put(-1.5,4.2){\tiny 311100}
\put(-4.5,2){\tiny 321100}
\put(-4.6,4){\tiny 221100}
\put(-7.4,1.6){\tiny 222100}
\put(-7.6,-0.2){\tiny 322100}
\put(-10,10.2){\tiny 211000}
\put(-11.2,7.8){\tiny 311000}
\put(-9,9.2){\tiny 111000}
\put(-9.2,5.9){\tiny 221000}
\put(-9.2,2.3){\tiny 222000}
\put(-11.5,4.5){\tiny 321000}
\put(-11,1.1){\tiny 322000}
\end{picture}
\caption{The $(4,6)$-Tamari lattice, bottom element is $322100$, and top is $000000$.}
\label{tamari46}
\end{center}
\end{figure}

For $m$ equal to $rn+1$, the following formulas (see~\cite{trivariate,bousquet,preville}) account for the number of \define{intervals}, pairs $(\nu,\mu)$ with $\nu\preceq\mu$ in the $(m,n)$-Tamari lattice, and the number  of pairs $(\nu,\pi)$ with $\pi$ a parking function on $\mu$:
\begin{align}
  &\sum_{\nu\preceq\mu\,\in \mathcal{T}_{rn+1,n}} 1=\frac {( r+1)}{ n\,( r n+1) } \binom{( r+1) ^{2}\,n+r}{ n-1}, \\ 
  &\sum_{\nu\preceq\mu\,\in \mathcal{T}_{rn+1,n}} \binom{n}{\rho(\mu)}=(r+1)^n(r\,n+1)^{n-2},
   \end{align}
 where both sums are over pairs $\nu\preceq\mu$ for $\nu$ and $\mu$ in the $(rn+1,n)$-Tamari lattice,
and where  $\rho(\mu):=(c_1,c_2,\ldots,c_k)$ is the \define{sequence of column heights} in the skew-shape $(\mu+1^n)/\mu$.
Both identities come out of the symmetric function formula
\begin{equation}\label{e_mn_3}
   \sum_{\nu\preceq \mu\,\in \mathcal{T}_{rn+1,n} } s_{(\mu+1^n)/\mu}(\bm{x})= \sum_{\lambda \vdash n} (-1)^{n-\ell(\lambda)}(r\, n+1)^{\ell(\lambda)-2}
                      {\textstyle\prod_{k\in\lambda} \binom{(r+1)\,k}{k}}\,\frac{p_\lambda(\bm{x})}{z_\lambda}.
 \end{equation}
 There is a similar formula for $m=rn-1$, namely
\begin{equation}\label{h_mn_3}
   \sum_{\nu\preceq \mu\,\in \mathcal{T}_{rn-1,n}} s_{(\mu+1^n)/\mu}(\bm{x})
=\sum_{\lambda\vdash n} (-1)^{n-\ell(\lambda)}\,((r+1)\,n-1)^{\ell(\lambda)-2}{\textstyle{\prod_{k\in\lambda} \binom{k\,(r+1)-1}{k}}}\,\frac{p_\lambda(\bm{x})}{z_\lambda}.
 \end{equation}
 These formulas extend in a natural manner (explained further below) formulas~\pref{e_mn_2} and \pref{e_mn_2} (see also~\pref{formules_ab}). They also afford extensions with three parameter $q$, $t$, and $r$, such that the case $r=0$ corresponds to the $\bm{e}_{mn}(q,t;\bm{x})$ of this chapter. In particular, formulas~\pref{e_mn_3} and~\pref{h_mn_3} are related in this way to~\pref{plus_moins_un}.
 
 The number of intervals in the $(m,n)$-Tamari lattices are given in table~\ref{tab3}, for $m\leq n\leq 7$. 
\begin{table}[ht]\renewcommand{\arraystretch}{1.5}
\begin{center}
\begin{tabular}{c||c|c|c|c|c|c|c|c|c|c|c|}
$n\setminus m$&1&2&3&4&5&6&7\\
\hline\hline
1 & 1 \\
2 & 1 & 3  \\
3 & 1 & 3 & 13 \\
4 & 1 & 6 & 13 & 68 \\
5 & 1 & 6 & 23 & 68 & 399  \\
6 & 1 & 10 & 58 & 161 & 399 & 2530  \\
7 & 1 & 10 & 58 & 248 & 866 & 2530 & 16965\\
\end{tabular}
\end{center}
\medskip
\caption{Number of intervals in the $(m,n)$-Tamari lattice.}\label{tab3}
\end{table}
 
Similarly table~\ref{tab4} gives the number of pairs $(\nu,\pi)$ (\define{decorated intervals}), with $\pi$ a parking function on $\mu\succeq\nu$.
 \begin{table}[ht]\renewcommand{\arraystretch}{1.5}
\begin{center}
\begin{tabular}{c||c|c|c|c|c|c|c|c|c|c|c|}
$n\setminus m$&1&2&3&4&5&6&7\\
\hline\hline
1 & 1 \\
2 & 1 & 4  \\
3 & 1 & 5 & 32  \\
4 & 1 & 17 & 49 & 400  \\
5 & 1 & 23 & 167 & 729 & 6912  \\
6 & 1 & 72 & 1048 & 4407 & 14641 & 153664  \\
7 & 1 & 102 & 1818 & 15626 & 90079 & 371293 & 4194304\\
\end{tabular}
\end{center}
\medskip
\caption{Number of decorated intervals in the $(m,n)$-Tamari lattice.}\label{tab4}
\end{table}

\section{Conclusion}
There is a global algebraic point of view (see~\cite{bergeron_several}) that unifies most of what we have seen at play in these notes, and ties it to the study of interesting representation theory problems, as well as many other subjects. The formulas in this broader approach take the form of linear combinations  $\mathcal{E}_{mn}(\bm{q};\bleu{\bm{x}})$, with positive integer coefficients, of products $s_\lambda(\bm{q})s_\mu(\bm{x})$, with $\bm{q}=(q_1,q_2,\cdots,q_k)$; and the various specific cases occurring in these notes are obtained by specialization of the $q_i$. 
We illustrate with an example, with $m=n=4$.

Let us set
\begin{align}
     \mathcal{E}_{44}(\bm{q};\bleu{\bm{x}})&:=\bleu{s_4(\bm{x})}\nonumber\\ 
           &\qquad +(s_1(\bm{q})+s_2(\bm{q})+s_3(\bm{q}))\,\bleu{s_{31}(\bm{x})}\nonumber\\ 
           &\qquad +(s_2(\bm{q})+s_{4}(\bm{q})+s_{11}(\bm{q}))\,\bleu{s_{22}(\bm{x})}\label{E_44}\\   
          &\qquad +(s_3(\bm{q})+s_{4}(\bm{q})+s_{5}(\bm{q})+s_{11}(\bm{q})+s_{21}(\bm{q})+s_{31}(\bm{q}))\,\bleu{s_{211}(\bm{x})}\nonumber\\   
         &\qquad+(s_6(\bm{q})+s_{41}(\bm{q})+s_{31}(\bm{q})+s_{111}(\bm{q}))\,\bleu{s_{1111}(\bm{x})}  \nonumber
  \end{align}
Then, we have
$$ \mathcal{E}_{44}(q;\bm{x})=H_n(q,\bm{x}),\qquad \mathcal{E}_{44}(q,t;\bm{x})=\nabla(e_4)=\bm{e}_{4,4}(q,t;\bm{x}),$$
and
$$\mathcal{E}_{44}(1,1,1;\bm{x})=\sum_{\nu\preceq \mu\,\in \mathcal{T}_{4,4} } s_{(\mu+1^n)/\mu}(\bm{x}).$$
In each of these specializations, some of the Schur functions  $s_\lambda(\bm{q})$ evaluate to $0$, simply because there are not enough variables (the number is less than $\ell(\lambda)$). Otherwise they all are essentially the ``same'' formula. It is interesting to observe that setting all the $q_i=1$ gives the formula
\begin{align*}
\mathcal{E}_{44}(k;\bm{x})&= s_{4} 
+ \left({\textstyle \binom{k}{3} + 3 \, \binom{k}{2} + 3 \, \binom{k}{1}}\right)s_{31}\\
 &\qquad + \left({\textstyle \binom{k}{4} + 5 \, \binom{k}{3} + 6 \, \binom{k}{2} + 2 \, \binom{k}{1}}\right)s_{22} \\
 &\qquad 
+ \left({\textstyle \binom{k}{5} + 8 \, \binom{k}{4} + 18 \, \binom{k}{3} + 15 \, \binom{k}{2} + 3 \, \binom{k}{1}}\right)s_{211} \\
 &\qquad 
+\left({\textstyle \binom{k}{6} + 9 \, \binom{k}{5} + 25 \, \binom{k}{4} + 29 \, \binom{k}{3} + 12 \, \binom{k}{2} + \binom{k}{1}}\right)s_{1111},
\end{align*}
with the convention that (or thinking of this as a plethysm with $k$ constant): 
   $$\mathcal{E}_{44}(k;\bm{x}):=\mathcal{E}_{44}(\underbrace{1,1,\cdots,1}_{k\ {\rm copies}};\bm{x}).$$
 In particular, we find that the formulas
\begin{align*}
   &\langle \mathcal{E}_{rn,n}(1;\bm{x}),e^n\rangle =  1, && \langle \mathcal{E}_{rn,n}(1;\bm{x}),p_1^n\rangle =  n!,\\
   &\langle \mathcal{E}_{rn,n}(2;\bm{x}),e^n\rangle =  \textstyle\frac{1}{rn+1}\binom{(r+1)n}{n}, && \langle \mathcal{E}_{rn,n}(2;\bm{x}),p_1^n\rangle =\textstyle(rn+1)^{n-1},\\
   &\langle \mathcal{E}_{rn,n}(3;\bm{x}),e^n\rangle =\textstyle \frac {( r+1)}{ n\,( r n+1) } \binom{( r+1) ^{2}\,n+r}{ n-1}, &&   \langle \mathcal{E}_{rn,n}(3;\bm{x}),p_1^n\rangle =\textstyle (r+1)^n (rn+1)^{n-2},
 \end{align*}
 come  out of the same ``molds''. For instance, these are respectively
$$ \langle \mathcal{E}_{4,4}(k;\bm{x}),e^4\rangle = \textstyle \binom{k}{1}+12 \, \binom{k}{2}+ 29 \, \binom{k}{3} + 25 \, \binom{k}{4} + 9 \, \binom{k}{5} +  \binom{k}{6} ,$$   
and
$$ \langle \mathcal{E}_{4,4}(k;\bm{x}),p_1^4 \rangle =\textstyle 1+ 23 \, \binom{k}{1}+ 78 \, \binom{k}{2} + 96 \, \binom{k}{3} + 51 \, \binom{k}{4}+ 12 \, \binom{k}{5}  + \binom{k}{6}   .$$   
  
  The expressions~$\mathcal{E}_{mn}(\bm{q};{\bm{x}})$ have many interesting structural properties and hidden symmetries that are currently being studied.

\section{Exercises and problems}
\begin{exer}\rm
\begin{enumerate}
\item[(a)] Prove that $\Cat_{rn+1,n}=\Cat_{rn,n}$ and $\mathcal{P}_{rn+1,n}=\mathcal{P}_{rn,n}$ by a simple combinatorial argument.
\item[(b)]    Prove  formula~\ref{cat_ab} via the ``cyclic approach''. For this, you encode a Dyck path as a sequence of $a+b$ letters either equal to ``$\bm{s}$'' (coding a ``south'' step in the path) or ``$\bm{e}$''  (coding an ``east'' step in the path). This sequence is periodically repeated ``infinitely'' in both directions. See that there are precisely $(a+b)$ different length $(a+b)$ subsequences of consecutive letters in this periodic infinite sequence. Each of these correspond to a path in the $a\times b$-rectangle. Show that there is a unique one that stays below the diagonal.
\item[(c)] Suitably adapt your previous argument to prove formula~\ref{park_ab}. You may be inspired by figure~\ref{fig_rect}.
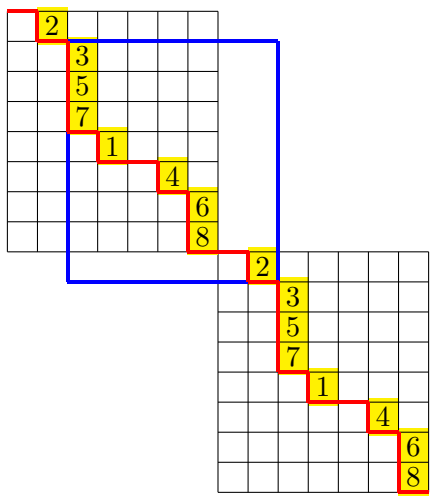
\begin{figure}[ht]
\setlength\unitlength{4mm}
\setlength{\carrelength}{5mm}
\def\jcarre{\put(0,0){\jaune{\linethickness{\carrelength}\line(1,0){1}}}}
\begin{center}
\begin{picture}(14,15)(0,1.5)
\put(7,-.48){ 
                     \multiput(1,8)(0,1){1}{\jcarre}
                    \multiput(2,7)(0,1){1}{\jcarre}
                    \multiput(2,5)(0,1){3}{\jcarre}
                    \multiput(3,4)(0,1){1}{\jcarre}
                    \multiput(5,3)(0,1){1}{\jcarre}
                    \multiput(6,1)(0,1){2}{\jcarre}}
\put(7,0){\multiput(0,0)(0,1){9}{\line(1,0){7}}
\multiput(0,0)(1,0){8}{\line(0,1){8}}}
\put(0,8){\put(0,-.48){ 
                    \multiput(1,8)(0,1){1}{\jcarre}
                    \multiput(2,7)(0,1){1}{\jcarre}
                    \multiput(2,5)(0,1){3}{\jcarre}
                    \multiput(3,4)(0,1){1}{\jcarre}
                    \multiput(5,3)(0,1){1}{\jcarre}
                    \multiput(6,1)(0,1){2}{\jcarre}}
\put(0,0){\multiput(0,0)(0,1){9}{\line(1,0){7}}
\multiput(0,0)(1,0){8}{\line(0,1){8}}}}
  \linethickness{.5mm}
 \put(2,7){\bleu{\put(0,0){\line(1,0){7}}
                        \put(7,0){\line(0,1){8}}
                        \put(0,0){\line(0,1){8}}
                        \put(0,8){\line(1,0){7}}}}
\put(7,0){\put(0,8){\rouge{\line(1,0){1}} }
\put(1,8){\rouge{\line(0,-1){1}}\put(.25,-.8){$2$}}
\put(1,7){\rouge{\line(1,0){1}} }
\put(2,7){\rouge{\line(0,-1){3}}\put(.25,-.8){$3$}\put(.25,-1.8){$5$}\put(.25,-2.8){$7$}}
\put(2,4){\rouge{\line(1,0){1}}}
\put(3,4){\rouge{\line(0,-1){1}}\put(.25,-.8){$1$}}
\put(3,3){\rouge{\line(1,0){2}}}
\put(5,3){\rouge{\line(0,-1){1}}\put(.25,-.8){$4$}}
\put(5,2){\rouge{\line(1,0){1}}}
\put(6,2){\rouge{\line(0,-1){2}}\put(.25,-.8){$6$}\put(.25,-1.8){$8$}}
\put(6,0){\rouge{\line(1,0){1}}}}
\put(0,8){
  \linethickness{.5mm}
\put(0,8){\rouge{\line(1,0){1}} }
\put(1,8){\rouge{\line(0,-1){1}}\put(.25,-.8){$2$}}
\put(1,7){\rouge{\line(1,0){1}} }
\put(2,7){\rouge{\line(0,-1){3}}\put(.25,-.8){$3$}\put(.25,-1.8){$5$}\put(.25,-2.8){$7$}}
\put(2,4){\rouge{\line(1,0){1}}}
\put(3,4){\rouge{\line(0,-1){1}}\put(.25,-.8){$1$}}
\put(3,3){\rouge{\line(1,0){2}}}
\put(5,3){\rouge{\line(0,-1){1}}\put(.25,-.8){$4$}}
\put(5,2){\rouge{\line(1,0){1}}}
\put(6,2){\rouge{\line(0,-1){2}}\put(.25,-.8){$6$}\put(.25,-1.8){$8$}}
\put(6,0){\rouge{\line(1,0){1}}}
}
\end{picture}\end{center}
\caption{A parking function ``periodically'' extended.}
\label{fig_rect}
\end{figure}
 \item[(d)] For a $(m,n)$-Dyck path $\mu$, let $\rho(\mu):=(c_1,c_2,\ldots,c_k)$ be the sequence of column heights in the skew-shape $(\mu+1^n)/\mu$. Show that the number of $(m,n)$-parking function is given by the combinatorial formula 
      \begin{equation}
          \mathcal{P}_{mn}=\sum_{\mu\subseteq \delta_{mn}} \binom{n}{\rho(\mu)},
      \end{equation}
  with the usual multinomial notation
        $$\binom{n}{\rho(\mu)}=\binom{n}{c_1,c_2,\ldots,c_k}=\frac{n!}{c_1!\,c_2!\cdots c_k!}.$$
\item[(e)] Extend this argument to semi-standard fillings of the shape $(\mu+1^n)/\mu$, to deduce that
     \begin{equation}
         \sum_{\mu\subseteq \delta_{a,b}} s_{(\mu+1^n)/\mu}(\bm{x}) = \frac{1}{a}e_b[a\,\bm{x}]
     \end{equation}
     prove that all expressions in~\pref{formules_ab} are equivalent.
\end{enumerate}

\end{exer}

\begin{exer}\rm
\begin{enumerate}
\item[(a)] Prove proposition~\ref{prop_constant_simple}. You may start your reflexion as follows. Associate to cells $(i,j)$ the weight $q\,z_{i}/z_{i-1}$, and globally weight a partition $\mu$ with the product of the weights of it cells. Here, one sets $z_{-1}:=1$. Find how one may get a constant term out of the expression, and link this to partitions that lie inside $\delta_{mn}$.
\item[(b)]  Formulate a similar statement for partitions lying inside any fixed partition.
\item[(c)] Assuming~\pref{formule_emn_dyck}, prove the constant term formula~(see \cite{negut})
 \begin{equation}
      \bm{e}_{mn}(q,1;\bm{x}):=\left(\frac{E[(z_0+\ldots+z_m)\,\bm{x}\,]}{\bm{z}_{mn}}  \prod_{i=1}^{m-1} \frac{1}{1-q\,z_i/z_{i+1}}\right)\Big|_{\bm{z}^0}
  \end{equation}
\end{enumerate}
\end{exer}

\begin{exer}\rm 
\begin{enumerate}
\item[(a)]
Show that
\begin{equation}
        \pi_n = \frac{qt}{qt-1}\, e_n\left[\frac{\bm{x}(qt-1)}{qt}\right].
   \end{equation}
\item[(b)] Use this to find that the dual basis of $\pi_\mu$ is given be the formula~\footnote{This is useful for the calculation of $\pi$-basis expansions of symmetric functions. }
   \begin{equation}
       \left({(qt-1)}/{qt}\right)^{\ell(\mu)}f_\mu[\bm{x}\,qt/(qt-1)].
   \end{equation}
  \item[(c)]  Show that
            \begin{displaymath}
              e_n(\bm{x})=\sum_{\mu\vdash n} m_\mu\!\Big[\frac{qt}{qt-1}\Big] e_\mu\!\Big[\frac{\bm{x}(qt-1)}{qt}\Big],
        \end{displaymath} 
        from which you may deduce the $\pi$-basis expansion of $e_n$.
      \item[(d)]  Find a plethystic formula for the $\pi$-expansion of $h_n$.
  \end{enumerate}
\end{exer}

\begin{exer}\rm 
  Show that, at $t=1$, formula~\pref{formule_compositional} specializes to formula~\pref{formule_emn_dyck}.
\end{exer}

\begin{exer}\rm
\begin{enumerate}
\item[(a)]
Show that the solutions of \pref{choix_matrice} may be indexed by the integer points that are closest (below) to the ``diagonal'' (from $(0,0)$ to $(m,n)$) in the $m\times n$ rectangle. See figure~\ref{fig_split}.
\begin{figure}[ht]
\begin{center}
\dessin{width=80pt}{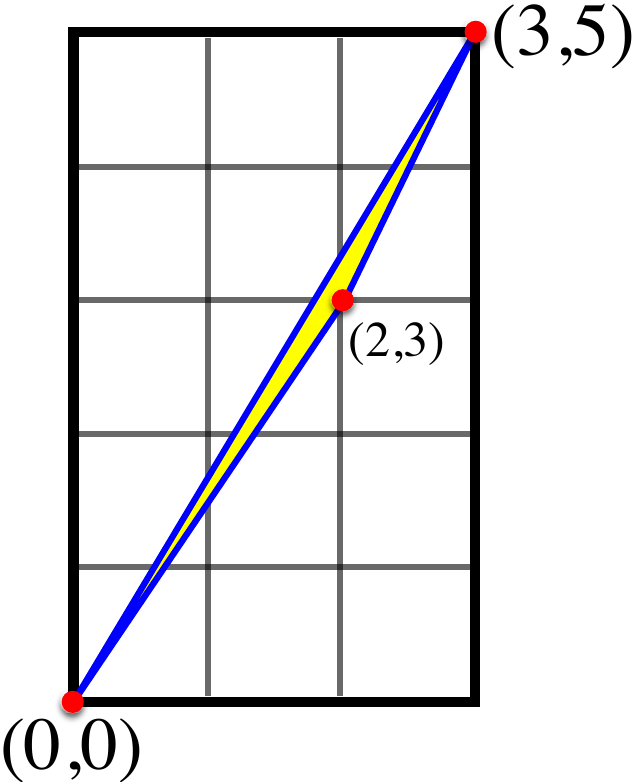}
\end{center}
\vskip-15pt
\caption{A splitting of $(3,5)$}
\label{fig_split}
\end{figure}
\item[(b)] Show that, when $(m,n)$ splits as $(r,s)+(u,v)$, then both the pairs $(r,s)$ and $(u,v)$ are co-prime.
\end{enumerate}
\end{exer}

\begin{exer}\rm See in Sage the function ``GeneralizedTamariLattice'' for the $(m,n)$-Tamari poset when $m>n$ in the coprime case.
\end{exer}


\end{chapter}

\clearpage
\addcontentsline{toc}{chapter}{\bleu{Bibliography}}

\end{document}